\documentclass[opre]{informs3aa} 

\OneAndAHalfSpacedXI 

\theoremstyle{TH}

\usepackage{mathrsfs}

\theoremstyle{TH}
\newtheorem{approximation}{Approximation}

\usepackage[colorlinks=true,breaklinks=true,bookmarks=true,urlcolor=blue,
     citecolor=blue,linkcolor=blue,bookmarksopen=false,draft=false]{hyperref}
     
\usepackage{multirow}
\usepackage[graphicx]{realboxes}
\usepackage{graphicx}
\usepackage{adjustbox}
\usepackage{bm}
\usepackage{dsfont}
\usepackage{comment}
\usepackage{longtable}
\newcommand{\qedwhite}{\hfill \ensuremath{\Box}}
\usepackage{algorithm}
\usepackage{algpseudocode}
\usepackage{wrapfig}

\usepackage{endnotes}
\let\footnote=\endnote
\let\enotesize=\normalsize
\def\notesname{Endnotes}%
\def\makeenmark{\hbox to1.225em{\theenmark.\enskip\hss}}
\def\enoteformat{\rightskip0pt\leftskip0pt\parindent=1.225em
  \leavevmode\llap{\makeenmark}}

\usepackage{natbib}
 \bibpunct[, ]{(}{)}{,}{a}{}{,}%
 \def\bibfont{\small}%
\TheoremsNumberedThrough     
\ECRepeatTheorems

\EquationsNumberedThrough    

\renewcommand{\theARTICLETOP}{}
\begin{document}


\RUNAUTHOR{Li and Mehrotra}

\RUNTITLE{Optimize Resource Allocation in Parallel Any-Scale Queues}

\TITLE{
Optimizing Equitable Resource Allocation in Parallel Any-Scale  Queues with Service Abandonment and its Application to Liver Transplant  
}


\ARTICLEAUTHORS{%
\AUTHOR{Shukai Li}
\AFF{Department of Industrial Engineering and Management Sciences, \\Northwestern University, Evanston, IL 60208, \EMAIL{shukaili2024@u.northwestern.edu}} 
\AUTHOR{Sanjay Mehrotra}
\AFF{Department of Industrial Engineering and Management Sciences, \\Northwestern University, Evanston, IL 60208, \EMAIL{mehrotra@northwestern.edu}} 
} 

\ABSTRACT{%
We study the problem of equitably and efficiently allocating an arriving resource to {\em multiple queues with customer abandonment}. 
The problem is motivated by the cadaveric liver allocation system of the United States, which includes a large number of small-scale (in terms of yearly arrival intensities) patient waitlists  with the possibility of patients abandoning (due to death) until the required service  is completed (matched donor liver arrives). 
We model each waitlist as a GI/MI/1+GI$_S$ queue, in which a virtual server receives a donor liver for the patient at the top of the waitlist, and  patients may abandon while waiting or during service. 
To evaluate the performance of each queue, we develop a finite approximation technique as an alternative to fluid or diffusion approximations, which are inaccurate unless the queue's arrival intensity is large. 
This finite approximation for hundreds of queues is used within an optimization model to optimally allocate donor livers to each waitlist. 
A piecewise linear approximation of the optimization model is shown to provide the desired accuracy.
Computational results show that solutions obtained in this way provide greater flexibility, and improve system performance when compared to solutions from the fluid models. 
Importantly, we find that appropriately increasing the proportion of livers allocated to waitlists with small scales or high mortality risks improves the allocation equity. 
This suggests a proportionately greater allocation of organs to smaller transplant centers and/or those with more vulnerable populations in an allocation policy.
While our motivation is from liver allocation, the solution approach developed in this paper is applicable in other operational contexts with similar modeling frameworks.
}



\maketitle
 
%


\section{Introduction}\label{sec:intro}
This paper is motivated by the need to adequately model patient waitlists and abandonment (death) in an organ transplant system, with the goal of improving organ allocation. 
Specifically, we are interested in studying if the scale of transplant centers and differences in patient mortality risks have an implication in designing an allocation policy.
Transplantation is often the only therapy for patients when organs such as the liver, kidney, heart, etc. fail to function adequately. 
Patients are waitlisted due to the shortage of donated organs, and often die on the waitlist. 

As a case study, we will focus on \emph{cadaveric liver} allocation in the United States.  
In 2021, 13,439 patients were added to the waitlists, while only 9,541 cadaveric livers were donated that year, highlighting the need to achieve equity in liver allocation while maintaining efficiency \citep{OPTN2021}.
Patients register at 140 transplant centers across the country. 
Patients have specific biological characteristics ($e.g.$, blood types), which are matched when a cadaveric liver becomes available. 
Consequently, based on patient characteristics, patients waiting at each transplant center are further divided into multiple waitlists. 
In the current system, the donor livers are procured from hospitals across the nation and subsequently given to a matched patient from a waitlist \citep{OPTN2023}. 

\subsection{Inequity in Liver Transplant}
Waitlist mortality and waiting time are key performance measures for an organ allocation system. 
The United States organ allocation system has suffered from geographical, gender, blood group, and race inequities \citep{jindal2005kidney,majhail2012racial,davis2014changes,davis2014extent,williams2015first}.
For example, Table \ref{table:regional disparity} shows geographical inequities in liver allocation.
In recent years, several operationalizable approaches have been proposed to reduce geographical inequity, including redistricting \citep{kong2010maximizing,gentry2013addressing,gentry2015gerrymandering}, expanding neighborhoods based on donor service areas \citep{kilambi2017improving,mehrotra2018concentric}, and continuous distribution based on a prioritization formula that uses a weighted sum of several attributes \citep{bertsimas2013fairness,mankowski2022designing}.  Continuous distribution is currently under further development at Organ Procurement and Transplantation Network (OPTN) \citep{OPTN2023b}.  
Central to an operationalizable policy is the question of supply-demand matching in an equitable way across patient characteristics.
Such an allocation policy should provide patients with equitable access to donor organs, regardless of their demographics, biological characteristics (blood group types), and waitlist survival risks. Transplant programs vary in the number of transplants performed each year, with lower-size programs generally in rural less densely populated areas.  Also, patients may have different mortality risks while waiting for transplant because of socio-economic reasons and the underlying decease resulting in liver failure (e.g., hepatocellular carcinoma, hepatitis, or fatty liver disease). 
This research is motivated by understanding if the scale of the transplant centers (measured by patient arrival intensity) and patient mortality risks have an implication in designing an allocation policy. 




\begin{table}[!htb]
\tiny
\centering
\caption{Livers' demand-supply imbalances in the United States in 2021. The transplant centers are grouped into 11 OPTN-defined regions. 
New donors, transplants performed, and deaths are expressed as a percentage of annual patient waitlist additions. Source:  \citealt{OPTN2021}. }
\begin{tabular}{lcccclllllcccc}                                                                                                                                                             \cline{1-5} \cline{10-14} 
            & \begin{tabular}[c]{@{}c@{}}Waitlist \\ additions\end{tabular} & \begin{tabular}[c]{@{}c@{}}New\\ donors\end{tabular} & \begin{tabular}[c]{@{}c@{}}Performed\\ transplants\end{tabular} & \begin{tabular}[c]{@{}c@{}}Annual\\ deaths\end{tabular} &  &  &  &  &             & \begin{tabular}[c]{@{}c@{}}Waitlist \\ additions\end{tabular} & \begin{tabular}[c]{@{}c@{}}New\\ donors\end{tabular} & \begin{tabular}[c]{@{}c@{}}Performed\\ transplants\end{tabular} & \begin{tabular}[c]{@{}c@{}}Annual\\ deaths\end{tabular} \\ \cline{1-5} \cline{10-14} 
Region  05  & 2,313                                                         & 65\%                                                 & 61\%                                                            & 10\%                                                     &  &  &  &  & Region  07  & 1,083                                                         & 59\%                                                 & 65\%                                                            & 5\%                                                     \\
Region  03  & 1,892                                                         & 88\%                                                 & 75\%                                                            & 7\%                                                    &  &  &  &  & Region  09  & 833                                                           & 45\%                                                 & 69\%                                                            & 8\%                                                     \\
Region  04  & 1,609                                                         & 57\%                                                 & 56\%                                                            & 8\%                                                     &  &  &  &  & Region  08  & 830                                                           & 85\%                                                 & 65\%                                                            & 5\%                                                     \\
Region  02  & 1,500                                                         & 75\%                                                 & 56\%                                                            & 9\%                                                    &  &  &  &  & Region  01  & 665                                                           & 39\%                                                 & 47\%                                                            & 18\%                                                    \\
Region  11  & 1,263                                                         & 82\%                                                 & 66\%                                                            & 10\%                                                     &  &  &  &  & Region  06  & 374                                                           & 109\%                                                & 67\%                                                            & 7\%                                                     \\ 
 \cline{10-14} 
Region  10  & 1,159                                                         & 81\%                                                 & 76\%                                                            & 6\%                                                    &  &  &  &  & Total & 13,439                                                        & 71\%                                                 & 64\%                                                            & 8\%                                                     \\
\cline{1-5} \cline{10-14} 
\end{tabular}
\label{table:regional disparity} 
\end{table}

\vspace*{0.1in}
\subsection{Modeling Framework, Basic Data Analysis, and Assumptions}
\begin{wrapfigure}{r}{0.5\textwidth}   
\vspace{-33pt}
\begin{center}
\includegraphics[width=.5\textwidth, height =180pt]{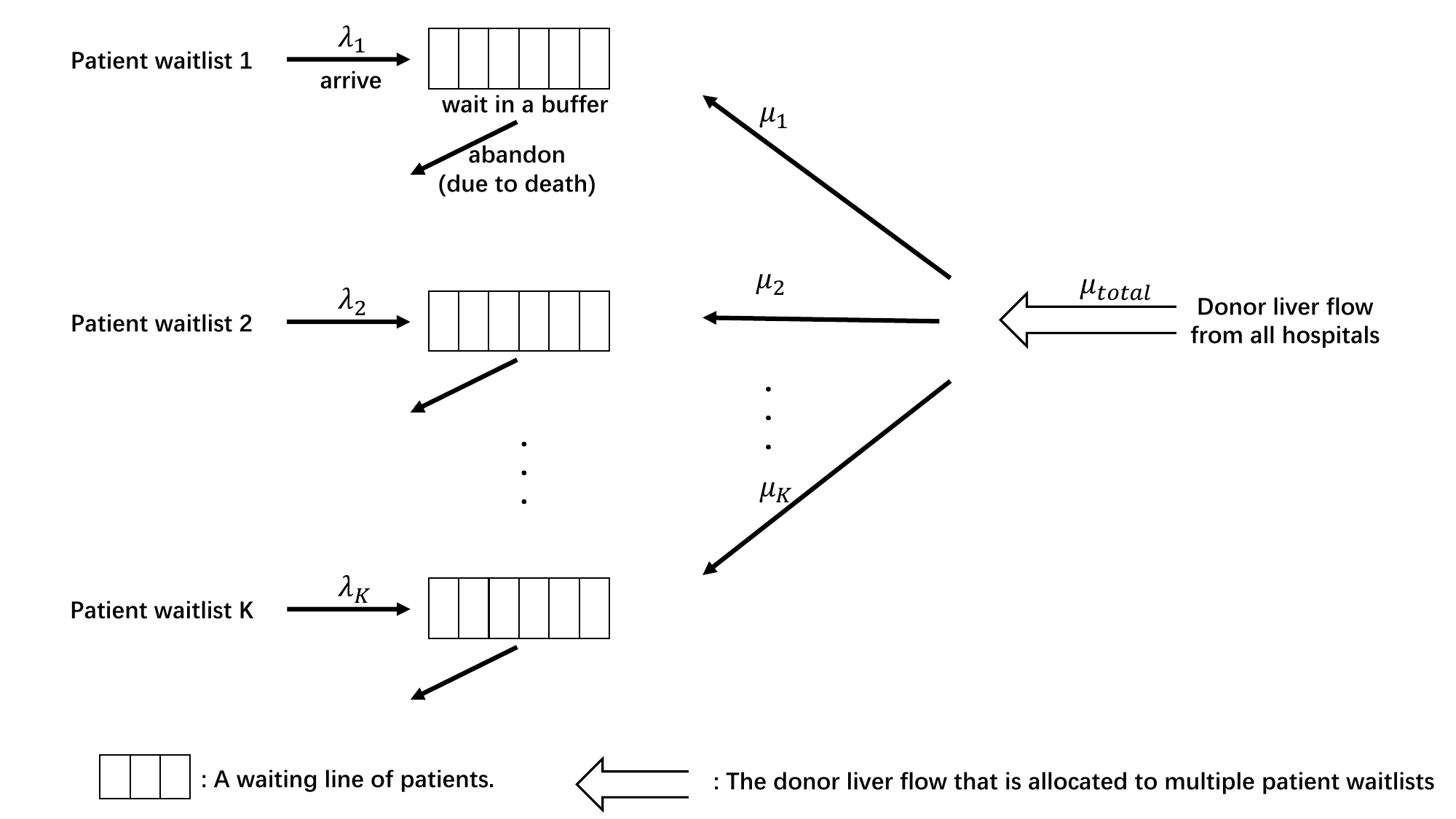}
\end{center}
{Figure 1.2 A schematic view of liver transplant queueing system.} 
\label{fig:1 to K}
\end{wrapfigure}

For patients of certain characteristics, Figure \ref{fig:1 to K} illustrates the allocation process  as modeled in this paper. Each patient waitlist is modeled as a queue. In this queue, the \emph{inter-arrival time} denotes the duration between two consecutive patient arrivals. \emph{The service time} is defined as \emph{the duration starting from when a patient reaches the top of the waitlist until a donor liver becomes available for that waitlist}. 

Patients depart a queue either because of receiving a donor liver, or death while waiting for a liver. 
A patient's death while in the queue is treated as an \emph{abandonment}. 
It is important to note that this abandonment may happen while the patient is at the top of the waitlist, i.e., it {\em may happen during service time}. 
The maximum time a patient can survive without transplant is called the \emph{patience time} of a patient. 
Donor livers are treated as coming from a virtual server. Due to limited shelf-life no inventories are kept. 

\begin{figure}[!htb]  
\centering 
{\includegraphics[width=1\textwidth, height =80pt]{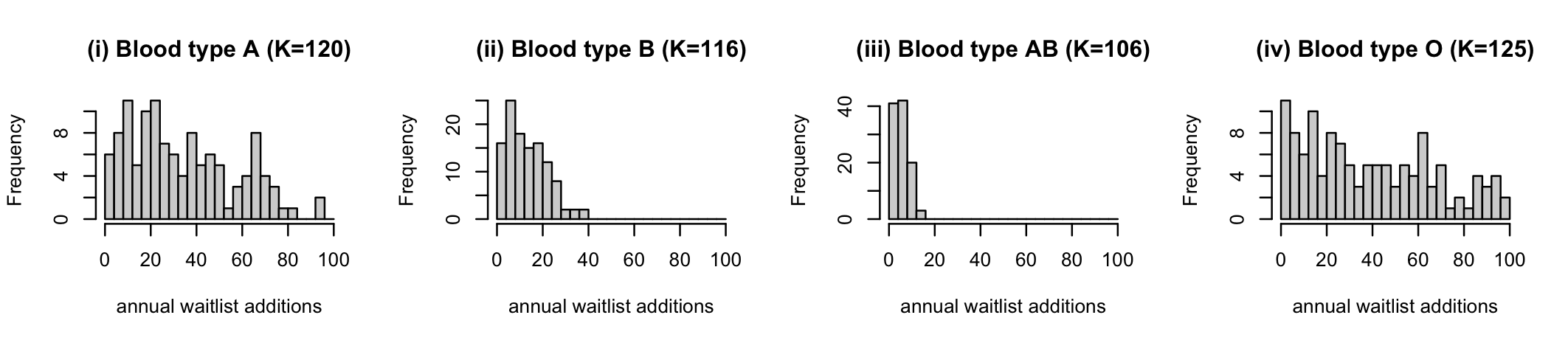}}
\vspace{-30pt}
\caption{Histogram of patient waitlists' arrival intensities by blood types. 
($K=$ number of waitlists).
}
\label{fig:addition dist}
\end{figure}

We analyzed the liver donation and patient arrival data (see Appendix \ref{sec:hypothesistest}). First, we found that more than half of the  waitlists are small-scale queues with patient arrival intensity of less than 12, while the average patience time across the entire patient population is $0.7$ years (see Figure \ref{fig:addition dist}). 
We also found that a mixture of exponential  distributions provides a significantly improved fit compared to an exponential distribution  when modeling the patient inter-arrival times (see Appendix \ref{sec:non-exp of IAT}). In the analysis we allow a general inter-arrival time distribution to specify the patient arrival process. 
We also found that the Poisson process fit is adequate when observing the liver donation process at the largest $50$ donor hospitals (see Appendix \ref{sec:exp of ST}). For the purpose of this study we assume that the donor livers are allocated to queues with a certain probability,  thus allowing for Poisson flow of livers  to patients in a queue.  

\begin{figure}[!htb]  
\centering 
{\includegraphics[width=1\textwidth, height = 75pt]{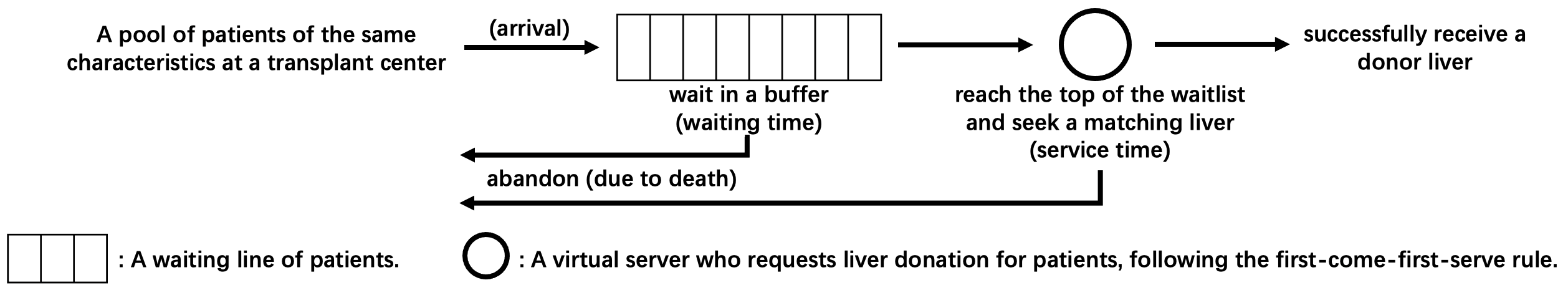}}
\caption{A queueing representation of a  waitlist for a pool of patients with the same characteristics.
For each patient, the waiting time represents the duration it takes to reach the top of the waitlist from the moment of arrival, while the service time represents the duration it takes to receive a donor liver from the moment of reaching the top of the waitlist.}
\label{fig:single queue}
\end{figure}

Although  patients may move randomly on the waitlist due to different disease dynamics, we assume a first-come-first-serve (FCFS) rule in our technical development. However, we let the patience time follow a general distribution with bounded support. We note that a key feature to be modeled is abandonment, and by allowing the patience time to follow a general distribution the random movement of patients on the waitlist can be incorporated in its estimation. 
This allows for the possibility of having sub-classes of patients with different disease progression dynamics within a queue in modeling the abandonment process.  
Our computations use a mixture of three exponential distributions based on three patient sub-classes within a waitlist identified using their disease severity score. 
We also assume that patients do not reject offered livers. 
This is reasonable as most viable donor livers get transplanted due to liver shortage. Overall, we describe each waitlist as a \emph{GI/MI/1+GI$_S$ queue} in which a virtual server receives a donor liver for the patient (customer) at the top of the waitlist, and patients may abandon while waiting or during service (denoted by the subscript ``\emph{$_S$}") as shown in Figure~\ref{fig:single queue}.  Technical assumptions on the inter-arrival time distributions and system performance measures are given in Assumptions~\ref{A:input distributions 1}--\ref{A:input distributions 3}. 

\vspace{.1in}
\noindent \textbf{Model justification.}\quad
Although we are motivated by liver allocation and use it as a case study, we draw some parallels from the literature on kidney allocation.
For instance, similar to \cite{ata2017organjet}, we decompose the organ allocation problem based on immunological types and assume that patients of different immunological types do not share organs. 
We model the organ flow into a waitlist as a virtual server whose function is to receive the donor liver as soon as it becomes available, accounting for the fact that cadaveric organs have short shelf-life and even the patient at the top of the waitlist has to wait for a donation. 
In literature, a transplant waitlist is often modeled as that an (aggregated) Poisson organ flow with rate $\mu$ is allocated to a Poisson patient flow with rate $\lambda$ and forces patients' departure from the waitlist; see, $e.g.$, \cite{zenios1999modeling,zenios2000dynamic,su2006recipient}.
Such a model is equivalent to a MI/MI/1+GI$_S$ queue 
where inter-arrival times and service times follow exponential distributions Exp($\lambda$) and Exp($\mu$), respectively.

\vspace{.1in}
\noindent\textbf{Optimization modeling framework for organ allocation.}\quad
Similar to \cite{zenios1999modeling}, we view the service rate at a patient waitlist as the organ allocation rate, $i.e.$, the number of livers that are allocated to a waitlist during a year. 
The equitable liver allocation problem is formulated as a capacity-sizing problem that optimizes the service rates at multiple (several hundred) queues within a general flexible optimization modeling framework that allows us to incorporate different objectives and constraints based on steady-state performance measures (see Subsubsection~\ref{sec:performance measures}).

\subsection{Contributions}
This paper studies the implications of patient arrival intensities (which translates into the transplant volume at a center) and patient mortality risks in designing an equitable policy framework. To the best of our knowledge, this has not been studied before. 
We develop a new method for solving the resource allocation models where an arriving resource needs to be distributed to multiple parallel queues of various scales with general inter-arrival time and customer patience time distributions. 
The following specific contributions are made. 

\vspace{.in}
\noindent {\bf Evaluating performance of queues with service abandonment.} \quad
We characterize the dynamics of a GI/MI/1+GI$_S$ queue with service abandonment as a Markov chain using the notion of a patient's {offered waiting time} (see Section~\ref{sec:MarkovChainDescription}). 
We show, for the first time, that steady-state performance measures that can be defined using the offered waiting time can be computed to any desired precision with a computable {\em deterministic bound} under mild  conditions (see Assumptions \ref{A:input distributions 1}--\ref{A:input distributions 2}).  
Examples of computable performance measures include average waiting time, average queue length, abandonment rate, the fraction of patients waiting beyond a threshold, 
etc. 
To achieve this, we build upon a \emph{finite approximation} approach introduced in \cite{li2023new}. 
Specifically, the transition kernel arising in the current context is analyzed afresh.

\vspace{.1in}
\noindent{\bf Capacity sizing to parallel queues of various scales.}
\quad
We apply finite approximation for each queue and establish a convergence theory for the general optimization model. 
Specifically, we show that the queue performance measures are continuous with respect to the service rate parameter under mild conditions (see Assumptions \ref{A:input distributions 1}--\ref{A:input distributions 2}). This enables us to use a piecewise linear approximation of the nonlinear constraints in the optimization model.
We subsequently show that the capacity sizing problem can be solved as a piecewise linear program to any desired accuracy. 
To the best of our knowledge, this is the first such methodological advancement for capacity sizing to queues with abandonment that have {\em any} (not necessarily large) arrival intensities, without resorting to simulation optimization. 
The limitations of using simulation, even for a single simpler queue, were demonstrated in \cite{li2023new}.

    
    

\vspace{.1in}
\noindent{\bf Novel findings on liver allocation based on waitlist features.}
\quad
First, we find that assuming  inter-arrival time and patience time distributions as exponential produces significantly different results  when compared to using the mixture of exponential distributions for inter-arrival and patience  times (see Appendix \ref{sec:error of MM1M}).
Our numerical results show that the solutions from the developed optimization model provide greater flexibility and improve resource utilization and equity compared to the solutions prescribed by a fluid model.
Importantly, we find that allocation equity ($resp.$ efficiency) is improved by appropriately increasing ($resp.$ decreasing) the proportion of livers allocated to waitlists with small scales or high mortality risks. 
This suggests a proportionately greater allocation of livers to smaller transplant centers and/or those with more vulnerable populations. We note that the size of a transplant program is also found relevant in risk aversion resulting from center evaluation \citep{mildebrath2021characterizing}.

\subsection{Literature Review}\label{sec:lit rev}
We identify three research articles allowing for abandonment until service ends. 
\cite{daley1965general} provides conditions for the existence of a stationary distribution of a GI/GI/1+GI$_S$ queue. 
They also provide expressions for the Laplace transform of the stationary distribution of the actual waiting time for the special case of exponential inter-arrival time and negative exponential patience time distributions. 
(Note that the actual waiting time equals the minimum between offered waiting time and patience time.)
The Laplace transform of actual waiting time distribution for deterministic patience time is also provided.  
\cite{movaghar2006queueing} considers state-dependent Markovian arrivals and exponential service time distributions. Abandonment rate and stationary distributions of the offer waiting time and queue length are provided. Finally, \cite{choi2001m} assumes Markovian arrivals with deterministic patience time while considering multiple servers.

Most of the other literature assumes abandonment until service starts. Moreover, literature often resorts to asymptotic analysis ($e.g.$, fluid or diffusion approximations), whose accuracy depends on large arrival intensities and service rates, which are often assumed to be greater than 100 to achieve satisfactory accuracy  \citep{whitt2006fluid,reed2008approximating,bassamboo2010accuracy}. 
 Literature has also explored the capacity sizing problems with customer abandonment  \citep{whitt2004efficiency,zeltyn2005call,ward2008asymptotically,mandelbaum2009staffing,bassamboo2010accuracy,lee2019pricing}. These models do not consider abandonment during service.

The transplant system operations literature that is most closely related to the current work is \cite{zenios2000dynamic,su2006recipient,akan2012broader,ata2021achievable,khademi2021asymptotically,hasankhani2021time,hasankhani2022proportionally}. 
These works are based on fluid approximation, which degenerates the original model into a deterministic fluid model by sending the arrival and total service rates to infinity. As remarked by \cite{zenios2000dynamic}, ``a fluid model is too crude to be used as a reliable performance analysis tool" for transplant system operations. 
The approximate solutions have no optimality guarantees; an exception is \cite{khademi2021asymptotically}, where the analysis is performed as a stochastic queuing control problem with time-varying Markovian arrivals.  We note that the gap in recommendations based on fluid approximation from that based on more accurate modeling is also not well understood.

Organs for transplant may be donated by deceased or living donors, $e.g.$, a patient's family member or an altruistic donation. In the United States, 
94\% of the liver transplants
were from deceased donors during the first half of 2022 \citep{OPTN2021}. Living donors typically make donations to a target recipient. A classical problem is how to maximize the organ exchange efficiency among patient-donor pairs \citep{roth2004kidney,anderson2015finding,ding2015non,agarwal2019market,ashlagi2021kidney}. For a more systematic review of managing living-donor transplant, we refer the readers to \cite{alagoz2004optimal,kaufman2017living}. The current work only considers deceased donors as living donations are often pre-arranged, and in case of liver only a small percent of transplants involve living donors. Nevertheless, patient abandonment due to finding a living donor while waiting on the queue can also be modeled as part of the abandonment distribution.



Other transplant system operations literature focuses on patient decision-making and explores ways to maximize individual or social welfare by considering factors such as whether to accept or decline an offered organ \citep{david1985time,su2004patient,su2005patient,su2006recipient,alagoz2007choosing,alagoz2007determining,batun2018optimal}, robust waiting time estimation \citep{bandi2019robust}, information transparency \citep{sandikcci2008estimating,sandikcci2013alleviating}, administration region design \citep{stahl2005methodological,kong2010maximizing, demirci2012exact}, multiple listing \citep{ata2017organjet}, medication management \citep{boloori2020data}, regulatory-induced risk aversion for transplant programs \citep{mildebrath2021characterizing}, etc.

\subsection{Organization and Basic Notation}
Section \ref{sec:model} provides the Markov chain representing a GI/MI/1+GI$_S$ queue allowing for abandonment until the end of service, discusses computations of performance measures, and provides a general capacity sizing model for multiple queues. 
Section \ref{sec:FAO} provides a two-step approximation for finding the optimal allocation.
We provide a convergence theory for our model in Section \ref{sec:main results}. 
We discuss our algorithm's computational advantages in Section \ref{sec:alg property}.
We validate the accuracy of our method via numerical experiments and show our solution outperforms that prescribed by a fluid model in Section \ref{sec:numerical}.
In Section \ref{sec:conclusion}, we discuss policy implications. 
Appendix \ref{Append:Notation} provides a summary of major notation.  
Appendix \ref{sec:hypothesistest} discusses statistical hypothesis tests for inter-arrival time and service time distributions. 
Appendice \ref{Append:Arrivial distribution}--\ref{append:proof} provide preliminaries and part of the proof for lemmas and theorems developed in this paper.
We define $f(x+):=\lim_{u\rightarrow x+}f(u)$ and $f(x-):= \lim_{u\rightarrow x-}f(u)$. 
$[x]_+:=\max\{x,0\}$. 
For $x\in\mathbb{R}^n$, dim$(x):=n$.  
``\emph{i.i.d.}" stands for ``independent and identically distributed". 
``\emph{w.r.t.}" stands for ``with respect to". 

\section{Model Formulation }\label{sec:model}
The system being studied is comprised of $K$ parallel queues, characterized by customer arrival intensities $\{\lambda_k\}_{k \in \mathscr{K}}$ and general customer inter-arrival time distributions, where $\mathscr{K}:= \{1,2,\dots,K\}$. 
A resource that arrives following a Poisson process needs to be equitably and efficiently allocated to these queues respectively using allocation rates $\{\mu_k\}_{k\in\mathcal{K}}$. 
The allocated resources to each queue also follow a Poisson process.  
Each arriving customer waits in the queue until they receive a resource or their patience time elapses.
Let $\pmb{\mu} := (\mu_1,\mu_2,\dots,\mu_K)$, which is to be determined. 
We assume $\pmb{\mu}\in[\mu_{\min},\mu_{\max}]^K$, where $\mu_{\max}>\mu_{\min}>0$. 

In the context of liver allocation, as shown in Figure \ref{fig:1 to K}, patients (customers) with different immunological, health, socio-economic, and/or geographic status form $K$ waitlists. 
Donor livers (resources) are then allocated to these waitlists to achieve equitable allocation while maintaining system efficiency.  
Each waitlist is described as a GI/MI/1+GI$_S$ queue, as shown in Figure \ref{fig:single queue}, each having a virtual server that follows the FCFS rule and receives donor livers for the patient at the top of the waitlist. 
The service time represents the duration starting from \emph{when the patient reaches the top of the waitlist} until \emph{a donor liver is available for that waitlist}.
In waitlist $k$, service times follow an exponential distribution Exp($\mu_k$) because donor livers allocated to waitlist $k$ follow a Poisson process with rate $\mu_k$, which is memory-less.
Each patient has a random amount of \emph{patience time} due to death and \emph{abandons} the system if their patience time elapses before the service is completed, $i.e.$, before a donor liver is successfully matched. 




\subsection{Markov Chain Representation of 
GI/MI/1+GI$_S$ Queue}\label{sec:model description} \label{sec:MarkovChainDescription}
We define \emph{offered waiting time}  as the duration it takes a patient to \emph{reach the top of the waitlist} 
from the time of their arrival should this patient have infinite patience time. 
Let $t_k^n$, $s_k^n$, $y_k^n$, and $\xi_k^n$ respectively be the inter-arrival time, service time, patience time, and offered waiting time of patient $n\in\mathbb{N}$. 
For all $n\in\mathbb{N}$, we have
\begin{align}
    \xi_k^{n+1} = \begin{cases}
        [\xi_k^n+s_k^{n}-t_k^{n+1}]_+, \hspace{20pt} \text{if\,} y_k^n > \xi_k^n+s_k^{n}, \\
        [\xi_k^n-t_k^{n+1}]_+, \hspace{43pt} \text{if\,} y_k^n < \xi_k^n,\\
        [y_k^n-t_k^{n+1}]_+, \hspace{43pt} \text{otherwise}.
    \end{cases} \label{eqn:original MC}
\end{align}
\eqref{eqn:original MC} describes three possible outcomes for patient $n+1$.  
$(a)$ If the patience time $y_k^n$ of patient $n$ is greater than the offered waiting time $\xi_k^n$ plus the service time $s_k^n$, patient $n$ is able to complete service, $i.e.$, to receive a donor liver.
Patients $\{0,1,\dots,n\}$ are all removed from the system 
after a duration of $\xi_k^n+s_k^n$ units of time since the arrival of patient $n$.  
Thus, the offered waiting time for patient $n+1$ will be $[\xi_k^n+s_k^n-t_k^{n+1}]_+$.
$(b)$ If the patience time $y_k^n$ of patient $n$ is smaller than the offered waiting time $\xi_k^n$, patient $n$ does not receive any service, $i.e.$, patient $n$ dies and never reaches the top of the waitlist.
Patients $\{0,1,\dots,n\}$ are all removed from the system after a duration of $\xi_k^n$ units of time since the arrival of patient $n$.
Thus, the offered waiting time for patient $n+1$ will be $[\xi_k^n-t_k^{n+1}]_+$.
$(c)$ If the patience time $y_k^n$ of patient $n$ is between $\xi_k^n$ and $\xi_k^n+s_k^n$, patient $n$ starts service but is unable to complete, $i.e.$, patient $n$ dies at the top of the waitlist  without receiving a liver. 
Patients $\{0,1,\dots,n\}$ are all removed from the system after a duration of $y_k^n$ units of time since the arrival of patient $n$.
Thus, the offered waiting time for patient $n+1$ is $[y_k^n-t_k^{n+1}]_+$.

We assume that, in waitlist $k\in\mathscr{K}$, the inter-arrival times $\{t^n_k\}_{n\in\mathbb{N}}$, service times $\{s^n_k\}_{n\in\mathbb{N}}$, and patience times $\{y^n_k\}_{n\in\mathbb{N}}$ form three \emph{i.i.d.} sequences that are independent of each other.
Then the offered waiting time sequence $\{\xi_k^n\}_{n\in\mathbb{N}}$ defined in \eqref{eqn:original MC} forms a (discrete-time) Markov chain supported on $\mathbb{R}_+$. 
This Markov chain is denoted by $\pmb{\Psi}_k(\mu)$ when service rate $\mu_k=\mu$. 
We let $A_k(x)$ be the distribution of inter-arrival time, $B(x;\mu):=1-\exp(-\mu x)$ be the exponential distribution of service time when service rate $\mu_k=\mu$, and $G_k(x)$ be the distribution of patience time. 
Particularly, the arrival intensity  $\lambda_k=1/\int_{\mathbb{R}_+}xdA_k(x)$.

\subsubsection{Computable Steady-State Performance Measures}\label{sec:performance measures}
We consider the  system's performance at steady state.
For waitlist $k\in\mathscr{K}$, we use $\pi_k(x;\mu)$ to denote the stationary distribution function of Markov chain $\pmb{\Psi}_k(\mu)$ that describes offered waiting time $\xi_k^n$ when service rate $\mu_k=\mu$. 
We express waitlist $k$'s steady-state performance measures as $\mathbb{E}_{\pi_k(\xi;\mu)} g(\xi,\mu)$.
Here
$\xi\in\mathbb{R}$ denotes a generic offered waiting time, $g(\xi,\mu)$ is a generic \emph{performance function} depending on offered waiting time $\xi$ and service rate $\mu$, and $\mathbb{E}_{\pi_k(\xi;\mu)} g(\xi,\mu):=\int_\mathbb{R} g(\xi,\mu) d\pi_k(\xi;\mu)$ is integral $w.r.t.$ $\xi$. 

For a better understanding of how to use the offered waiting time distribution to construct steady-state performance measures, we provide some illustrative examples.
Following \cite{movaghar2006queueing}, we define \emph{offered sojourn time} as the duration it takes a patient to \emph{receive a donor liver} from the time of their arrival should this patient have infinite patience time.
For waitlist $k$ with service rate $\mu_k=\mu$, patients' average offered sojourn time can be written as $\int_{\mathbb{R}_+}\int_{\mathbb{R}_+}(\xi+s) dB(s;\mu)d\pi_k(\xi;\mu)=\mathbb{E}_{\pi_k(\xi;\mu)}\big(\xi+\frac{1}{\mu}\big)$.  
The abandonment probability, $i.e.,$ the fraction of patients who die without receiving a donor liver, can be written as $\int_{\mathbb{R}_+}\int_{\mathbb{R}_+}\int_{\mathbb{R}_+}\mathds{1}\{y\leqslant \xi+s\} dG_k(y)dB(s;\mu)d\pi_k(\xi;\mu)=\mathbb{E}_{\pi_k(\xi;\mu)}\big[\int_{\mathbb{R}_+}\mu\exp(-\mu s)G_k(\xi+s) ds\big]$. 
The fraction of patients having to wait for more than $\bar{\xi}$ time to receive a donor liver can be written as $\int_{\mathbb{R}_+}\int_{\mathbb{R}_+}\mathds{1}\{\xi+s > \bar\xi\} dB(s;\mu)d\pi_k(\xi;\mu)=\mathbb{E}_{\pi_k(\xi;\mu)}\big[\int_{\mathbb{R}_+}\mu\exp(-\mu s)\mathds{1}\{\xi+s > \bar\xi\} ds\big]$.
The average waitlist length can be written as $\int_{\mathbb{R}_+}\int_{\mathbb{R}_+}\int_{\mathbb{R}_+}\lambda_k\min\{y,\xi+s\} dG_k(y)dB(s;\mu)d\pi_k(\xi;\mu)=\mathbb{E}_{\pi_k(\xi;\mu)}\big[\int_{\mathbb{R}_+}\int_{\mathbb{R}_+}\lambda_k\mu\exp(-\mu s)\min\{y,\xi+s\} dG_k(y)ds\big]$. 
Lastly, we consider QALYs gained as a possible evaluation metric. 
We let $h^{\text{QALY}}_k(x)$ be the average QALYs gained by a patient who entered waitlist $k$ and received a donor liver after a duration of $x$ units of time since arrival. The
function $h^{\text{QALY}}_k$ can be estimated using data; see, $e.g.$, \cite{axelrod2018economic}.
Then $g^{\text{QALY}}_{k}(\xi,\mu):= \int_{\mathbb{R}_+} \bar{G}_k(\xi+s) h^{\text{QALY}}_k(\xi+s) dB(s;\mu)=\int_{\mathbb{R}_+} \mu\exp(-\mu s)\bar{G}_k(\xi+s) h^{\text{QALY}}_k(\xi+s) ds$ represents the average QALYs gained by a patient with offered waiting time $\xi$.
Here $\bar{G}_k(x):= 1-G_k(x)$. 
Thus, $\mathbb{E}_{\pi_k(\xi;\mu)} g^{\text{QALY}}_k(\xi,\mu)$ represents the average QALYs gained per patient in waitlist $k$ when service rate $\mu_k=\mu$. 

The performance evaluation can be done using the implicit stationary distribution $\pi_k(x;\mu)$ of Markov chain $\pmb{\Psi}_k(\mu)$. 
Stationary distribution function $\pi_k(x;\mu)$ is a solution to the following infinite-dimensional balance equations
\begin{align}
    \pi_k(x;\mu)= \int_{\mathbb{R}_+} \tau_k(x,u;\mu) d\pi_k(u;\mu), \quad x \in \mathbb{R}_+, \label{eqn:balance eq}
\end{align}
where the integral is defined $w.r.t.$ $u$ and the transition kernel $\tau_k(x,u;\mu):=\mathbb{P}[\xi_k^{n+1}\leqslant x | \xi_k^{n}= u]$ returns the probability that the next patient's offered waiting time is less than $x$, given that the current patient's offered waiting time is $u$, under service rate $\mu_k=\mu$. 
By the definition of transition kernel $\tau_k$ and the Markov chain dynamics in \eqref{eqn:original MC}, we have the following: 
\begin{proposition}
The transition kernel $\tau_k(x,u;\mu)$ can be expressed as 
\begin{align}
    \tau_k(x,u;\mu) = \int_{\mathbb{R}_+} \int_{\mathbb{R}_+} \mu\exp(-\mu s)A^*_k(\min\{u+s, \max\{y,u\}\}-x)  dG_k(y)ds,&\nonumber\\ 
    (x,u)\in \mathbb{R}_+^2,\mu\in[\mu_{\min},\mu_{\max}].& \label{eqn:def kernel} 
\end{align}
Here $A^*_k(x):=1-A_k(x-)$.
\end{proposition}
\begin{proof}{Proof:}
    By definition of $\tau_k$ and the equation of $\xi_k^{n+1}$ in \eqref{eqn:original MC}, we have that, for all $x,u\in\mathbb{R}_+$, 
    \begin{align*}
        \tau_k(x,u;\mu)=& \mathbb{P}[\xi_k^{n+1}\leqslant x | \xi_k^{n}= u] = \mathbb{P}\Big\{\big[\min\{\xi_k^n+s_k^n,\max\{\xi_k^n,y_k^n\}\}-t_k^{n+1}\big]_+ \leqslant x \Big| \xi_k^{n}= u\Big\}\\
        =& \mathbb{P}\Big\{\min\{u+s_k^n,\max\{y_k^n,u\}\}-t_k^{n+1} \leqslant x\Big\}=A_k^*(\min\{u+s_k^n,\max\{y_k^n,u\}\}- x)\\ 
        = & \int_{\mathbb{R}_+} \int_{\mathbb{R}_+} \mu\exp(-\mu s)\cdot A^*_k(\min\{u+s, \max\{y,u\}\}-x)  dG_k(y)ds.\hspace{65pt} \qedwhite
    \end{align*}
\end{proof}
Kernel function $\tau_k$ can be evaluated given the service rate $\mu$ and the distributions $A_k$ and $G_k$. 
How to solve balance equations \eqref{eqn:balance eq} and obtain stationary distribution $\pi_k$ to any desired precision is unknown. 

\subsection{Capacity Sizing with Flexible Objectives and Constraints}\label{sec:opt problem formulation}
We consider a general class of optimization models over parallel queues formulated as follows:
\begin{align}
    \min_{\pmb{\mu},\pmb{w}} \quad & \pmb{p}^\text{T} \pmb{w} \label{opt:objective}\tag{P} \\
    \text{s.t.}\quad & \mathbb{E}_{\pi_k(\xi;\mu_k)}g_l(\xi,\mu_k) \leqslant w_{k,l},\quad (k,l)\in\mathscr{K}\times\mathscr{L}, \label{opt:constraint relax} \tag{\ref{opt:objective}-a}\\
    & \pmb{M} \pmb{w} \leqslant \pmb{d}, \label{opt:constraint} \tag{\ref{opt:objective}-b}\\
    & \pmb{\mu}\in [\mu_{\min},\mu_{\max}]^{K},\, \pmb{w}\in\mathbb{R}^{K\times L}.
    \label{opt:constraint lu bound} \tag{\ref{opt:objective}-c}
\end{align}
Here a set $\mathscr{L}$ of performance measures are considered and $L:=|\mathscr{L}|$.
$\mathbb{E}_{\pi_k(\xi;\mu_k)}g_l(\xi,\mu_k)$ represents waitlist $k$'s steady-state performance under measure $l$ and service rate $\mu_k$.
Variable vector $\pmb{w} = (w_{k,l})_{k\in \mathscr{K},l\in\mathscr{L}} $ collects all performance measures that we consider to optimize the system.
Constant vector $\pmb{p} \in \mathbb{R}^{\text{dim}(\pmb{w})} $ collects cost coefficients associated with performance measures. 
Constant vector $\pmb{d}$ and matrix $\pmb{M}\in \mathbb{R}^{\text{dim}(\pmb{d})\times\text{dim}(\pmb{w}) }$ construct additional constraints.
When $\pmb{p}$ and $\pmb{M}$ are non-negative, objective \eqref{opt:objective} and constraint \eqref{opt:constraint} ensure that constraint \eqref{opt:constraint relax} is tightly satisfied by any optimal solution, $i.e.$, $w_{k,l}$ returns the exact value rather than merely an upper bound of $\mathbb{E}_{\pi_k(\xi;\mu_k)}g_l(\xi,\mu_k)$. Assumptions on $g_l(\xi,\mu_k)$ are stated in Assumption~\ref{A:input distributions 2}. 

\vspace{.1in}
\noindent\textbf{Equitable liver allocation problem formulation.}\quad
We examine two liver allocation models that serve as special cases of \eqref{opt:objective}. 
In model \eqref{model:ost}, we maximize the equity of offered sojourn time across waitlists.
The offered sojourn time is an important performance measure as it represents the duration it takes a patient to receive a donor liver from the time of their arrival.
Specifically, we solve the following optimization problem: 
\begin{align}
    \min_{\pmb{\mu},\pmb{w},\bar{w}, \pmb{z}, Z } \quad & Z \label{model:ost} \tag{P-OST}\\
    \text{s.t.}\qquad & \mathbb{E}_{\pi_k(\xi;\mu_k)}\big(\xi+\frac{1}{\mu_k}\big) = w_{k},\quad k\in\mathscr{K}, \label{model:ost-1}  \tag{\ref{model:ost}-a}\\
    & \pmb{1}^\text{T}\pmb{\lambda}\cdot\Bar{w} = \pmb{\lambda}^{\text{T}} \pmb{w}, \quad  - \pmb{z} \leqslant \pmb{w} -\Bar{w} \leqslant \pmb{z},\quad KZ= \pmb{1}^\text{T}\pmb{z}, \label{model:ost-2} \tag{\ref{model:ost}-b}\\
    & \bar{w}\leqslant \varsigma,\quad \pmb{1}^\text{T}\pmb{\mu}\leqslant \mu_{\text{total}}, \label{model:ost-3}\tag{\ref{model:ost}-c}\\
    & \mu_k\in[\vartheta_L \lambda_k, \vartheta_U \lambda_k],k\in\mathcal{K}; \pmb{w}\in\mathbb{R}_+^{K}, \bar{w}\in\mathbb{R}_+, \pmb{z}\in\mathbb{R}_+^{K}, Z\in \mathbb{R}_+. \label{model:ost-4} \tag{\ref{model:ost}-d}
\end{align}
In \eqref{model:ost-1}, $w_k$ represents patients' expected offered sojourn time in waitlist $k$.
In \eqref{model:ost-2}, the efficiency measure $\Bar{w}$ represents the weighted average value of expected offered sojourn time across all waitlists and the equity measure $Z$ represents the mean absolute deviation of expected offered sojourn time across all waitlists.
\eqref{model:ost-3} includes an efficiency constraint and a resource constraint: $\bar{w}\leqslant \varsigma$ represents that the system's offered  sojourn time should be less than $\varsigma$ on average and $\pmb{1}^\text{T}\pmb{\mu}\leqslant \mu_{\text{total}}$ specifies an upper bound for the total liver allocation rates.
\eqref{model:ost-4} requires that, for each waitlist, the liver supply must align with the demand within a specified proportion range, denoted by $[\vartheta_L,\vartheta_U]$.
As the efficiency constraint relaxes ($i.e.$, as $\varsigma$ increases), allocation equity is expected to improve. 

We also use model \eqref{model:ab} to maximize the equity of abandonment probability across waitlists. 
The abandonment probability is an important performance measure as it represents  patient mortality on the waitlist.
The interpretation of model \eqref{model:ab} is the same as that of \eqref{model:ost} except that ``offered sojourn time" is replaced by ``abandonment probability".
\begin{align}
    \min_{\pmb{\mu},\pmb{w},\bar{w}, \pmb{z}, Z } \quad & Z \label{model:ab}\tag{P-Ab} \\
    \text{s.t.}\qquad & \mathbb{E}_{\pi_k(\xi;\mu_k)}\Big[\int_{\mathbb{R}_+}\mu_k\exp(-\mu_k s)G_k(\xi+s) ds\Big] = w_{k},\quad k\in\mathscr{K}, \label{model:ab-1} \tag{\ref{model:ab}-a}\\
    & \pmb{1}^\text{T}\pmb{\lambda}\cdot\Bar{w} = \pmb{\lambda}^{\text{T}} \pmb{w}, \quad  - \pmb{z} \leqslant \pmb{w} -\Bar{w} \leqslant \pmb{z},\quad KZ= \pmb{1}^\text{T}\pmb{z},\label{model:ab-2}\tag{\ref{model:ab}-b}\\
    & \bar{w}\leqslant \varsigma,\quad \pmb{1}^\text{T}\pmb{\mu}\leqslant \mu_{\text{total}}, \label{model:ab-3}\tag{\ref{model:ab}-c}\\
    & \mu_k\in[\vartheta_L \lambda_k, \vartheta_U \lambda_k],k\in\mathcal{K}; \pmb{w}\in\mathbb{R}_+^{K}, \bar{w}\in\mathbb{R}_+, \pmb{z}\in\mathbb{R}_+^{K}, Z\in \mathbb{R}_+. \label{model:ab-4}\tag{\ref{model:ab}-d}
\end{align}


\section{Solution using Finite and Piecewise Linear Approximation} 
\label{sec:FAO}
The capacity sizing problem \eqref{opt:objective} involves evaluating steady-state performance measures of continuous-state Markov chains $\big\{\pmb{\Psi}_k(\mu)\big\}_{k\in\mathscr{K}}$, which are implicit functions of service rates $\pmb{\mu}$. 
Evaluating these measures is a challenging task since the infinite-dimensional balance equations \eqref{eqn:balance eq} do not have an analytical solution for the stationary distribution. 
We now present an approximation approach for finding a solution of \eqref{opt:objective} to any desired precision.
First, we construct a finite-state Markov chain to approximate the original continuous-state Markov chain $\pmb{\Psi}_k(\mu)$, which enables us to obtain approximate performance measures (Subsection \ref{sec:construct MC}).
Then we construct piecewise linear (PWL) constraints to approximate constraints \eqref{opt:constraint relax} (Subsection \ref{sec:construct pwl constraints}). 
Model \eqref{opt:objective} is now solved using a deterministic PWL program (see Algorithm \ref{alg1}).



\subsection{Finite-State Markov Chain Approximation}
\label{sec:construct MC}
We now present the procedure for constructing an approximate finite-state Markov chain for a given instance of waitlist $k\in\mathscr{K}$, with service rate $\mu_k=\mu\in[\mu_{\min},\mu_{\max}]$.
Suppose the patience times $\{y_k^n\}_{n\in\mathbb{N}}$ 
 are bounded by constant $\bar{y}_k$ (see Assumption \ref{A:input distributions 1} of Section \ref{sec:main results}). 
Then the offered waiting times $\{\xi_k^n\}_{n\in\mathbb{N}}$ defined in \eqref{eqn:original MC} are uniformly bounded by  $\bar{y}_k$; that is,  Markov chain $\pmb{\Psi}_k(\mu)$ is supported on $[0,\bar{y}_k]$. 
This is because, for any patient $n$, the patience times of preceding patients are all bounded by $\bar{y}_k$, and patient $n$ is guaranteed to reach the top of the waitlist within a duration of $\bar{y}_k$ units of time since their arrival. 

We utilize a truncation approach to generate a sequence of finite-state Markov chains, denoted by $\big\{\pmb{\Psi}^{(r)}_k(\mu)\big\}_{r\in\mathbb{N}}$, to approximate the orignal continuous-state Markov chain $\pmb{\Psi}_k(\mu)$.
The $r$-th Markov chain $\pmb{\Psi}^{(r)}_k(\mu)$ has $J^{(r)}=2^r+1$ states.
The states and transition matrix are respectively denoted by $\{c^{(r)}_{k_i}\}_{i=1}^{J^{(r)}}$ and $\pmb{Q}_k^{(r)}(\mu)=\big\{q^{(r)}_{k_{i,j}}(\mu)\big\}_{i,j=1,2,\dots,J^{(r)}}$, where 
\begin{align}
    c^{(r)}_{k_i} := \frac{(i-1)\bar{y}_k}{2^r},\qquad q_{k_{i,j}}^{(r)}(\mu):= \tau_k\big(c^{(r)}_{k_j},c^{(r)}_{k_i};\mu\big)-\tau_k\big(c^{(r)}_{k_{j-1}},c^{(r)}_{k_i};\mu\big),\,\qquad  i,j=1,2,\dots,J^{(r)}.
    \label{eqn:def approximate MC}
\end{align}
Particularly, $c^{(r)}_{k_0}:= -\infty$. 

We use $\pi_k^{(r)}(x;\mu)$ to denote the discrete stationary distribution under transition matrix $\pmb{Q}_k^{(r)}(\mu)$ and use $\pmb{v}^{(r)}_k(\mu) = \big({v}^{(r)}_{k_1}(\mu),{v}^{(r)}_{k_2}(\mu),\dots,{v}^{(r)}_{k_{J^{(r)}}}(\mu)\big)^\text{T}$ to denote the stationary probability vector, $i.e.$, $\pmb{v}^{(r)}_k(\mu)$ is the vector of jump sizes of stationary distribution $\pi^{(r)}_k(x;\mu)$. 
By definition, $\pmb{v}^{(r)}_k(\mu)$ can be obtained by solving a \emph{finite} set of balance equations
\begin{align}
    &\Big( \pmb{Q}_k^{(r)}(\mu)-\pmb{I}\Big)^\text{T} \pmb{v}^{(r)}_k(\mu) = \mathbf{0},\qquad \mathbf{1}^\text{T} \pmb{v}^{(r)}_k(\mu) = 1.
    \label{eqn:approx stationary vector}
\end{align}
We will show in Theorem \ref{theorem:MC compute} that 
\eqref{eqn:approx stationary vector} can be solved exactly by matrix inversion or numerically by iterations.
Then distribution $\pi_k^{(r)}$ can be obtained by
\begin{align}
    &\pi_k^{(r)}(x;\mu) =  \sum_{i=1}^{J^{(r)}} \mathds{1}\{x\geqslant {c}^{(r)}_{k_i}\}{v}^{(r)}_{k_i}(\mu).\label{eqn:approx stationary dist}
\end{align}

For a performance function $g(\xi,\mu)$, we approximate the performance measure $\mathbb{E}_{\pi_k(\xi;\mu)}g(\xi,\mu)$ by
\begin{align}
    \mathbb{E}_{\pi_k(\xi;\mu)}g(\xi,\mu)\approx \mathbb{E}_{\pi^{(r)}_k(\xi;\mu)}g(\xi,\mu)=\sum_{i=1}^{J^{(r)}}g(c_{k_i}^{(r)},\mu)v_{k_i}^{(r)}(\mu).\label{eqn:approx performance measure}
\end{align}
We call \eqref{eqn:approx performance measure} a \emph{finite approximation}, and the resulting performance measure a \textit{finite approximate performance measure}. 
Here the original stationary distribution $\pi_k$ of continuous-state Markov chain $\pmb{\Psi}_k(\mu)$ is replaced by an approximate stationary distribution $\pi^{(r)}_k$ of a finite-state Markov chain $\pmb{\Psi}^{(r)}_k(\mu)$.
Lemma \ref{lemma:fa uniform convergence} shows that the accuracy of approximation in \eqref{eqn:approx performance measure} improves with increasing $r$, i.e., by adding more finite approximation states.

\subsection{Piecewise Linear Approximation of Constraints} \label{sec:construct pwl constraints}
We now describe a procedure for constructing a PWL constraint to approximate an instance of constraint $(k,l)\in\mathscr{K}\times\mathscr{L}$ in \eqref{opt:constraint relax}.
We generate a sequence of PWL functions $\big\{\phi_{k,l}^{(r,m)}(\mu)\big\}_{m\in\mathbb{N}}$, where $r\in\mathbb{N}$ specifies the order of finite approximation.
For $m\in\mathbb{N}$, we use  $N^{(m)}:=2^m+1$ equally spaced knots:
\begin{align*}
    \mu_i^{(m)} :=\mu_{\min} + \frac{i}{2^m} (\mu_{\max}-\mu_{\min}),\quad i=0,1,\dots,2^{m}.
\end{align*}
Then we define a PWL function:
\begin{align}
    \small
    \phi^{(r,m)}_{k,l}(\mu) := 
    \begin{cases}
        \mathbb{E}_{\pi^{(r)}_k(\xi;\mu^{(m)}_{i})}g_l(\xi,\mu^{(m)}_{i}), \hspace{122pt} \mu=\mu^{(m)}_{i},i\in\{0,1,\dots,2^m\};\\[-10pt]\\
        \phi^{(r,m)}_{k,l}(\mu^{(m)}_{i-1}) + \frac{\phi^{(r,m)}_{k,l}(\mu^{(m)}_i) -\phi^{(r,m)}_{k,l}(\mu^{(m)}_{i-1}) }{\mu^{(m)}_i- \mu^{(m)}_{i-1}}(\mu- \mu^{(m)}_{i-1} ),\quad\, \mu \in (\mu^{(m)}_{i-1},\mu^{(m)}_i), i\in\{1,2,\dots, 2^m\}.
    \end{cases} \label{eqn:pwl function}
\end{align}
For any $(r,m)\in\mathbb{N}^2$, $\phi^{(r,m)}_{k,l}(\mu)$ as a function of $\mu\in[\mu_{\min},\mu_{\max}]$ can be regarded as a PWL approximator to the finite approximate performance measure function $\mathbb{E}_{\pi^{(r)}_k(\xi;\mu)}g_l(\xi,\mu)$ as well as the exact performance measure function $\mathbb{E}_{\pi_k(\xi;\mu)}g_l(\xi,\mu)$. We note that, for a specific $g_l(\xi,\mu)$, it may be possible to use a PWL approximation with fewer knots in the optimization model as $\mathbb{E}_{\pi^{(r)}_k(\xi;\mu)}g_l(\xi,\mu)$ is pre-computable.

\subsection{Summary of Approximation Algorithm.}
Our algorithm for solving the capacity sizing problem \eqref{opt:objective} is summarized in Algorithm \ref{alg1}.
Here $\varepsilon>0$ is an arbitrarily small number and controls the optimality gap; see Theorem \ref{theorem: FAO consistency}. Subsection \ref{sec:consistency} shows that, for an arbitrary $\varepsilon>0$, the solution obtained from Algorithm \ref{alg1} is $\varepsilon$-optimal to problem \eqref{opt:objective} as long as the finite approximation order $r$ and the PWL approximation order $m$ are sufficiently large. 

\begin{algorithm}[!htb]
	\caption{A finite and piecewise linear approximation algorithm for the constrained capacity sizing problem \eqref{opt:objective}.
	} \label{alg1}
	\begin{algorithmic}
		\State \textbf{Input}: some $r,m\in\mathbb{N}$ and $\varepsilon>0$.
		\State \textbf{Output}: an approximate solution $(\pmb{\mu}^{(r,m,\varepsilon)},\pmb{w}^{(r,m,\varepsilon)})$.
		\State \emph{Step 1}. For all $k\in\mathscr{K}$ and  $\mu_i^{(m)}(i=0,1,\dots,2^m)$, compute the finite-support stationary distribution $\pi^{(r)}_k(x;\mu_i^{(m)})$: \\
		\qquad Define a finite set of states $\{c_{k_j}^{(r)}\}_{j=1}^{J^{(r)}}$ and a transition matrix $\pmb{Q}_k^{(r)}(\mu_i^{(m)})$ via \eqref{eqn:def approximate MC}.\\
		\qquad Compute the stationary probability vector $\pmb{v}_k^{(r)}(\mu_i^{(m)})$ via \eqref{eqn:approx stationary vector}.\\
		\qquad Obtain the stationary distribution $\pi^{(r)}_k(x;\mu_i^{(m)})= \sum_{j=1}^{J^{(r)}}\mathds{1}\{x \geqslant c_{k_j}^{(r)} \} {v}_{k_j}^{(r)}(\mu_i^{(m)})$ via \eqref{eqn:approx stationary dist}.
		\State \emph{Step 2}. Construct the piecewise linear functions $\phi^{(r,m)}_{k,l}(\mu)$ for all $k\in\mathscr{K},l\in\mathscr{L}$ via \eqref{eqn:approx performance measure} and \eqref{eqn:pwl function}.
		\State \emph{Step 3}. Solve the piecewise linear program and obtain the solution $(\pmb{\mu}^{(r,m,\varepsilon)},\pmb{w}^{(r,m,\varepsilon)})$:
		\begin{align}
            \min_{\pmb{\mu},\pmb{w}} \quad & \pmb{p}^\text{T} \pmb{w} \label{opt:pwl objective}\tag{PWLP} \\
            \text{s.t.}\quad & \phi_{k,l}^{(r,m)}(\mu_k) \leqslant w_{k,l}+\frac{\varepsilon}{2},\quad (k,l)\in\mathscr{K}\times\mathscr{L}, \label{opt:pwl constraint relax} \tag{\ref{opt:pwl objective}-a}\\
            & \pmb{M} \pmb{w} \leqslant \pmb{d}, \label{opt:pwl constraint} \tag{\ref{opt:pwl objective}-b}\\
            & \pmb{\mu}\in [\mu_{\min},\mu_{\max}]^{K},\, \pmb{w}\in\mathbb{R}^{K\times L}.
            \label{opt:pwl constraint lu bound} \tag{\ref{opt:pwl objective}-c}
        \end{align}
	    \end{algorithmic} 
\end{algorithm}

\section{Main Results}\label{sec:main results}
The approximate solutions obtained from Algorithm \ref{alg1} are near-optimal in the following sense.
\begin{definition}[$\varepsilon$-optimal solution]\label{def:near optimal solution}
Consider  problem \eqref{opt:objective} and assume that it has an optimal solution $(\pmb{\mu}^*,\pmb{w}^*)$.
A vector $(\pmb{\mu}^\varepsilon, \pmb{w}^\varepsilon)$ with
$\pmb{\mu}^\varepsilon=(\mu^\varepsilon_k)_{k\in\mathscr{K}}\in[\mu_{\min},\mu_{\max}]^{K}$ and $\pmb{w}^\varepsilon=(w^\varepsilon_{k,l})_{k\in\mathscr{K},l\in\mathscr{L}}\in\mathbb{R}^{K\times L}$ is an $\varepsilon$-optimal solution to problem \eqref{opt:objective} if $(i)$ $\pmb{p}^\text{T} \pmb{w}^\varepsilon \leqslant \pmb{p}^\text{T} \pmb{w}^*$,
$(ii)$ $\mathbb{E}_{\pi_k(\xi;\mu^\varepsilon_k)}g_l(\xi,\mu^\varepsilon_k) \leqslant w_{k,l}^\varepsilon+\varepsilon,\, \forall (k,l)\in\mathscr{K}\times\mathscr{L}$, and
$(iii)$ $\pmb{M} \pmb{w}^\varepsilon \leqslant \pmb{d}$.
\end{definition}
In other words, an $\varepsilon$-optimal solution achieves an objective value no worse than the exact optimal solution by relaxing the stochastic constraints \eqref{opt:constraint relax} with a predefined gap $\varepsilon$.
Our convergence results are proved under the following assumptions.

\begin{assumption}[appropriate distributions] 
For all waitlist $k\in\mathscr{K}$, the associated distributions $A_k(x)$ and $G_k(x)$ satisfy the following: 
\begin{itemize}
    \item[$(i)$] $\exists a'_{k},a''_k>0$ such that $\big|\frac{dA_k(x) }{dx}\big|\leqslant a'_{k},$ 
    $\big|\frac{d^2A_k(x) }{dx^2}\big|\leqslant a''_{k}$ for all $x\in\mathbb{R}_{+}$.
    \item[$(ii)$] $\exists \bar{y}_k>0$ such that $G_k(\bar{y}_k)=1$.
\end{itemize}
\label{A:input distributions 1}
\end{assumption}

\begin{assumption}[appropriate performance functions]
For all $g_l(\xi,\mu)\,(l\in \mathscr{L})$ and  $\mu\in[\mu_{\min},\mu_{\max}]$, we have the following:
\begin{itemize}
    \item[$(i)$] $g_l(\xi,\mu)=g_l(0,\mu)$ for all $\xi<0$, and $g_l(\xi,\mu)=g_l(\bar{y},\mu)$ for all $\xi> \bar{y}$, where $\bar{y} := \max_{k\in\mathscr{K}}\bar{y}_k$.
    \item[$(ii)$] $V_\xi\big(g_l(\xi,\mu)\big)<\infty$, $i.e.$, $g_l(\xi,\mu)$ has a bounded variation w.r.t. $\xi$.
    \item[$(iii)$] $\lim_{\mu'\rightarrow \mu}\sup_{\xi\in\mathbb{R}_+}\big|g_l(\xi,\mu')-g_l(\xi,\mu)\big|=0$, $i.e.$, $g_l(\xi,\mu)$ is continuous w.r.t. $\mu$ under a supremum norm. 
\end{itemize}
\label{A:input distributions 2}
\end{assumption}
Here operation $V(f)$ denotes total variation of a function  $f:\mathbb{R} \mapsto\mathbb{R}$:
\begin{align}
    V(f):= \sup_{x_1<x_2<\dots<x_n, x_1,x_2,\dots,x_n\in\mathbb{R};n\in\mathbb{N}, n\geqslant 2 } \sum_{i=2}^n|f(x_i)-f(x_{i-1})|, \nonumber
\end{align}
where the supremum is over all partitions of $\mathbb{R}$. 
When $f'(x)=\frac{df(x)}{dx}$ is well defined, $V(f)= \int_{\mathbb{R}} |\frac{df(x)}{dx}|dx = \int_{\mathbb{R}} |f'(x)|dx.$
Paricularly, for $f(x,u):\mathbb{R}^n\times \mathbb{R}\mapsto\mathbb{R}$, we use $V_u(f(x,u))$ to denote the total variation of $f(x,u)$ as a single variable function of $u$ with a fixed $x\in \mathbb{R}^n$.

In Assumption~\ref{A:input distributions 1}, $(i)$ assumes the smoothness of inter-arrival time distribution.
Examples include mixture of exponential distributions, Erlang distributions, Gamma distributions with a shape parameter no smaller than 2, and distributions taking a polynomial form; Appendix \ref{Append:Arrivial distribution} summarizes common distributions satisfying this assumption. Item $(ii)$ assumes patience time is upper bounded, which is a natural constraint given that patients have bounded life expectancy.



In Assumption~\ref{A:input distributions 2}, $(i)$ is not restrictive because the offered waiting time $\xi$ is supported on $[0,\bar{y}]$ and $g_l(\xi,\mu)$ can be redefined to take any value for $\xi$ outside support.
Item $(ii)$ has two different sufficient conditions which can be verified easily. 
The first condition is that $g_l(\xi,\mu)$ is monotonic or convex/concave $w.r.t.$ $\xi\in[0,\bar{y}]$. 
The second condition is that $\frac{dg_l(\xi,\mu)}{d\xi}$ is bounded for $\xi\in[0,\bar{y}]$, which suggests   $V_\xi\big(g_l(\xi,\mu)\big)=\int_{0}^{\bar{y}} |\frac{dg_l(\xi,\mu)}{d\xi}|d\xi\leqslant \bar{y}\cdot \sup_{\xi\in[0,\bar{y}]}|\frac{dg_l(\xi,\mu)}{d\xi}|<\infty.$
Item $(iii)$ also has two different sufficient conditions.
The first condition is that $g_l(\xi,\mu)$ is independent of $\mu$, which suggests $\sup_{\xi\in\mathbb{R}_+}\big|g_l(\xi,\mu')-g_l(\xi,\mu)\big|=0$. 
The second condition is that $\frac{dg_l(\xi,\mu)}{d\mu}$ is bounded for $\xi\in[0,\bar{y}]$ and $\mu\in[\mu_{\min},\mu_{\max}]$, which suggests $\lim_{\mu'\rightarrow\mu}\sup_{\xi\in\mathbb{R}_+}\big|g_l(\xi,\mu')-g_l(\xi,\mu)\big|=\lim_{\mu'\rightarrow\mu}\sup_{\xi\in[0,\bar{y}]}\big|g_l(\xi,\mu')-g_l(\xi,\mu)\big|\leqslant
\lim_{\mu'\rightarrow\mu} |\mu'-\mu|\cdot\sup_{\xi\in[0,\bar{y}],\tilde\mu\in[\mu_{\min},\mu_{\max}]}\big|\frac{dg_l(\xi,\tilde\mu)}{d\tilde\mu}\big|=0.$

\subsection{Consistency Results }\label{sec:consistency}
\begin{theorem}[Consistency of Algorithm \ref{alg1}]\label{theorem: FAO consistency}
Under Assumptions \ref{A:input distributions 1}--\ref{A:input distributions 2}, for all $\varepsilon>0$, there exist $r_0$ and $ m_0$ such that the solution $(\pmb{\mu}^{(r,m,\varepsilon)},\pmb{w}^{(r,m,\varepsilon)})$ obtained from Algorithm \ref{alg1} is an $\varepsilon$-optimal solution to capacity sizing problem \eqref{opt:objective} for all $r>r_0,m>m_0$.
\end{theorem}

To establish this theorem, we analyze the solution gap between  the original problem \eqref{opt:objective} and the approximate problem \eqref{opt:pwl objective}. We decompose the gap into two components: finite approximation error and PWL approximation error. We leverage the following two lemmas, whose proofs are
relegated to Appendix \ref{append:proof}. 
Lemma \ref{lemma:fa uniform convergence} below shows that the finite approximate performance measures uniformly converge to the exact performance measures. This lemma suggests that the finite approximation error diminishes as the number of states used increases, $i.e.$, as $r$ goes to $+\infty$.
Lemma \ref{lemma:pm continuity} shows that the exact performance measures are continuous functions of the service rates, and furthermore, these functions are uniformly continuous.
As more knots are used in the PWL approximation, $i.e.$, as $m$ goes to $+\infty$, the difference between two neighboring knots decreases. 
Then the generated PWL function is closer to the function of exact performance measures due to the later's uniform continuity from Lemma \ref{lemma:pm continuity}.
\begin{lemma}[Uniform convergence of finite approximation]\label{lemma:fa uniform convergence}
Under Assumptions \ref{A:input distributions 1}--\ref{A:input distributions 2}, for all $\varepsilon>0$, there exists $r_0\in\mathbb{N}$ such that 
\begin{align*}
    \sup_{\mu\in[\mu_{\min},\mu_{\max}]}\Big|\mathbb{E}_{\pi_k(\xi;\mu)}g_l(\xi,\mu)- \mathbb{E}_{\pi^{(r)}_k(\xi;\mu)}g_l(\xi,\mu)\Big|\leqslant \varepsilon,\quad   k\in\mathscr{K},l\in\mathscr{L}, r> r_0.
\end{align*}
\end{lemma}

\begin{lemma}[Continuity of performance measures]\label{lemma:pm continuity}
Under Assumptions \ref{A:input distributions 1}--\ref{A:input distributions 2}, for all $k\in\mathscr{K}, l\in\mathscr{L}$, the exact performance measure $\mathbb{E}_{\pi_k(\xi;\mu)}g_l(\xi,\mu)$ as a function of rate $\mu$ is continuous on $[\mu_{\min},\mu_{\max}]$.
\end{lemma}

Now we formally complete the proof of Theorem \ref{theorem: FAO consistency} based on Lemmas \ref{lemma:fa uniform convergence}--\ref{lemma:pm continuity}.
\begin{proof}{Proof of Theorem \ref{theorem: FAO consistency}:}
Consider an arbitrarily small $\varepsilon>0$.
We first analyze the gap between constraints \eqref{opt:constraint relax} in the original problem \eqref{opt:objective} and constraints \eqref{opt:pwl constraint relax} in the approximate problem \eqref{opt:pwl objective}. 
By Lemma \ref{lemma:pm continuity}, for all $k\in\mathscr{K},l\in\mathscr{L}$, performance measure $\mathbb{E}_{\pi_k(\xi;\mu)}g_l(\xi,\mu)$ is a continuous function of $\mu\in[\mu_{\min},\mu_{\max}].$
Because the domain $[\mu_{\min},\mu_{\max}]$ is compact, $\mathbb{E}_{\pi_k(\xi;\mu)}g_l(\xi,\mu)$ is further uniformly continuous. 
Therefore, there exists $\delta_{k,l}>0$ such that for all $\mu,\mu'\in [\mu_{\min},\mu_{\max}]$, 
\begin{align}
    |\mu'-\mu|\leqslant \delta_{k,l}\quad  \Rightarrow \quad \Big|\mathbb{E}_{\pi_k(\xi;\mu')}g_l(\xi,\mu')-\mathbb{E}_{\pi_k(\xi;\mu)}g_l(\xi,\mu)\Big|\leqslant \frac{\varepsilon}{4}.
    \label{eqn:uniform continuity}
\end{align}
By Lemma \ref{lemma:fa uniform convergence}, there exists $r_0$ such that for all $k\in\mathscr{K}, l\in\mathscr{L}, r>r_0$, 
\begin{align}
    \sup_{\mu\in[\mu_{\min},\mu_{\max}]}\Big|\mathbb{E}_{\pi_k(\xi;\mu)}g_l(\xi,\mu)-\mathbb{E}_{\pi^{(r)}_k(\xi;\mu)}g_l(\xi,\mu)\Big|\leqslant \frac{\varepsilon}{4}.
    \label{eqn:uniform convergence}
\end{align}
Define $m_0:=\min\{m\in\mathbb{N}|\frac{\mu_{\max}-\mu_{\min}}{2^m}\leqslant \delta_{k,l},\forall k\in\mathscr{K},l\in\mathscr{L}\}.$
Then for all $k\in\mathscr{K}, l\in\mathscr{L}, r>r_0, m>m_0, \mu\in [\mu^{(m)}_{i-1},\mu^{(m)}_{i}], i\in \{1,2,\dots,2^m\}$, we have
\begin{align}
    &\Big| \phi_{k,l}^{(r,m)}(\mu) - \mathbb{E}_{\pi_k(\xi;\mu)}g_l(\xi,\mu) \Big| \nonumber \\
    \leqslant &\max \Bigg\{\Big|\phi_{k,l}^{(r,m)}(\mu^{(m)}_{i-1}) - \mathbb{E}_{\pi_k(\xi;\mu)}g_l(\xi,\mu)\Big|, \Big|\phi_{k,l}^{(r,m)}(\mu^{(m)}_i) - \mathbb{E}_{\pi_k(\xi;\mu)}g_l(\xi,\mu)\Big|\Bigg\} \nonumber\\
    & \text{(due to $\mu\in [\mu^{(m)}_{i-1},\mu^{(m)}_{i}]$ and the linearity of $\phi_{k,l}^{(r,m)}(\mu)$ on $[\mu^{(m)}_{i-1},\mu^{(m)}_{i}]$)}\nonumber\\
    = &  \max\Bigg\{\Big|\mathbb{E}_{\pi^{(r)}_k(\xi;\mu^{(m)}_{i-1})}g_l(\xi,\mu^{(m)}_{i-1}) - \mathbb{E}_{\pi_k(\xi;\mu)}g_l(\xi,\mu)\Big|, \Big|\mathbb{E}_{\pi^{(r)}_k(\xi;\mu^{(m)}_i)}g_l(\xi,\mu^{(m)}_i) - \mathbb{E}_{\pi_k(\xi;\mu)}g_l(\xi,\mu)\Big|\Bigg\}  \nonumber\\ 
    & \text{(by equality \eqref{eqn:pwl function})} \nonumber \\
    \leqslant &  \max\Bigg\{\Big|\mathbb{E}_{\pi_k(\xi;\mu^{(m)}_{i-1})}g_l(\xi,\mu^{(m)}_{i-1}) - \mathbb{E}_{\pi_k(\xi;\mu)}g_l(\xi,\mu)\Big|+\Big|\mathbb{E}_{\pi^{(r)}_k(\xi;\mu^{(m)}_{i-1})}g_l(\xi,\mu^{(m)}_{i-1}) - \mathbb{E}_{\pi_k(\xi;\mu^{(m)}_{i-1})}g_l(\xi,\mu^{(m)}_{i-1}) \Big|, \nonumber\\ 
    &  \Big|\mathbb{E}_{\pi_k(\xi;\mu^{(m)}_i)}g_l(\xi,\mu^{(m)}_i)  - \mathbb{E}_{\pi_k(\xi;\mu)}g_l(\xi,\mu)\Big|+\Big|\mathbb{E}_{\pi^{(r)}_k(\xi;\mu^{(m)}_i)}g_l(\xi,\mu^{(m)}_i)  - \mathbb{E}_{\pi_k(\xi;\mu^{(m)}_i)}g_l(\xi,\mu^{(m)}_i)\Big|\Bigg\} \nonumber \\
    &  \text{(by the triangle inequality)} \nonumber \\
    \leqslant &  \max\Bigg\{\Big|\mathbb{E}_{\pi_k(\xi;\mu^{(m)}_{i-1})}g_l(\xi,\mu^{(m)}_{i-1}) - \mathbb{E}_{\pi_k(\xi;\mu)}g_l(\xi,\mu)\Big|+\frac{\varepsilon}{4}, \Big|\mathbb{E}_{\pi_k(\xi;\mu^{(m)}_i)}g_l(\xi,\mu^{(m)}_i)  - \mathbb{E}_{\pi_k(\xi;\mu)}g_l(\xi,\mu)\Big|+\frac{\varepsilon}{4}\Bigg\} \nonumber\\ 
    & \text{(due to inequality \eqref{eqn:uniform convergence} and $r>r_0$)} \nonumber \\
    \leqslant &  \max \{\frac{\varepsilon}{4} + \frac{\varepsilon}{4}, \frac{\varepsilon}{4} + \frac{\varepsilon}{4} \} =\frac{\varepsilon}{2}. \label{eqn:FA PWLA  error crude} \\
    &\text{(due to inequality \eqref{eqn:uniform continuity}, $|\mu-\mu^{(m)}_{i-1}|\leqslant \frac{\mu_{\max}-\mu_{\min}}{2^m}\leqslant \delta_{k,l}$, and $|\mu-\mu^{(m)}_{i}|\leqslant \frac{\mu_{\max}-\mu_{\min}}{2^m}\leqslant \delta_{k,l}$)}\nonumber
\end{align}
Because $[\mu_{\min},\mu_{\max}]=\cup_{i=1,2,\dots,2^m}[\mu^{(m)}_{i-1},\mu^{(m)}_{i}]$, inequality \eqref{eqn:FA PWLA  error crude} can also be stated as
\begin{align}
    \sup_{\mu\in[\mu_{\min},\mu_{\max}]} \Big| \phi_{k,l}^{(r,m)}(\mu) - \mathbb{E}_{\pi_k(\xi;\mu)}g_l(\xi,\mu) \Big| \leqslant \frac{\varepsilon}{2}, \quad  k\in\mathscr{K}, l\in\mathscr{L}, r>r_0, m>m_0. \label{eqn:FA PWLA  error}
\end{align}

Now we use \eqref{eqn:FA PWLA  error} to prove $\varepsilon$-optimality of solution $(\pmb{\mu}^{(r,m,\varepsilon)},\pmb{w}^{(r,m,\varepsilon)})$ obtained by Algorithm \ref{alg1} with $r>r_0$ and $m>m_0$; that is, we prove $(\pmb{\mu}^{(r,m,\varepsilon)},\pmb{w}^{(r,m,\varepsilon)})$ satisfies conditions $(i)$--$(iii)$ of Definition \ref{def:near optimal solution}.

$(i)$. Let $(\pmb{\mu}^*,\pmb{w}^*)=\big(({\mu}_k^*)_{k\in\mathscr{K}},({w}^*_{k,l})_{k\in\mathscr{K},l\in\mathscr{L}}\big)$ be the optimal solution of problem \eqref{opt:objective} and $(\pmb{\mu}^{(r,m,\varepsilon)},\pmb{w}^{(r,m,\varepsilon)})=\big(({\mu}_k^{(r,m,\varepsilon)})_{k\in\mathscr{K}},({w}^{(r,m,\varepsilon)}_{k,l})_{k\in\mathscr{K},l\in\mathscr{L}}\big)$ be the solution of problem \eqref{opt:pwl objective}.  
Then $(\pmb{\mu}^*,\pmb{w}^*)$ satisfies \eqref{opt:constraint relax}, and we have $\mathbb{E}_{\pi_k(\xi;\mu^*_k)}g_l(\xi,\mu^*_k) \leqslant w^*_{k,l}$ for all $(k,l)\in\mathscr{K}\times\mathscr{L}.$
Then inequality \eqref{eqn:FA PWLA  error} suggests that  $\phi_{k,l}^{(r,m)}(\mu^*_k)  \leqslant \mathbb{E}_{\pi_k(\xi;\mu^*_k)}g_l(\xi,\mu^*_k) +\frac{\varepsilon}{2} \leqslant w_{k,l}^*+\frac{\varepsilon}{2}$ for all $(k,l)\in\mathscr{K}\times\mathscr{L}.$
Thus, $(\pmb{\mu}^*,\pmb{w}^*)$ also satisfies \eqref{opt:pwl constraint relax}. 
Because $(\pmb{\mu}^*,\pmb{w}^*)$ satisfies \eqref{opt:constraint} and \eqref{opt:constraint lu bound} in problem \eqref{opt:objective} by definition, it naturally satisfies  \eqref{opt:pwl constraint} and \eqref{opt:pwl constraint lu bound} in problem \eqref{opt:pwl objective} as these constraints are identical.
Thus, $(\pmb{\mu}^*,\pmb{w}^*)$ is a feasible solution to problem \eqref{opt:pwl objective}.
Since $(\pmb{\mu}^{(r,m,\varepsilon)},\pmb{w}^{(r,m,\varepsilon)})$ is the optimal solution to problem \eqref{opt:pwl objective}, we have $\pmb{p}^{\text{T}} \pmb{w}^{(r,m,\varepsilon)} \leqslant \pmb{p}^{\text{T}} \pmb{w}^*$.

$(ii)$. Because $(\pmb{\mu}^{(r,m,\varepsilon)},\pmb{w}^{(r,m,\varepsilon)})$ is a solution to problem \eqref{opt:pwl objective}, it satisfies  \eqref{opt:pwl constraint relax} and we have  $\phi_{k,l}^{(r,m)}(\mu^{(r,m,\varepsilon)}_k)  \leqslant w_{k,l}^{(r,m,\varepsilon)}+\frac{\varepsilon}{2},\, (k,l)\in\mathscr{K}\times\mathscr{L}.$
Then inequality \eqref{eqn:FA PWLA  error} suggests that $\mathbb{E}_{\pi_k(\xi;\mu^{(r,m,\varepsilon)}_k)}g_l(\xi,\mu^{(r,m,\varepsilon)}_k)  \leqslant \phi_{k,l}^{(r,m)}(\mu^{(r,m,\varepsilon)}_k)  +\frac{\varepsilon}{2} \leqslant w_{k,l}^{(r,m,\varepsilon)}+\frac{\varepsilon}{2}+\frac{\varepsilon}{2} = w_{k,l}^{(r,m,\varepsilon)}+\varepsilon$, $(k,l)\in\mathscr{K}\times\mathscr{L}$.

$(iii)$ also holds because it is exactly constraint \eqref{opt:pwl constraint}. 
\qedwhite
\end{proof}

\begin{remark}[Comparison with \citealt{li2023new}]
Our Lemmas \ref{lemma:fa uniform convergence}--\ref{lemma:pm continuity} significantly extend the findings of \cite{li2023new}, who focus solely on performance evaluation of a fixed Markov chain and do not account for varying decision parameters.
Lemma \ref{lemma:fa uniform convergence} on the uniform convergence ($w.r.t.$ service rate $\mu$) of finite approximate performance measures  is an extension of the \emph{pointwise} convergence developed in \cite{li2023new} and is significantly stronger because we construct an infinite-dimensional approximation, $i.e.$, over all $\mu\in[\mu_{\min},\mu_{\max}]$.
We also provide a brief overview of how our results build upon \cite{li2023new}.   
First, we show that Assumption \ref{A:input distributions 1} establishes sufficient conditions for the pointwise convergence of finite approximate performance measures, building on the work of \cite{li2023new}. 
Next, Assumption \ref{A:input distributions 2} indicates equicontinuity ($w.r.t.$ service rate $\mu$) of finite approximate performance measures, which enhances the pointwise convergence of finite approximate performance measures to uniform convergence and is \emph{not} covered by \cite{li2023new}.
Then the equicontinuity indicates the continuity of finite approximate performance measures 
and thereby indicates the continuity of exact performance measures as limit functions, which is the result of Lemma \ref{lemma:pm continuity}.
\end{remark}

\section{Validity of Finite Approximation and  Computational  Complexity}
\label{sec:alg property}
In this section, we examine the validity and computational complexity of Algorithm \ref{alg1} to achieve a desired precision $\varepsilon$. 
As a reminder, Algorithm \ref{alg1} has three steps for a given set of parameters $(r,m,\varepsilon)$. 
Step 1 employs the finite approximation, using $J^{(r)}=2^r+1$ states. 
It involves computing stationary probability vectors for $KN^{(m)}$ transition matrices, each comprising $J^{(r)}$ states. 
Step 2 employs the PWL approximation, using $N^{(m)}=2^m+1$ knots.
It involves computing inner products for $KLN^{(m)} $ pairs of $J^{(r)}$-dimension vectors. 
Step 3 solves a program with $KL$ PWL constraints, each having $N^{(m)}$ knots.
Overall, the computation complexity of Algorithm \ref{alg1} is determined by $J^{(r)}$ and $N^{(m)}$, $i.e.$, the numbers of finite approximation states and PWL approximation knots. 
Subsection \ref{sec:compute stationary} optimizes the computation of Step 1 through matrix products, which can be readily parallelized on multiple processors. 
Subsection \ref{sec:FAO gap} outlines how to numerically evaluate the system performance under the solution service rate vector $\pmb{\mu}^{(r,m,\varepsilon)}$ obtained from Algorithm \ref{alg1}.
Additionally, it details how to obtain an error bound associated with the numerical evaluation.
Subsection \ref{sec:convergence rate} addresses the question of determining the minimum number of finite approximation states and PWL approximation knots required to achieve a given optimality gap $\varepsilon>0$.
This is important for Algorithm \ref{alg1}, as it impacts the solution accuracy and the computational resources needed.


\subsection{Validity of Finite Approximation}\label{sec:compute stationary}
\begin{theorem}\label{theorem:MC compute}
Under Assumption \ref{A:input distributions 1}, there exists $r^*\in\mathbb{N}$ such that for all $k\in\mathscr{K}$, $\mu\in[\mu_{\min},\mu_{\max}]$, and $r> r^*$, the finite-state Markov chain $\pmb{\Psi}^{(r)}_k(\mu)$, whose states $\{c^{(r)}_{k_i}\}_{i=1}^{J^{(r)}}$ and transition matrix $\pmb{Q}_k^{(r)}(\mu)$ are defined in \eqref{eqn:def approximate MC}, has a unique stationary probability vector $\pmb{v}^{(r)}_k(\mu)$ satisfying \eqref{eqn:approx stationary vector}. 
Moreover, $\pmb{v}^{(r)}_k(\mu) = \lim_{n\rightarrow \infty} \big(\pmb{Q}^{(r)^{\textnormal{T}}}_k(\mu)\big)^n\pmb{v}_0$, where $\pmb{v}_0\in\mathbb{R}_+^{J^{(r)}}$ is an arbitrary probability vector.
\end{theorem}
\begin{proof}{Proof:}
Assumption \ref{A:input distributions 1} indicates that for all waitlist $k\in\mathscr{K}$, there exists $\varepsilon_k>0$ such that
\begin{align}
   \delta_k:=\mathbb{P}[t_k^{n+1}-s_k^n\geqslant \varepsilon_k]>0
    \label{ieq:arrival time gap}
\end{align}
for all service rate $\mu_k=\mu\in[\mu_{\min},\mu_{\max}]$. 
Here $n\in\mathbb{N}$  indexes an arbitrary patient.
We define 
\begin{align}
    r^*:=\min\Big\{r\in\mathbb{N}\,\Big|\, \frac{\bar{y}_k}{2^r}\leqslant \varepsilon_k,k\in\mathscr{K}\Big\}. \label{eqn:min r for MC}
\end{align}
We will show  $r^*$ satisfies the requirements of Theorem \ref{theorem:MC compute}. 
Specifically, we prove that, for all $k\in\mathscr{K}$, $\mu\in[\mu_{\min},\mu_{\max}]$, and $r>r^*$, the finite-state Markov chain $\pmb{\Psi}^{(r)}_k(\mu)$ supported on $\{c^{(r)}_{k_i}\}_{i=1}^{J^{(r)}}$ with transition matrix $\pmb{Q}^{(r)}_k(\mu)=\{q^{(r)}_{k_{i,j}}(\mu)\}_{i,j= 1,2,\dots,J^{(r)}}$ satisfies the following two properties:
\begin{itemize}
    \item[$(a)$] State $0$ is accessible for all states $\{c^{(r)}_{k_i}\}_{i=1}^{J^{(r)}}.$
    \item[$(b)$] The communicating class $S\subseteq\{c^{(r)}_{k_i}\}_{i=1}^{J^{(r)}}$ that includes state $0$ ($i.e.$, the states that are mutually accessible with state 0) corresponds to an ergodic transition matrix. 
    In other words, $\{q^{(r)}_{k_{i,j}}\}_{i,j\in S}$ is an irreducible, positive recurrent, and aperiodic transition matrix.
\end{itemize}
Properties $(a)$ and $(b)$ indicate that the stationary probability vector of Markov chain $\pmb{\Psi}^{(r)}_k(\mu)$ is equivalent to the stationary probability vector of transition matrix $\{q^{(r)}_{k_{i,j}}(\mu)\}_{i,j\in S}$. 
Since transition matrix $\{q^{(r)}_{k_{i,j}}(\mu)\}_{i,j\in S}$ is ergodic, it has a unique stationary probability vector, which is also the unique stationary probability vector  $\pmb{v}^{(r)}_k(\mu)$ of Markov chain $\pmb{\Psi}^{(r)}_k(\mu)$. 
Because state $0$ is accessible for all states of $\pmb{\Psi}^{(r)}_k(\mu)$, any initial vector $\pmb{v}_0$ will be  absorbed onto  $S$ and converge to $\pmb{v}^{(r)}_k(\mu)$, $i.e.$ $\pmb{v}^{(r)}_k(\mu) = \lim_{n\rightarrow \infty} \big(\pmb{Q}^{(r)^{\textnormal{T}}}_k(\mu)\big)^n\pmb{v}_0$. Now we formally prove properties (a) and (b).

$(a)$. 
Recall that for Markov chain $\pmb{\Psi}_k(\mu)$ described in \eqref{eqn:original MC}, we have inequality \eqref{ieq:arrival time gap}.
In other words, with probability $\delta_k>0$, the current patient's service time is less than the inter-arrival time by $\varepsilon_k$.
Event $\{t_k^{n+1}-s_k^n\geqslant \varepsilon_k\}$ indicates 
$\{ \xi_k^{n+1} \leqslant [\xi_k^{n} -\varepsilon_k]_+\}$.
In other words, when the current patient's service time is less than the inter-arrival time by $\varepsilon_k$, the next patient's offered waiting time will be less than that of the current patient by $\varepsilon_k$.
This is because the next patient's offered waiting time equals the current patient's offered waiting time, plus the time for which the current patient occupies the server, and minus the inter-arrival time. 
Moreover, the current patient's occupation time is bounded by the required service time.
Thus, conditional on any given $\xi^n_k=u\in [0,\bar{y}_k]$, 
\begin{align}
    \mathbb{P}\big[\xi_k^{n+1} \leqslant [\xi_k^{n} -\varepsilon_k]_+\,|\,\xi^n_k=u\big] \geqslant \delta_k,\quad u\in [0,\bar{y}_k]. \label{eqn:drop prob 1}
\end{align}
Now consider the finite-state Markov chain $\pmb{\Psi}^{(r)}_k(\mu)$, which is obtained by truncating the original Markov chain $\pmb{\Psi}_k(\mu)$.
We have that, for all $i\in\{2,3,\dots,J^{(r)}\}$, 
\begin{align}
    \sum_{j=1}^{i-1} q^{(r)}_{k_{i,j}}(\mu) = &  \tau_k(c^{(r)}_{k_{i-1}},c^{(r)}_{k_{i}};\mu) \qquad \text{(due to definition \eqref{eqn:def approximate MC})} \nonumber\\
    = & \mathbb{P}[ \xi_k^{n+1}  \leqslant c^{(r)}_{k_{i-1}} \,|\, \xi_k^{n} = c^{(r)}_{k_{i}} ]   \qquad \text{(by definition of transition kernel $\tau_k$)} \nonumber\\
    \geqslant  & \mathbb{P}\big[\xi_k^{n+1} \leqslant [\xi_k^{n} -\varepsilon_k]_+\,|\, \xi_k^{n} = c^{(r)}_{k_{i}}\big] \qquad \text{(due to $|c^{(r)}_{k_{i}}-c^{(r)}_{k_{i-1}}|=\frac{\bar{y}_k}{2^r}\leqslant \varepsilon_k$ by definition \eqref{eqn:min r for MC})}\nonumber\\
    \geqslant & \delta_k >0. \qquad \text{(due to \eqref{eqn:drop prob 1})} \label{eqn:transit left}
\end{align}
In other words, with probability at least $\delta_k>0$, Markov chain $\pmb{\Psi}^{(r)}_k(\mu)$ will transit (leftwards) to states with strictly smaller index $i$. 
Thus, state 0 is accessible for all states $\{c^{(r)}_{k_i}\}_{i=1}^{J^{(r)}}.$

$(b)$. $\{q^{(r)}_{k_{i,j}}(\mu)\}_{i,j\in S}$ is an irreducible transition matrix because $S$ is a communicating class and $0\in S$ is accessible for all states $\{c^{(r)}_{k_i}\}_{i=1}^{J^{(r)}}$.
$\{q^{(r)}_{k_{i,j}}(\mu)\}_{i,j\in S}$ is positive recurrent because it is irreducible and has finite states.
Then we only need to prove aperiodicity. 
Consider state $0$.  
Using a similar argument to \eqref{eqn:transit left}, we have $q^{(r)}_{k_{1,1}}(\mu) =   \tau_k(0,0;\mu)
    =  \mathbb{P}[ \xi_k^{n+1}  \leqslant 0 |\, \xi_k^{n} = 0 ]  
    =   \mathbb{P}\big[\xi_k^{n+1} \leqslant [\xi_k^{n} -\varepsilon_k]_+\,|\, \xi_k^{n} = 0\big]
    \geqslant \delta_k >0.$
In other words,  with probability at least $\delta_k>0$, Markov chain $\pmb{\Psi}^{(r)}_k(\mu)$ will stay at state $c^{(r)}_{k_1}=0$ if it is currently at state $0$.
Therefore, state 0 has a period of $1$ and the transition matrix $\{q^{(r)}_{k_{i,j}}(\mu)\}_{i,j\in S}$ is aperiodic.
\qedwhite
\end{proof}

\subsection{Error Bounds for Evaluating System Performance}
\label{sec:FAO gap}
The following theorem provides an error bound on a waitlist's finite approximate performance measure. 
The error bound is computable by solving a linear program.
\begin{theorem}[computable error bound on performance measure]\label{theorem:FAO optimality gap}
Under Assumptions \ref{A:input distributions 1}--\ref{A:input distributions 2}, 
we have the following error bound for all $k\in\mathscr{K}, l\in\mathscr{L}$, $\mu\in[\mu_{\min},\mu_{\max}]$, and $r>r^*$. 
\begin{align}
    & \Big | \mathbb{E}_{\pi_k(\xi;\mu)}g_l(\xi,\mu) -  \mathbb{E}_{\pi^{(r)}_k(\xi;\mu)}g_l(\xi,\mu) \Big |\leqslant \frac{\bar{y}_ka'_k}{2^{r-2}}\cdot V_\xi\big(g_l(\xi,\mu)\big)\cdot  e(k,r,\mu),  \label{eqn:li error bound}
\end{align}
where $e(k,r,\mu):= {\max_{s=0,1,\dots,J^{(r)}} \frac{1}{z^*_s}}$ and $z^*_s$ is obtained from the following linear program:
\begin{align}
    &\min_{z_s\in\mathbb{R},\{a_j\}_{j=1}^{J^{(r)}} \in[-1,1]^{J^{(r)}} } \quad\, z_s  \nonumber\\
    &s.t.  \quad \,  
    -z_s  \leqslant a_j+ \eta \delta_{js}(1-a_j) -\sum_{i=1}^{J^{(r)}} \big(1+\sum_{t=1}^{j}q^{(r)}_{k_{i,t}}(\mu)\big)\cdot(a_i-a_{i-1}) \leqslant z_s,\,\quad j=0,1,\dots,J^{(r)}, \eta\in \{0,1\}.  \nonumber
\end{align} 
Here $\{q^{(r)}_{k_{i,j}}(\mu)\}_{i,j=1,2,\dots,J^{(r)}}$ are defined in \eqref{eqn:def approximate MC}, $a_0\equiv 0$, $\delta_{js}:=\mathds{1}\{j=s\}$, and $\sum_{t=1}^0x_t:=0$ for all $x_t$.
\end{theorem}
The proof is provided in Appendix \ref{append:proof}. 
For any waitlist $k\in\mathscr{K}$, under performance function $g_l\,(l\in\mathscr{L})$ and service rate $\mu^{(r,m,\epsilon)}_k$, we can approximate the waitlist's steady-state performance measure as $\mathbb{E}_{\pi_k\left(\xi;\mu^{(r,m,\epsilon)}_k\right)}g_l(\xi,\mu^{(r,m,\epsilon)}_k)\approx \mathbb{E}_{\pi^{(r)}_k\left(\xi;\mu^{(r,m,\epsilon)}_k\right)}g_l(\xi,\mu^{(r,m,\epsilon)}_k)$.
By plugging $\mu=\mu^{(r,m,\epsilon)}_k$ into \eqref{eqn:li error bound}, we can obtain the associated approximation error bound.

\subsection{Convergence Rate} \label{sec:convergence rate}
We make an extra assumption  for the result stated in this subsection. This assumption requires that the performance function is Lipschitz continuous $w.r.t.$ its service rate parameter.
\begin{assumption}[Lipschitz continuous performance functions]
For each performance function $g_l$ $(l\in\mathscr{L})$, we have $(i)$ $\sup_{\mu\in[\mu_{\min},\mu_{\max}]}V_\xi(g_l(\xi,\mu))<\infty$, and $(ii)$ there exists $\zeta_l>0$ such that  $\big|g_l(\xi,\mu')-g_l(\xi,\mu)\big|\leqslant \zeta_l|\mu'-\mu|$ for all $\xi\in\mathbb{R}_+$ and $\mu,\mu'\in[\mu_{\min},\mu_{\max}]$.
\label{A:input distributions 3}
\end{assumption}
 We note that the examples mentioned for $(iii)$ of Assumption \ref{A:input distributions 2} satisfy $(ii)$ of Assumption \ref{A:input distributions 3}. 
One sufficient condition for $(i)$ of Assumption \ref{A:input distributions 3} is that $\frac{dg_l(\xi,\mu)}{d\xi}$ is bounded for $\xi\in\mathbb{R}_+$ and $\mu\in[\mu_{\min},\mu_{\max}].$  

The following theorem shows that using  $O(\varepsilon^{-1})$ finite approximation states and $O(\varepsilon^{-1})$ PWL approximation knots in Algorithm \ref{alg1} is sufficient to yield an $\varepsilon$-optimal solution.
\begin{theorem}[convergence rate of Algorithm \ref{alg1}]\label{theorem: FAO convergence rate}
Under Assumptions \ref{A:input distributions 1}--\ref{A:input distributions 3}, there exist constants $C_F$ and $C_P$ such that for all $\varepsilon>0$, the solution $(\pmb{\mu}^{(r,m,\varepsilon)},\pmb{w}^{(r,m,\varepsilon)})$ obtained from Algorithm \ref{alg1}  is an $\varepsilon$-optimal solution to problem \eqref{opt:objective} as long as the number of finite approximation states $J^{(r)}=2^r+1\geqslant \frac{C_F}{\varepsilon}$ and the number of PWL approximation knots $N^{(m)}=2^m+1\geqslant \frac{C_P}{\varepsilon}$.
\end{theorem}
The proof is provided in Appendix \ref{append:proof}. 
We discuss the underlying intuitions.
As previously analyzed in Subsection \ref{sec:consistency}, the optimality gap depends on the finite approximation error and the PWL approximation error.
We expect the finite approximation error to be in the order of $O(\frac{1}{J^{(r)}})$ due to the use of a finite-support probability distribution with evenly distributed states to approximate the original stationary distribution that is predominantly continuous on its support.
Similarly, the PWL approximation error is anticipated to be in the order of $O(\frac{1}{N^{(m)}})$  because a PWL function is used to approximate the original continuous constraint function, and the knots are also evenly distributed. 
Therefore, we can expect that using $O(\varepsilon^{-1})$ finite approximation states and $O(\varepsilon^{-1})$ PWL approximation knots will be sufficient to achieve an optimality gap of $\varepsilon$.

\section{Computational Results}
\label{sec:numerical}
We now present results from numerical experiments using recent liver transplant data to assess the effectiveness of our finite approximation method in evaluating waitlist performance and optimizing liver allocation. Specifically, we examine the optimal liver allocation rates for waitlists of various scales and patient mortality risks. 
All computations were conducted on an 8-core Intel platform with a 3.6GHz processor and 32G memory.

\subsection{Data and Model Setup}
\label{sec:exp data and parameter}
We use data obtained from OPTN for waitlisted patients and liver transplants performed in the United States between January 1st, 2015 and December 31st, 2020.
Table \ref{table:chracteristic} summarizes the demographic and clinical characteristics of patients and donors. During the study horizon, 68,557 patients registered for liver transplants at 140 transplant centers. Only 60,243 liver donations were made during this period, thus only about  87.9\% of the demand was met. 

\begin{table}[!htb]
\scriptsize
\centering
\caption{Demographic and clinical profile of study patients and donors from 2015 to 2020. Percentages are reported by ignoring missing values.}
\label{table:chracteristic}
\begin{tabular}{lcllllc}
\multicolumn{2}{c}{(a) Patients}                                                              &  &  &                       & \multicolumn{2}{c}{(b) Donors}                                          \\ \cline{1-2} \cline{6-7} 
\multicolumn{1}{|l}{Characteristic/Variable}  & \multicolumn{1}{l|}{Number of patients (\%)} &  &  & \multicolumn{1}{l|}{} & Characteristic/Variable & \multicolumn{1}{l|}{Number of donors (\%)} \\ \cline{1-2} \cline{6-7} 
\multicolumn{1}{|l}{Total patients}           & \multicolumn{1}{c|}{68,557}                  &  &  & \multicolumn{1}{l|}{} & Total donors            & \multicolumn{1}{c|}{60,243}                  \\
\multicolumn{1}{|l}{Gender}                   & \multicolumn{1}{c|}{}                        & \hspace{20pt} &  & \multicolumn{1}{l|}{} & Gender                  & \multicolumn{1}{c|}{}                        \\
\multicolumn{1}{|l}{\qquad Female}                   & \multicolumn{1}{c|}{25,479 (37.2\%)}         &  &  & \multicolumn{1}{l|}{} & \qquad Female                  & \multicolumn{1}{c|}{23,689 (39.3\%)}         \\
\multicolumn{1}{|l}{\qquad Male}                     & \multicolumn{1}{c|}{43,078 (62.8\%)}         &  &  & \multicolumn{1}{l|}{} & \qquad Male                    & \multicolumn{1}{c|}{36,554 (60.7\%)}         \\
\multicolumn{1}{|l}{Blood type}               & \multicolumn{1}{c|}{}                        &  &  & \multicolumn{1}{l|}{} & Blood type              & \multicolumn{1}{c|}{}                        \\
\multicolumn{1}{|l}{\qquad A}                        & \multicolumn{1}{c|}{25,612 (37.4\%)}         &  &  & \multicolumn{1}{l|}{} & \qquad A                       & \multicolumn{1}{c|}{22,274 (37.0\%)}         \\
\multicolumn{1}{|l}{\qquad B}                        & \multicolumn{1}{c|}{8,400(12.3\%)}           &  &  & \multicolumn{1}{l|}{} & \qquad B                       & \multicolumn{1}{c|}{6,908 (11.5\%)}          \\
\multicolumn{1}{|l}{\qquad AB}                       & \multicolumn{1}{c|}{2,725 (4.0\%)}           &  &  & \multicolumn{1}{l|}{} & \qquad AB                      & \multicolumn{1}{c|}{1,740 (2.9\%)}           \\
\multicolumn{1}{|l}{\qquad O}                        & \multicolumn{1}{c|}{31,820 (46.4\%)}         &  &  & \multicolumn{1}{l|}{} & \qquad O                       & \multicolumn{1}{c|}{29,321 (48.7\%)}         \\ \cline{6-7} 
\multicolumn{1}{|l}{MELD score}               & \multicolumn{1}{c|}{}                        &  &  &                       &                         & \multicolumn{1}{l}{}                         \\
\multicolumn{1}{|l}{$\qquad 6\sim20$}            & \multicolumn{1}{c|}{17,326 (42.2\%)}          &  &  &                       &                         & \multicolumn{1}{l}{}                         \\
\multicolumn{1}{|l}{\qquad $21\sim30$}               & \multicolumn{1}{c|}{12,186 (29.7\%)}         &  &  &                       &                         & \multicolumn{1}{l}{}                         \\
\multicolumn{1}{|l}{\qquad $31\sim 40$}                      & \multicolumn{1}{c|}{11,560 (28.1\%)}         &  &  &                       &                         & \multicolumn{1}{l}{}                         \\
\multicolumn{1}{|l}{Total transplant centers} & \multicolumn{1}{c|}{140}                     &  &  &                       &                         & \multicolumn{1}{l}{}                         \\ \cline{1-2}
\end{tabular}
\end{table}

Each of the 140 transplant centers is further divided into four patient waitlists based on blood types (A, B, AB, and O). As discussed before, each waitlist is modeled as a GI/MI/1+GI$_S$ queue as in Figure \ref{fig:single queue}. For each waitlist, the arrival intensity represents the annual waitlist additions of patients. A mixture of two exponential distributions is employed to model the patient inter-arrival times: the exponential distribution across all transplant centers has an average Kolmogorov-Smirnov distance of 0.14, which reduces to 0.07 when using a mixture of two exponential distributions. For further details and an illustrative example, see  Appendix \ref{sec:non-exp of IAT}.
Figure \ref{fig:addition dist} provides the histogram of waitlist arrival intensities for different blood types. 
Across all waitlists, the average arrival intensity is 23.6 patients per year. 

\begin{figure}[!htb]  
\centering 
{\includegraphics[width=1\textwidth, height = 80pt]{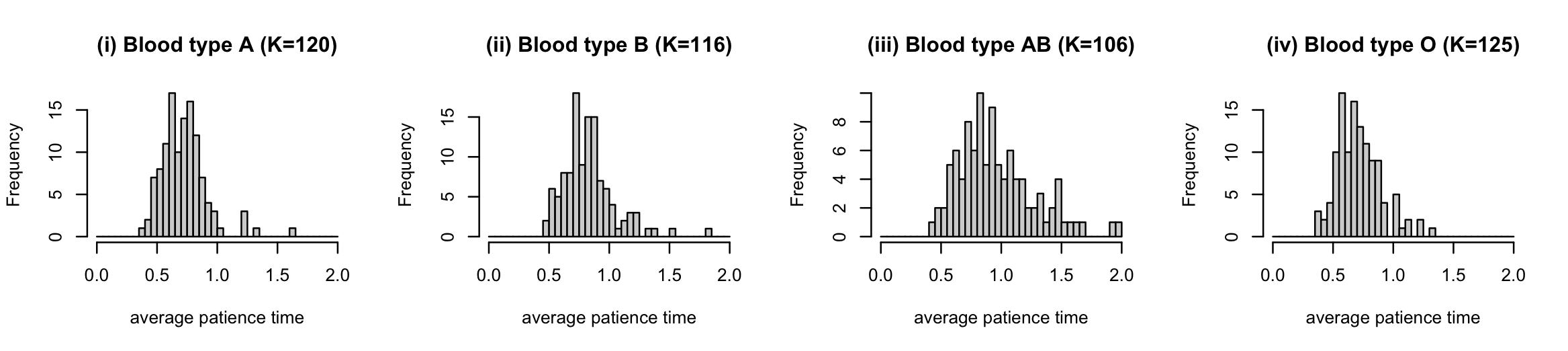}}
\caption{Histogram of patient waitlists' average patience times (years) by blood types. 
($K=$ number of waitlists.)
}
\label{fig:patience time dist}
\end{figure}

Patients' patience time is defined as the survival time without receiving a transplant from the date of registration.
To estimate the patience time distribution for a particular waitlist, the Model for End-Stage Liver Disease (MELD) score distribution of this waitlist and the three-month mortality rates associated with MELD scores, as described in \cite{myers2011gender}, are considered.
We use a mixture of three
exponential distributions based on three patient sub-classes within a waitlist identified using their MELD scores ($6\sim 20$, $21\sim 30$, $31\sim 40$).
The mixture distribution is further truncated to a maximum of 25 years to satisfy Assumption \ref{A:input distributions 1}.
The histograms of waitlists' average patience times for different blood types are given in Figure \ref{fig:patience time dist}. 
The average patience time of the entire patient population is estimated to be 0.7 years.

\subsection{Evaluating Steady-State Performance}
\label{sec:exp evaluation}
We first compare the accuracy of finite approximation $vs.$ fluid/diffusion approximation in evaluating the patient waitlists' long-term performance. 
For each waitlist, we compute the (average) offered sojourn time and the abandonment probability.  
In these performance evaluation experiments, the inter-arrival, service, and patience time distributions are estimated as in Subsection \ref{sec:exp data and parameter}, with the service rate being defined as the product of the arrival intensity and 0.879, which is the ratio of donors over patients in our data.
To implement finite approximation, a 4,097-state Markov chain is constructed for each waitlist. The time required to evaluate a single waitlist ranges between $75\sim90$ seconds. 
To implement fluid approximation, we apply the results of \cite{whitt2006fluid}. 
We note \cite{whitt2006fluid} assumes that the customer does not abandon after service begins. 
While this assumption deviates from the transplant scenario in the current work, it is considered inconsequential in fluid approximation since the method relies on high arrival and service rates, which results in negligible occurrence of abandonment during service. 
To implement diffusion approximation, we apply the second-order method in \cite{braverman2022high}.  We note that the evaluation of the transition kernel in \eqref{eqn:def kernel} involves two nested integrals, and the subsequent steps in this method require an additional two nested integrals.
Therefore, the computation of diffusion approximation involves four nested integrals. 
The average time required to evaluate a single waitlist instance using this method is  about 90 seconds. 
A steady-state simulation with a sample size of $10^7$ is employed as a benchmark to compare the results obtained from the above methods.

\begin{figure}[!htb]  
\centering 
{\includegraphics[width=.9\textwidth, height =180pt]{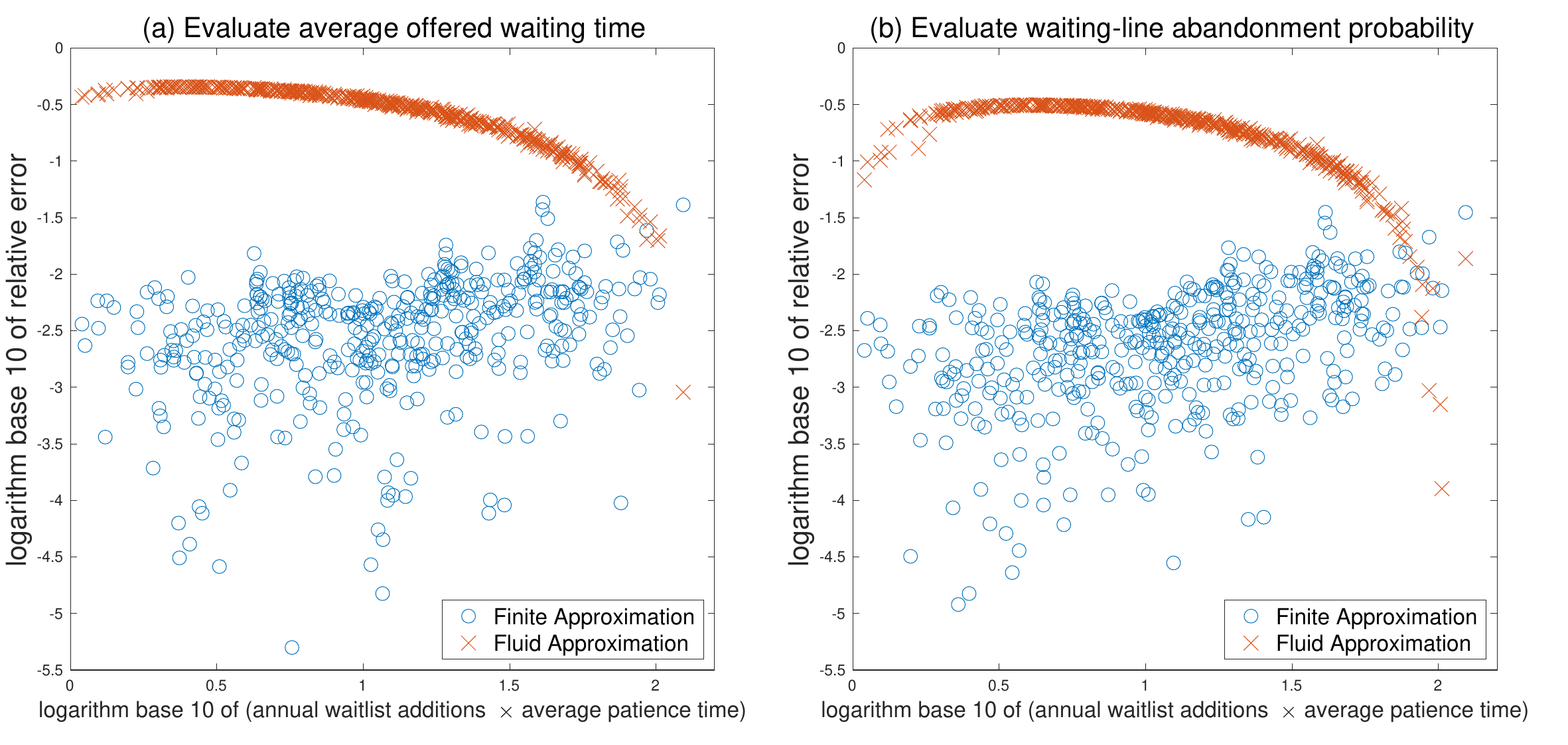}}
\caption{Relative errors of finite approximation and fluid approximation in evaluating average offered sojourn time and abandonment probability for each waitlist. (Number of waitlists = 467.  Each waitlist is evaluated twice respectively using finite approximation and fluid approximation and yields two points in each of plots (a) and (b).)
}
\label{fig:error dist}
\end{figure}

\begin{figure}[!htb]  
\centering 
{\includegraphics[width=.9\textwidth, height =180pt]{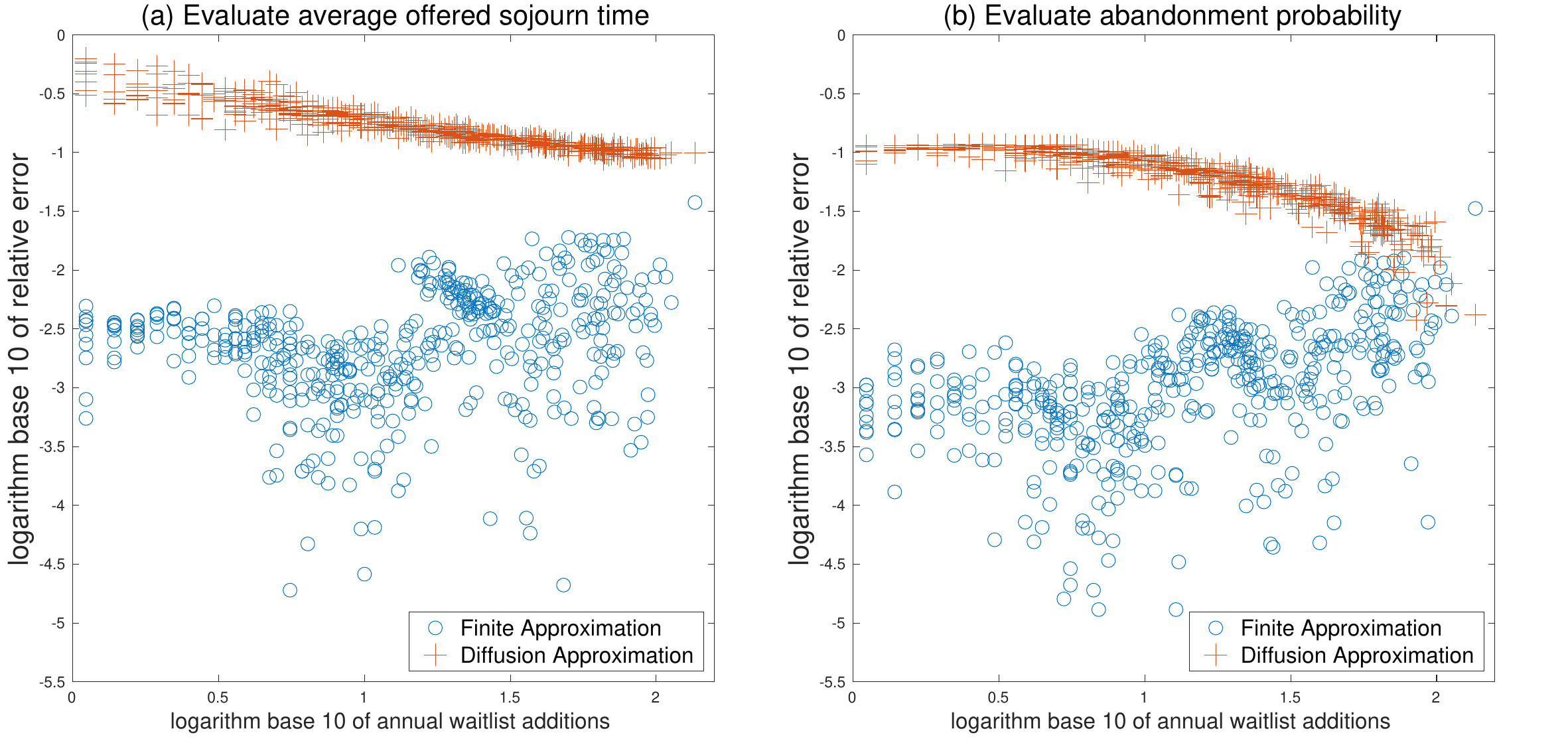}}
\caption{Relative errors of finite approximation and diffusion approximation in evaluating average offered sojourn time and abandonment probability for each waitlist. (Number of waitlists = 467.  Each waitlist is evaluated twice respectively using finite approximation and diffusion approximation and yields two points in each of plots (a) and (b).)
}
\label{fig:error dist da}
\end{figure}

The resulting relative error distributions are presented in Figures \ref{fig:error dist}--\ref{fig:error dist da} and Table \ref{table:error all table}. 
Out of the total of 467 waitlists analyzed, finite approximation has a relative error of less than 1\% for over 92\% of the waitlists.
In contrast, fluid approximation produces a relative error greater than 25\% for more than 93\% of the waitlists. Diffusion approximation produces a relative error greater than 10\% for over 90\% ($resp.$ 20\%) of the waitlists when evaluating offered sojourn time ($resp.$ abandonment probability). 
These findings are consistent with our observation in Figure \ref{fig:addition dist} that most waitlists fail to meet the large market (arrival intensity) assumption required by asymptotic analysis. We also recall that the accuracy of the finite approximation approach can be improved, as desired, by increasing the number of states in the Markov chain approximation.


\begin{table}[!htb]
\tiny
\centering
\caption{The relative error distributions of finite approximation, fluid approximation, and diffusion approximation in evaluating average offered sojourn time and abandonment probability for all waitlists.}
\label{table:error all table}
\begin{tabular}{|l|ccccc|ccccc|}
\hline
\multirow{2}{*}{\begin{tabular}[c]{@{}l@{}}Performance evaluation\\ relative error\end{tabular}} & \multicolumn{5}{c|}{(a) Evaluating average offered sojourn time}                                                                                                                                                                                                                       & \multicolumn{5}{c|}{(b) Evaluating abandonment probability}                                                                                                                                                                                                                \\
& \textless{}1\%                                       & 1\%$\sim$5\%                                       & 5\%$\sim$10\%                                        & 10\%$\sim$25\%                                         & \textgreater{}25\%                                     & \textless{}1\%                                       & 1\%$\sim$5\%                                       & 5\%$\sim$10\%                                         & 10\%$\sim$25\%                                         & \textgreater{}25\%                                     \\ \hline
\begin{tabular}[c]{@{}l@{}}Finite approximation,\\ waitlist number (\%)\end{tabular}             & \begin{tabular}[c]{@{}c@{}}434\\ (92.9\%)\end{tabular} & \begin{tabular}[c]{@{}c@{}}33\\ (7.1\%)\end{tabular} & \begin{tabular}[c]{@{}c@{}}0\\ (0)\end{tabular}  & \begin{tabular}[c]{@{}c@{}}0\\ (0)\end{tabular}        & \begin{tabular}[c]{@{}c@{}}0\\ (0)\end{tabular}        & \begin{tabular}[c]{@{}c@{}}458\\ (98.1\%)\end{tabular} & \begin{tabular}[c]{@{}c@{}}9\\ (1.9\%)\end{tabular}  & \begin{tabular}[c]{@{}c@{}}0\\ (0)\end{tabular}       & \begin{tabular}[c]{@{}c@{}}0\\ (0)\end{tabular}        & \begin{tabular}[c]{@{}c@{}}0\\ (0)\end{tabular}        \\ \hline
\begin{tabular}[c]{@{}l@{}}Fluid approximation,\\ waitlist number (\%)\end{tabular}              & \begin{tabular}[c]{@{}c@{}}0\\ (0)\end{tabular}    & \begin{tabular}[c]{@{}c@{}}0\\ (0)\end{tabular} & \begin{tabular}[c]{@{}c@{}}0\\ (0)\end{tabular} & \begin{tabular}[c]{@{}c@{}}8\\ (1.7\%)\end{tabular} & \begin{tabular}[c]{@{}c@{}}459\\ (98.3\%)\end{tabular} & \begin{tabular}[c]{@{}c@{}}0\\ (0)\end{tabular}   & \begin{tabular}[c]{@{}c@{}}0\\ (0)\end{tabular} & \begin{tabular}[c]{@{}c@{}}0\\ (0)\end{tabular} & \begin{tabular}[c]{@{}c@{}}30\\ (6.4\%)\end{tabular} & \begin{tabular}[c]{@{}c@{}}437\\ (93.6\%)\end{tabular} \\ \hline
\begin{tabular}[c]{@{}l@{}}Diffusion approximation,\\ waitlist number (\%)\end{tabular}              & \begin{tabular}[c]{@{}c@{}}0\\ (0)\end{tabular}    & \begin{tabular}[c]{@{}c@{}}0\\ (0)\end{tabular} & \begin{tabular}[c]{@{}c@{}}26\\ (5.6\%)\end{tabular} & \begin{tabular}[c]{@{}c@{}}355\\ (76.0\%)\end{tabular} & \begin{tabular}[c]{@{}c@{}}86\\ (18.4\%)\end{tabular} & \begin{tabular}[c]{@{}c@{}}6\\ (1.3\%)\end{tabular}   & \begin{tabular}[c]{@{}c@{}}146\\ (31.3\%)\end{tabular} & \begin{tabular}[c]{@{}c@{}}219\\ (46.9\%)\end{tabular} & \begin{tabular}[c]{@{}c@{}}96\\ (20.6\%)\end{tabular} & \begin{tabular}[c]{@{}c@{}}0\\ (0)\end{tabular} \\ \hline
\end{tabular}
\end{table}

\subsection{Optimizing Liver Allocation}
\label{sec:exp optimization}
We now discuss results from the two optimization models introduced in Section~\ref{sec:opt problem formulation}. Recall that these models maximize allocation equity while enforcing efficiency constraints. 
In model \eqref{model:ost} both efficiency and equity are measured by the average offered sojourn time, while in model \eqref{model:ab}, efficiency and equity are measured by the abandonment probabilities. 
Recall that in both models, the decision variables are liver allocation rates (service rates) for each waitlist, which measure the number of livers distributed to that waitlist per year. 
We assume that the sum of liver allocation rates across all waitlists should be no more than the number of donors per year; thus, we set constant $\mu_{\text{total}}=\pmb{1}^\text{T}\pmb{\lambda}\cdot 0.879$.
We also require that, for each waitlist, the allocated donor livers  should align with 65\% to 95\% of the demand; thus, we set constants $(\vartheta_L, \vartheta_U)=(0.65,0.95)$.
We respectively apply finite approximation (Algorithm \ref{alg1}) and fluid approximation (by substituting fluid approximation for finite approximation in Algorithm \ref{alg1}) to determine the optimal liver allocation schemes. We note that models \eqref{model:ost} and \eqref{model:ab} satisfy all assumptions made in this paper (Assumptions \ref{A:input distributions 1}--\ref{A:input distributions 3}). 
To assess the accuracy of the solutions obtained through finite and fluid approximations, we conduct a steady-state simulation with a sample size of $10^6$.
It takes finite approximation 11.6 hours to solve all optimization problems (a total of 64 PWL programs for four blood types, two models \eqref{model:ost} and \eqref{model:ab}, and eight different values of efficiency enforcement parameter $\varsigma$). 
The number of PWL approximation knots is seven.

\begin{figure}[!htb]  
\centering 
{\includegraphics[width=1.0\textwidth, height =280pt]{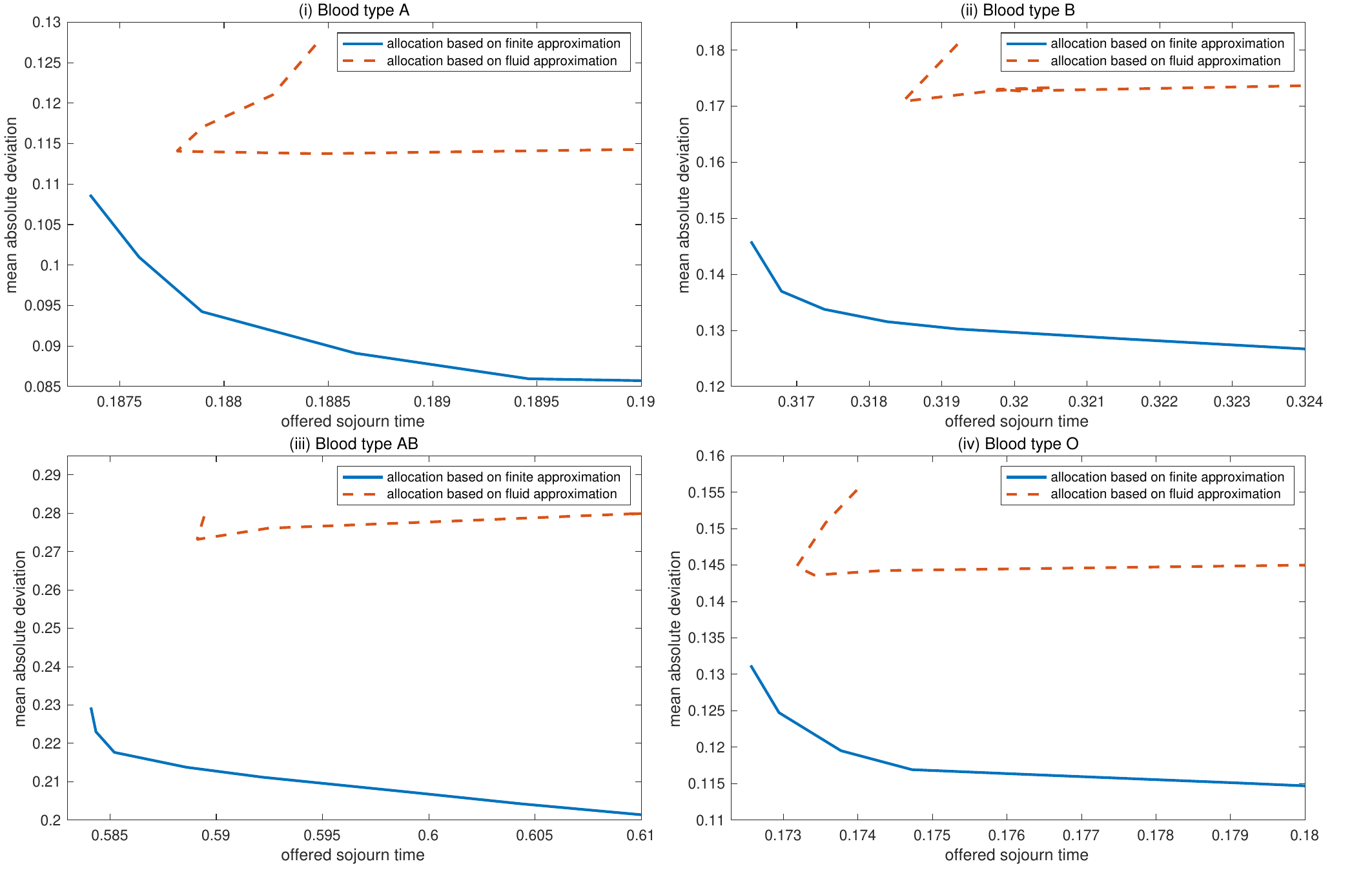}}
\caption{Performance of finite approximation $vs.$ fluid approximation in improving allocation equity (measured by mean absolute deviation across waitlists) of expected offered sojourn time under constraints of efficiency (measured by weighted average across waitlists)  in model \eqref{model:ost}.}
\label{fig:pareto ost}
\end{figure}

\begin{figure}[!htb]  
\centering 
{\includegraphics[width=1.0\textwidth, height =280pt]{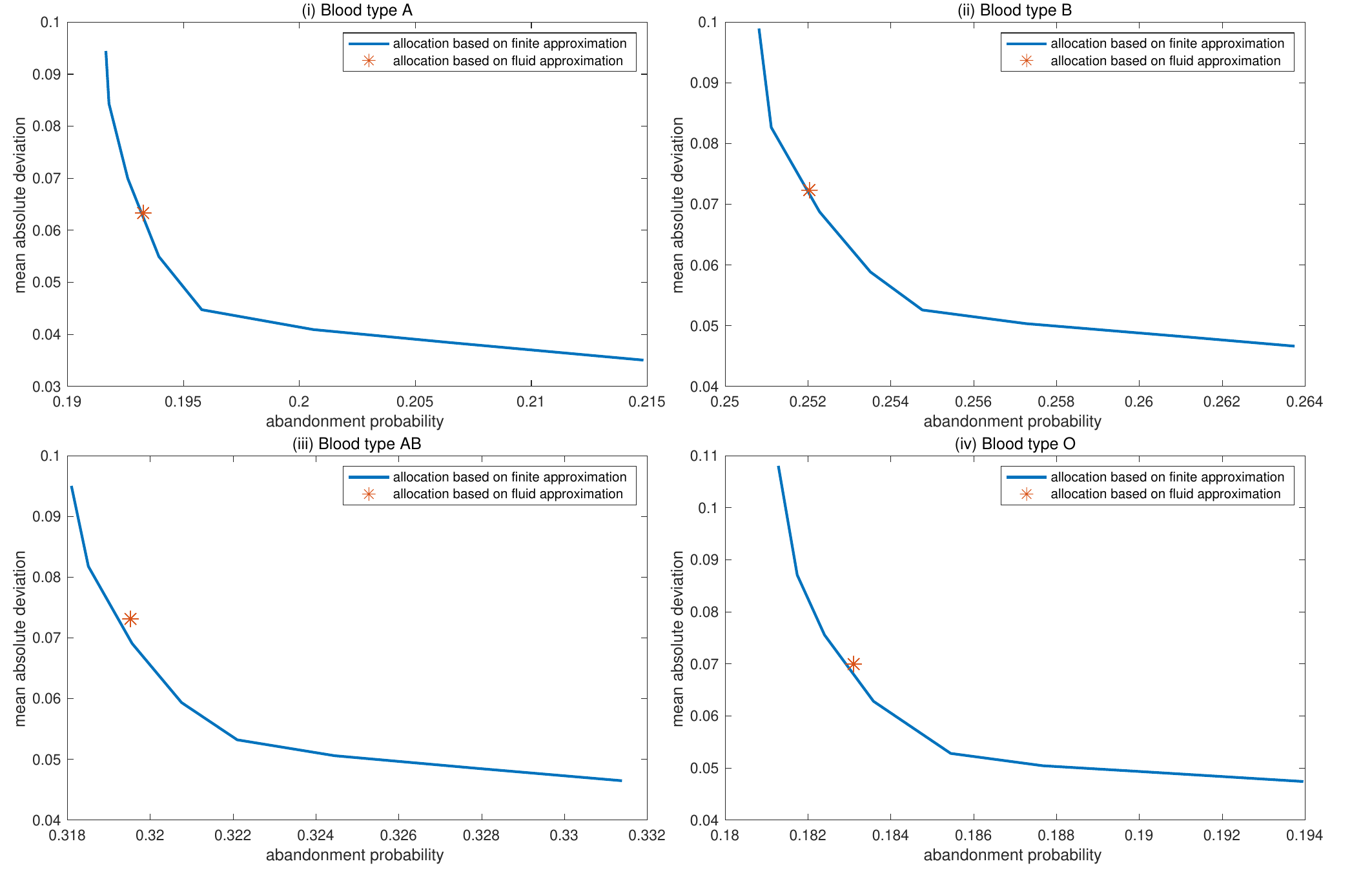}}
\caption{Performance of finite approximation $vs.$ fluid approximation in improving allocation equity (measured by mean absolute deviation across waitlists) of abandonment probability under constraints of efficiency (measured by weighted average across waitlists) in model \eqref{model:ab}. 
Note that, based on fluid approximation, the most efficient and meantime most equitable allocation scheme is to allocate in proportion to demands. 
As a result, fluid approximation has a unique solution plotted in each of (i)--(iv). }
\label{fig:pareto ab}
\end{figure}

\vspace{.1in}
\noindent{\bf Improved Resource Utilization.}\quad 
As shown in Figures \ref{fig:pareto ost}-\ref{fig:pareto ab}, solutions produced by finite approximation achieve greater equity in liver allocation as the efficiency constraint is relaxed. 
In contrast, solutions produced by fluid approximation exhibit suboptimal resource utilization and cannot further improve allocation equity after reaching a certain threshold.
In Figure \ref{fig:pareto ost}, for each blood type group, finite approximation is able to achieve the maximal efficiency fluid approximation can reach and meantime further improve equity by at least 14\% less deviation than the minimal deviation fluid approximation can reach.  
(iii) of Figure \ref{fig:pareto ost} highlights that fluid approximation can result in lower efficiency and inferior equity simultaneously. 
This can be attributed to the relatively small scale of the waitlists, as shown in Figure \ref{fig:addition dist}, and the consequential substantial evaluation errors of fluid approximation.
Figure \ref{fig:pareto ab} illustrates that fluid approximation lacks flexibility in balancing efficiency with equity when using abandonment probabilities as performance measures. 
These findings support the assertion made by \cite{zenios2000dynamic} that ``a fluid model is too crude to be used as a reliable performance analysis tool", and highlight the limitations of managerial insights derived from a fluid model.

\vspace{.1in}
\noindent{\bf Implications of transplant center scale and patient mortality risks.}\quad
We also explored how the optimal liver allocation rates vary as we relax efficiency constraints and enhance equity. 
Our analysis based on finite approximation reveals that a more equitable scheme is linked to a proportionally higher allocation of livers to waitlists with smaller arrival intensities or greater mortality risks ($i.e.$, average patience time).
To demonstrate our findings, we divide the waitlists into nine clusters based on small/medium/large arrival intensity and low/medium/high mortality risks.
Here ``small" and ``low" denote the first 33\% in quantile, while ``large" and ``high" denote the last 33\%.
All allocation rates are expressed as the percentage of arrival intensities (demand). 
The average liver allocation rates for each cluster are presented in Tables \ref{table:allocation rates ost}-\ref{table:allocation rates ab}.
We compare the most efficient and the most equitable allocation schemes and find that the cluster with large arrival intensity and/or low mortality risk ($resp.$ with small arrival intensity and/or high mortality risk) has a relatively smaller ($resp.$ larger) allocation rate in an equitable scheme. 

\begin{table}[!htb]
\caption{Cluster-Average liver allocation rates when optimizing offered sojourn time in model \eqref{model:ost}. (Liver allocation rates are expressed as a percentage of patient arrival intensity.)}
\label{table:allocation rates ost}
\centering
\tiny
\begin{tabular}{|l|ccc|ccc|ccc|}
\hline
\multirow{2}{*}{\small \textit{Blood type A}} & \multicolumn{3}{c|}{\bf \scriptsize Most efficient}                                                                                                                                        & \multicolumn{3}{c|}{\bf \scriptsize Most equitable}                                                                                                                                        & \multicolumn{3}{c|}{\bf \scriptsize Difference}                                                                                                                                            \\ \cline{3-3} \cline{6-6} \cline{9-9}
                              & \begin{tabular}[c]{@{}c@{}}Low \\ death risk\end{tabular} & \begin{tabular}[c]{@{}c@{}}Medium \\ death risk\end{tabular} & \begin{tabular}[c]{@{}c@{}}High \\ death risk\end{tabular} & \begin{tabular}[c]{@{}c@{}}Low \\ death risk\end{tabular} & \begin{tabular}[c]{@{}c@{}}Medium \\ death risk\end{tabular} & \begin{tabular}[c]{@{}c@{}}High \\ death risk\end{tabular} & \begin{tabular}[c]{@{}c@{}}Low \\ death risk\end{tabular} & \begin{tabular}[c]{@{}c@{}}Medium \\ death risk\end{tabular} & \begin{tabular}[c]{@{}c@{}}High \\ death risk\end{tabular} \\ \hline
Small arrival intensity       & 94.6\%                                                    & 90.6\%                                                       & 86.4\%                                                     & 95.0\%                                                    & 95.0\%                                                       & 94.7\%                                                     & 0.4\%                                                     & 4.4\%                                                        & 8.4\%                                                      \\
Medium arrival intensity      & 92.2\%                                                    & 86.4\%                                                       & 79.7\%                                                     & 92.0\%                                                    & 88.4\%                                                       & 84.3\%                                                     & -0.3\%                                                    & 2.0\%                                                        & 4.5\%                                                      \\
High arrival intensity         & 92.5\%                                                    & 87.7\%                                                       & 80.8\%                                                     & 88.4\%                                                    & 83.2\%                                                       & 79.5\%                                                     & -4.1\%                                                    & -4.6\%                                                       & -1.3\%                                                     \\ \hline \hline
\multirow{2}{*}{\small \textit{Blood type B}} & \multicolumn{3}{c|}{\bf \scriptsize Most efficient}                                                                                                                                        & \multicolumn{3}{c|}{\bf \scriptsize Most equitable}                                                                                                                                        & \multicolumn{3}{c|}{\bf \scriptsize Difference}                                                                                                                                            \\ \cline{3-3} \cline{6-6} \cline{9-9}
                              & \begin{tabular}[c]{@{}c@{}}Low \\ death risk\end{tabular} & \begin{tabular}[c]{@{}c@{}}Medium \\ death risk\end{tabular} & \begin{tabular}[c]{@{}c@{}}High \\ death risk\end{tabular} & \begin{tabular}[c]{@{}c@{}}Low \\ death risk\end{tabular} & \begin{tabular}[c]{@{}c@{}}Medium \\ death risk\end{tabular} & \begin{tabular}[c]{@{}c@{}}High \\ death risk\end{tabular} & \begin{tabular}[c]{@{}c@{}}Low \\ death risk\end{tabular} & \begin{tabular}[c]{@{}c@{}}Medium \\ death risk\end{tabular} & \begin{tabular}[c]{@{}c@{}}High \\ death risk\end{tabular} \\ \hline
Small arrival intensity       & 95.0\%                                                    & 93.8\%                                                       & 92.5\%                                                     & 95.0\%                                                    & 95.0\%                                                       & 95.0\%                                                     & 0.0\%                                                     & 1.3\%                                                        & 2.5\%                                                      \\
Medium arrival intensity      & 92.8\%                                                    & 90.1\%                                                       & 83.2\%                                                     & 91.9\%                                                    & 90.7\%                                                       & 89.3\%                                                     & -0.8\%                                                    & 0.5\%                                                        & 6.1\%                                                      \\
High arrival intensity         & 92.5\%                                                    & 86.3\%                                                       & 78.8\%                                                     & 88.1\%                                                    & 83.2\%                                                       & 81.7\%                                                     & -4.4\%                                                    & -3.0\%                                                       & 3.0\%                                                      \\ \hline \hline
\multirow{2}{*}{\small \textit{Blood type AB} } & \multicolumn{3}{c|}{\bf \scriptsize Most efficient}                                                                                                                                        & \multicolumn{3}{c|}{\bf \scriptsize Most equitable}                                                                                                                                        & \multicolumn{3}{c|}{\bf \scriptsize Difference}                                                                                                                                            \\ \cline{3-3} \cline{6-6} \cline{9-9}
                              & \begin{tabular}[c]{@{}c@{}}Low \\ death risk\end{tabular} & \begin{tabular}[c]{@{}c@{}}Medium \\ death risk\end{tabular} & \begin{tabular}[c]{@{}c@{}}High \\ death risk\end{tabular} & \begin{tabular}[c]{@{}c@{}}Low \\ death risk\end{tabular} & \begin{tabular}[c]{@{}c@{}}Medium \\ death risk\end{tabular} & \begin{tabular}[c]{@{}c@{}}High \\ death risk\end{tabular} & \begin{tabular}[c]{@{}c@{}}Low \\ death risk\end{tabular} & \begin{tabular}[c]{@{}c@{}}Medium \\ death risk\end{tabular} & \begin{tabular}[c]{@{}c@{}}High \\ death risk\end{tabular} \\ \hline
Small arrival intensity       & 95.0\%                                                    & 95.0\%                                                       & 92.7\%                                                     & 95.0\%                                                    & 95.0\%                                                       & 94.9\%                                                     & 0.0\%                                                     & 0.0\%                                                        & 2.1\%                                                      \\
Medium arrival intensity      & 95.0\%                                                    & 91.3\%                                                       & 82.5\%                                                     & 93.9\%                                                    & 89.4\%                                                       & 81.8\%                                                     & -1.1\%                                                    & -1.8\%                                                       & -0.7\%                                                     \\
High arrival intensity         & 92.7\%                                                    & 87.5\%                                                       & 76.7\%                                                     & 89.6\%                                                    & 81.9\%                                                       & 78.0\%                                                     & -3.1\%                                                    & -5.6\%                                                       & 1.4\%                                                      \\ \hline \hline
\multirow{2}{*}{\small \textit{Blood type O}} & \multicolumn{3}{c|}{\bf \scriptsize Most efficient}                                                                                                                                        & \multicolumn{3}{c|}{\bf \scriptsize Most equitable}                                                                                                                                        & \multicolumn{3}{c|}{\bf \scriptsize Difference}                                                                                                                                            \\ \cline{3-3} \cline{6-6} \cline{9-9}
                              & \begin{tabular}[c]{@{}c@{}}Low \\ death risk\end{tabular} & \begin{tabular}[c]{@{}c@{}}Medium \\ death risk\end{tabular} & \begin{tabular}[c]{@{}c@{}}High \\ death risk\end{tabular} & \begin{tabular}[c]{@{}c@{}}Low \\ death risk\end{tabular} & \begin{tabular}[c]{@{}c@{}}Medium \\ death risk\end{tabular} & \begin{tabular}[c]{@{}c@{}}High \\ death risk\end{tabular} & \begin{tabular}[c]{@{}c@{}}Low \\ death risk\end{tabular} & \begin{tabular}[c]{@{}c@{}}Medium \\ death risk\end{tabular} & \begin{tabular}[c]{@{}c@{}}High \\ death risk\end{tabular} \\ \hline
Small arrival intensity       & 94.7\%                                                    & 88.9\%                                                       & 85.8\%                                                     & 95.0\%                                                    & 95.0\%                                                       & 94.8\%                                                     & 0.3\%                                                     & 6.1\%                                                        & 9.0\%                                                      \\
Medium arrival intensity      & 92.0\%                                                    & 87.1\%                                                       & 80.7\%                                                     & 92.9\%                                                    & 89.1\%                                                       & 84.2\%                                                     & 0.9\%                                                     & 2.0\%                                                        & 3.5\%                                                      \\
High arrival intensity         & 93.6\%                                                    & 88.1\%                                                       & 79.4\%                                                     & 88.3\%                                                    & 83.8\%                                                       & 80.9\%                                                     & -5.3\%                                                    & -4.4\%                                                       & 1.5\%                                                      \\ \hline
\end{tabular}
\end{table}

\begin{table}[!htb]
\caption{Cluster-Average liver allocation rates when optimizing abandonment probability in model \eqref{model:ab}. (Liver allocation rates are expressed as a percentage of patient arrival intensity.)}
\label{table:allocation rates ab}
\tiny
\centering
\begin{tabular}{|l|ccc|ccc|ccc|}
\hline
\multirow{2}{*}{\small \textit{Blood type A}} & \multicolumn{3}{c|}{\bf \scriptsize Most efficient}                                                                                                                                        & \multicolumn{3}{c|}{\bf \scriptsize Most equitable}                                                                                                                                        & \multicolumn{3}{c|}{\bf \scriptsize Difference}                                                                                                                                            \\ \cline{3-3} \cline{6-6} \cline{9-9}
                              & \begin{tabular}[c]{@{}c@{}}Low \\ death risk\end{tabular} & \begin{tabular}[c]{@{}c@{}}Medium \\ death risk\end{tabular} & \begin{tabular}[c]{@{}c@{}}High \\ death risk\end{tabular} & \begin{tabular}[c]{@{}c@{}}Low \\ death risk\end{tabular} & \begin{tabular}[c]{@{}c@{}}Medium \\ death risk\end{tabular} & \begin{tabular}[c]{@{}c@{}}High \\ death risk\end{tabular} & \begin{tabular}[c]{@{}c@{}}Low \\ death risk\end{tabular} & \begin{tabular}[c]{@{}c@{}}Medium \\ death risk\end{tabular} & \begin{tabular}[c]{@{}c@{}}High \\ death risk\end{tabular} \\ \hline
Small arrival intensity       & 75.4\%                                                    & 78.8\%                                                       & 70.1\%                                                     & 94.0\%                                                    & 95.0\%                                                       & 95.0\%                                                     & 18.6\%                                                    & 16.3\%                                                       & 24.9\%                                                     \\
Medium arrival intensity      & 86.7\%                                                    & 85.2\%                                                       & 83.4\%                                                     & 84.1\%                                                    & 87.6\%                                                       & 90.7\%                                                     & -2.6\%                                                    & 2.3\%                                                        & 7.2\%                                                      \\
High arrival intensity        & 91.7\%                                                    & 90.0\%                                                       & 89.8\%                                                     & 82.1\%                                                    & 81.6\%                                                       & 82.4\%                                                     & -9.5\%                                                    & -8.4\%                                                       & -7.3\%                                                     \\ \hline  \hline
\multirow{2}{*}{\small \textit{Blood type B}} & \multicolumn{3}{c|}{\bf \scriptsize Most efficient}                                                                                                                                        & \multicolumn{3}{c|}{\bf \scriptsize Most equitable}                                                                                                                                        & \multicolumn{3}{c|}{\bf \scriptsize Difference}                                                                                                                                            \\ \cline{3-3} \cline{6-6} \cline{9-9}
                              & \begin{tabular}[c]{@{}c@{}}Low \\ death risk\end{tabular} & \begin{tabular}[c]{@{}c@{}}Medium \\ death risk\end{tabular} & \begin{tabular}[c]{@{}c@{}}High \\ death risk\end{tabular} & \begin{tabular}[c]{@{}c@{}}Low \\ death risk\end{tabular} & \begin{tabular}[c]{@{}c@{}}Medium \\ death risk\end{tabular} & \begin{tabular}[c]{@{}c@{}}High \\ death risk\end{tabular} & \begin{tabular}[c]{@{}c@{}}Low \\ death risk\end{tabular} & \begin{tabular}[c]{@{}c@{}}Medium \\ death risk\end{tabular} & \begin{tabular}[c]{@{}c@{}}High \\ death risk\end{tabular} \\ \hline
Small arrival intensity       & 75.4\%                                                    & 78.1\%                                                       & 70.8\%                                                     & 95.0\%                                                    & 95.0\%                                                       & 95.0\%                                                     & 19.6\%                                                    & 16.9\%                                                       & 24.2\%                                                     \\
Medium arrival intensity      & 87.8\%                                                    & 86.5\%                                                       & 82.3\%                                                     & 85.6\%                                                    & 89.6\%                                                       & 94.5\%                                                     & -2.2\%                                                    & 3.0\%                                                        & 12.2\%                                                     \\
High arrival intensity        & 90.9\%                                                    & 90.6\%                                                       & 90.0\%                                                     & 79.3\%                                                    & 84.0\%                                                       & 86.8\%                                                     & -11.7\%                                                   & -6.6\%                                                       & -3.2\%                                                     \\ \hline \hline
\multirow{2}{*}{\small \textit{Blood type AB}} & \multicolumn{3}{c|}{\bf \scriptsize Most efficient}                                                                                                                                        & \multicolumn{3}{c|}{\bf \scriptsize Most equitable}                                                                                                                                        & \multicolumn{3}{c|}{\bf \scriptsize Difference}                                                                                                                                            \\ \cline{3-3} \cline{6-6} \cline{9-9}
                              & \begin{tabular}[c]{@{}c@{}}Low \\ death risk\end{tabular} & \begin{tabular}[c]{@{}c@{}}Medium \\ death risk\end{tabular} & \begin{tabular}[c]{@{}c@{}}High \\ death risk\end{tabular} & \begin{tabular}[c]{@{}c@{}}Low \\ death risk\end{tabular} & \begin{tabular}[c]{@{}c@{}}Medium \\ death risk\end{tabular} & \begin{tabular}[c]{@{}c@{}}High \\ death risk\end{tabular} & \begin{tabular}[c]{@{}c@{}}Low \\ death risk\end{tabular} & \begin{tabular}[c]{@{}c@{}}Medium \\ death risk\end{tabular} & \begin{tabular}[c]{@{}c@{}}High \\ death risk\end{tabular} \\ \hline
Small arrival intensity       & 82.7\%                                                    & 72.3\%                                                       & 73.6\%                                                     & 93.4\%                                                    & 95.0\%                                                       & 95.0\%                                                     & 10.7\%                                                    & 22.7\%                                                       & 21.4\%                                                     \\
Medium arrival intensity      & 88.3\%                                                    & 83.8\%                                                       & 82.1\%                                                     & 82.4\%                                                    & 93.0\%                                                       & 94.1\%                                                     & -5.9\%                                                    & 9.3\%                                                        & 11.9\%                                                     \\
High arrival intensity        & 93.0\%                                                    & 90.3\%                                                       & 90.0\%                                                     & 78.3\%                                                    & 81.6\%                                                       & 86.9\%                                                     & -14.7\%                                                   & -8.6\%                                                       & -3.1\%                                                     \\ \hline \hline
\multirow{2}{*}{\small \textit{Blood type O}} & \multicolumn{3}{c|}{\bf \scriptsize Most efficient}                                                                                                                                        & \multicolumn{3}{c|}{\bf \scriptsize Most equitable}                                                                                                                                        & \multicolumn{3}{c|}{\bf \scriptsize Difference}                                                                                                                                            \\ \cline{3-3} \cline{6-6} \cline{9-9}
                              & \begin{tabular}[c]{@{}c@{}}Low \\ death risk\end{tabular} & \begin{tabular}[c]{@{}c@{}}Medium \\ death risk\end{tabular} & \begin{tabular}[c]{@{}c@{}}High \\ death risk\end{tabular} & \begin{tabular}[c]{@{}c@{}}Low \\ death risk\end{tabular} & \begin{tabular}[c]{@{}c@{}}Medium \\ death risk\end{tabular} & \begin{tabular}[c]{@{}c@{}}High \\ death risk\end{tabular} & \begin{tabular}[c]{@{}c@{}}Low \\ death risk\end{tabular} & \begin{tabular}[c]{@{}c@{}}Medium \\ death risk\end{tabular} & \begin{tabular}[c]{@{}c@{}}High \\ death risk\end{tabular} \\ \hline
Small arrival intensity       & 72.1\%                                                    & 75.0\%                                                       & 67.4\%                                                     & 94.9\%                                                    & 95.0\%                                                       & 95.0\%                                                     & 22.8\%                                                    & 20.0\%                                                       & 27.6\%                                                     \\
Medium arrival intensity      & 86.0\%                                                    & 83.4\%                                                       & 83.0\%                                                     & 87.0\%                                                    & 90.2\%                                                       & 91.4\%                                                     & 1.0\%                                                     & 6.7\%                                                        & 8.4\%                                                      \\
High arrival intensity        & 91.7\%                                                    & 90.5\%                                                       & 88.9\%                                                     & 84.2\%                                                    & 84.2\%                                                       & 84.6\%                                                     & -7.5\%                                                    & -6.3\%                                                       & -4.3\%                                                     \\ \hline
\end{tabular}
\end{table}

\begin{figure}[!htb]  
\centering 
\scriptsize 
\textbf{(i) Changes of allocation when improving equity of offered sojourn time for patients of blood type A.}
{\includegraphics[width=1\textwidth, height =160pt]{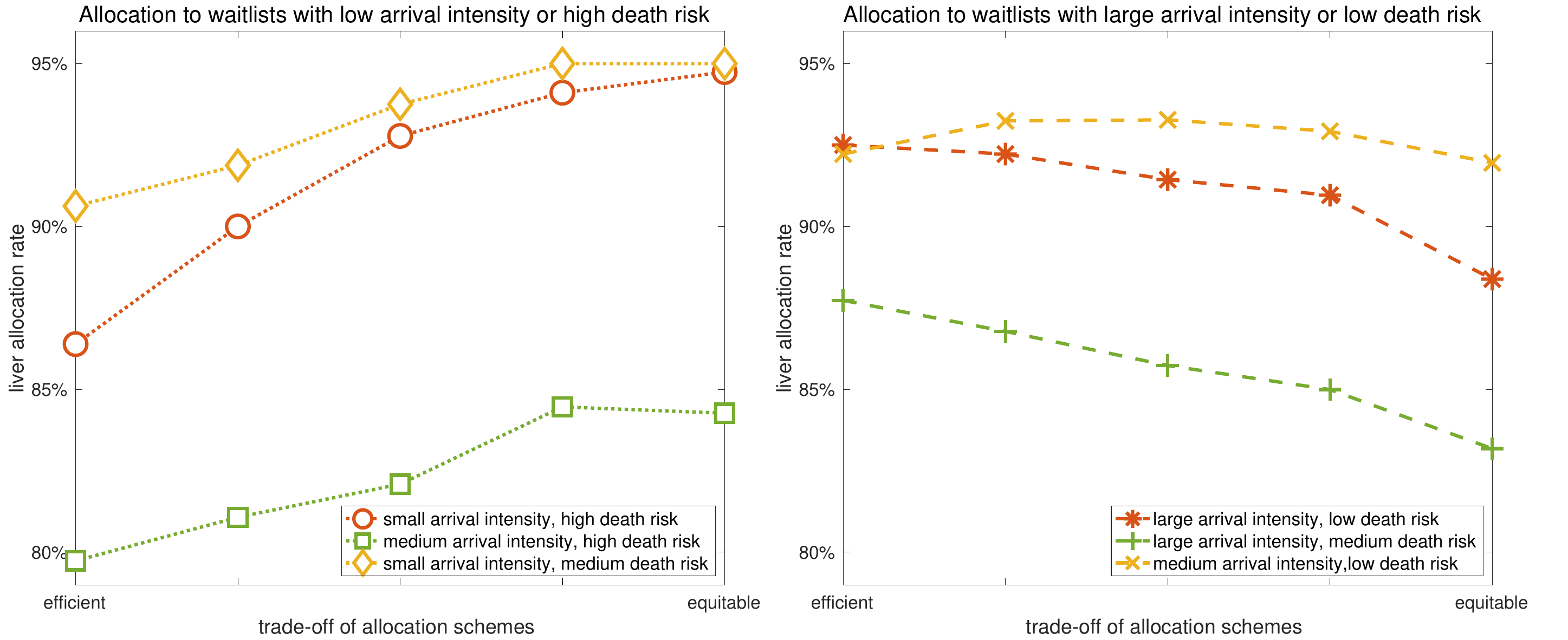}}
\textbf{(ii) Changes of allocation when improving equity of abandonment probability for patients of blood type B.}
{\includegraphics[width=1\textwidth, height =160pt]{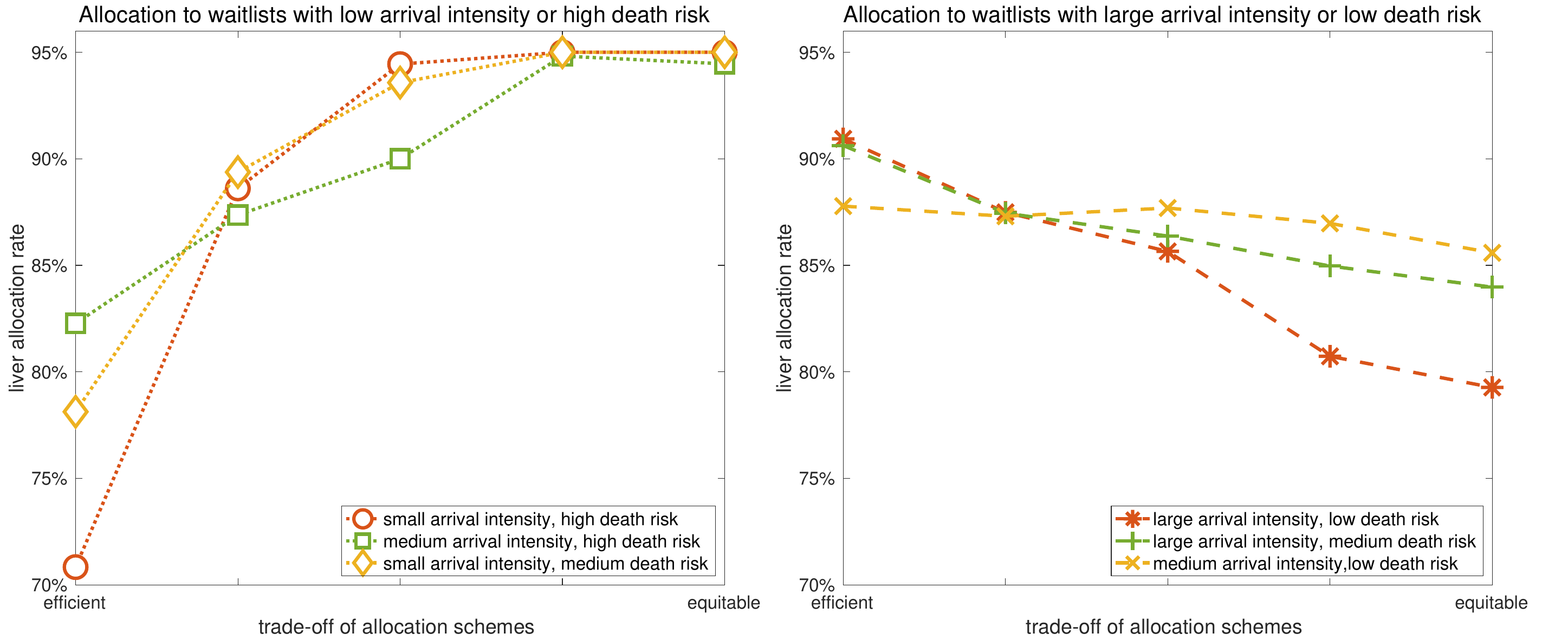}}
\caption{Waitlists with a small arrival intensity or a high death risk receive more livers in an equitable allocation scheme. (Liver allocation rates are expressed as a percentage of patient arrival intensity.) 
}
\label{fig:cluster trend}
\end{figure}

\begin{figure}[!htb]  
\centering 
\scriptsize 
\textbf{(i) Optimizing offered sojourn time for patients of blood type A.}
{\includegraphics[width=1\textwidth, height =140pt]{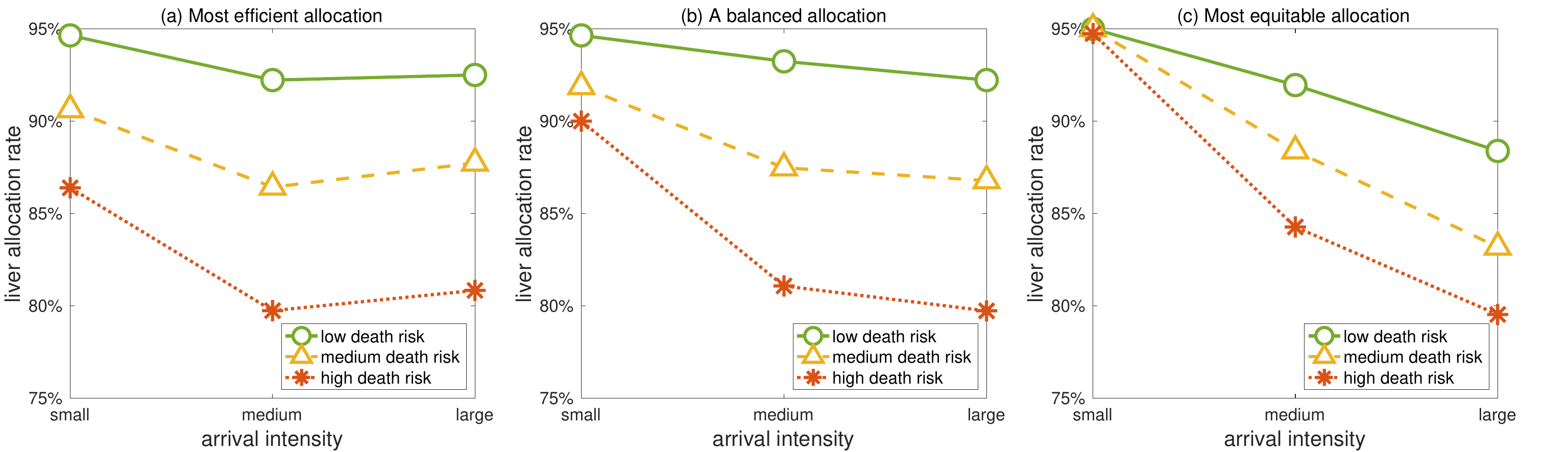}}
\textbf{(ii) Optimizing  abandonment probability for patients of blood type B.}
{\includegraphics[width=1\textwidth, height =140pt]{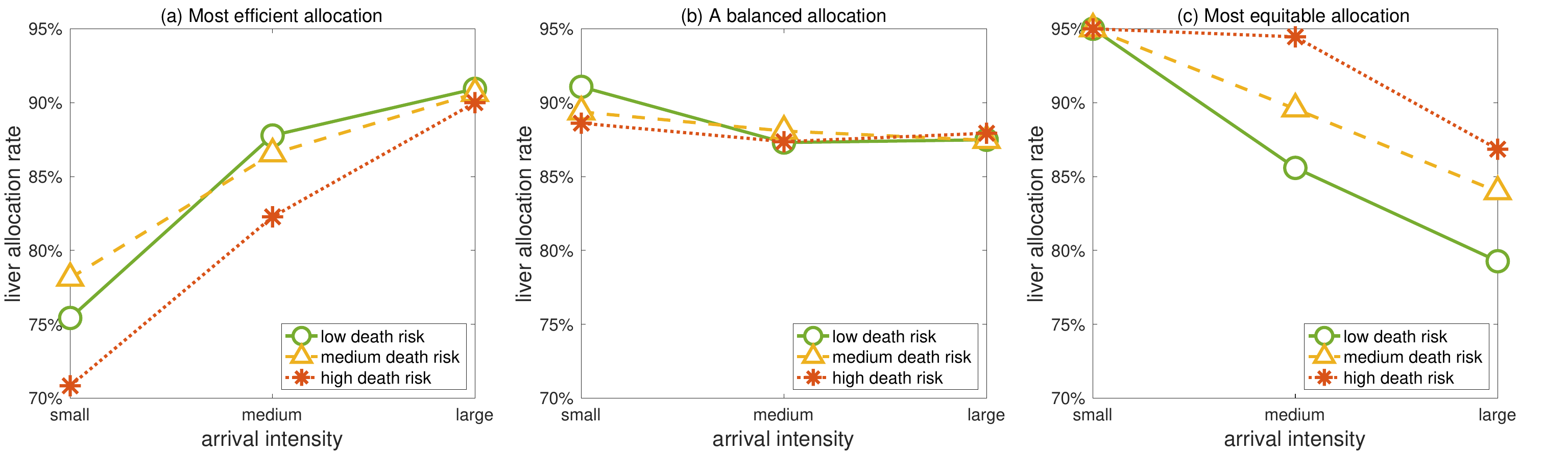}}
\caption{Cluster-Average liver allocation rates based on optimal solutions prescribed by finite approximation. (Liver allocation rates are expressed as a percentage of patient arrival intensity.)}
\label{fig:cluster plot}
\end{figure}

For a better understanding of the indications of arrival intensities and mortality risks, we further examine the allocation schemes for blood type A when optimizing offered sojourn time and the allocation schemes for blood type B when optimizing abandonment probability. 
Figure \ref{fig:cluster trend} shows that, for the three clusters with small arrival intensity or high mortality risk, the liver allocation rates increase by at least 4\% ($resp.$ 12\%) in an equitable allocation scheme in terms of offered sojourn time ($resp.$ abandonment probability). 
In contrast, the three clusters with large arrival intensity or low mortality risk are associated with lower allocation rates as equity improves.
The cluster-average liver allocation rates are plotted in Figure \ref{fig:cluster plot}.


\section{Concluding Remarks}\label{sec:conclusion}
This work investigates the equitable and efficient allocation of resources, modeled as service rates, to multiple parallel and any-scale queues. 
The study is motivated by the need to model liver allocation in the national transplant system of the United States, which experiences a scarcity of donor organs and includes a large number of small-scale waitlists. 
The waitlists are modeled as non-fluid queues with service abandonment. 
A finite approximation technique is developed to evaluate system performance as an alternative to fluid or diffusion approximations which rely on a large market assumption for accuracy. 
An algorithm is also presented that uses finite approximation to optimally allocate the arriving resource (donated organs) to each waitlist. 
Results from numerical experiments based on the OPTN database indicate that finite approximation provides more accurate solutions that significantly improve resource utilization compared to those prescribed by a fluid model. 
Notably, allocation equity (respectively, efficiency) is improved by appropriately increasing (respectively, decreasing) the proportion of organs allocated to waitlists with small scales or higher abandonment (death) risks.
Our findings confirm that a policy that prioritizes sicker patients is appropriate. Note that the current policy gives priority to patients based on MELD score and additional sickness statuses (e.g., Status 1A and  1B patients).  
However, to enhance allocation equity further, our results suggest allocating a proportionately greater number of organs to smaller transplant centers, and more vulnerable patients whose waitlist survival is not sufficiently captured in the MELD score should be considered in future policy revisions.

\bigskip




\bibliographystyle{informs2014} 
\vspace{.1in}
\def\bibfont{\scriptsize}
\bibliography{reference} 



\newpage
\begin{APPENDICES}
\section{Summary of Major Notation}\label{Append:Notation}
\smallskip
Table \ref{table:notation summary} summarizes the major mathematical notation used in the manuscript.
\scriptsize
\begin{longtable}{c|l}
\caption{Major notation and definitions}\\
\hline 
\label{table:notation summary}
\textbf{Notation} & \textbf{Definition}\\
\hline
\endfirsthead
\multicolumn{2}{c}%
{\tablename\ \thetable\ -- \textit{Continued from previous page}} \\
\hline
\textbf{Notation} & \textbf{Definition} \\
\hline
\endhead
\hline \multicolumn{2}{r}{\textit{Continued on next page}} \\
\endfoot
\hline
\endlastfoot
\\[-10pt]
${A_k}$ & Probability distribution function of inter-arrival times in the $k$-th waitlist.  \\ \hline \\[-10pt]
${a'_k,a_k''}$ & Bound of first- and second-order derivatives of $A_k$.  \\ \hline \\[-10pt]
${B}$ & Exponential distribution function of service times.  \\ \hline \\[-10pt]
$\{c_{k_i}^{(r)}\}_{i=1}^{J^{(r)}} $ & States of approximate finite-state Markov chain for the $k$-th waitlist.  \\ \hline \\[-10pt]
${G_k}$ & Probability distribution function of patience time in the $k$-th waitlist.  \\ \hline \\[-10pt]
${g_l}$ & Generic performance function.  \\ \hline \\[-10pt]
${J^{(r)}}$ & Number of states of the $r$-th approximate finite-state Markov chain. \\ \hline \\[-10pt]
$k,K,\mathscr{K}$ & Index, total number, and index set of patient waitlists.  \\ \hline \\[-10pt]
$l,L,\mathscr{L}$ & Index, total number, and index set of constraints.  \\ \hline \\[-10pt]
${\lambda_k}$ & Patient arrival intensity in the $k$-th waitlist.  \\ \hline \\[-10pt]
$(m)$ & Index of the piecewise linear approximation sequence. \\ \hline \\[-10pt]
${\mu_k},\mu_{\max},\mu_{\min}$ & Service (liver allocation) rate at the $k$-th waitlist and its upper and lower bounds. \\ \hline \\[-10pt]
${N^{(m)}}$ & Number of knots used in the $m$-th approximate piecewise linear  constraint.  \\ \hline \\[-10pt]
$\pmb{p},\pmb{w},\pmb{M},\pmb{d}$  & Coefficients and variables in the stochastic program for liver allocation.  \\ \hline \\[-10pt]
${\pi_k},\pi_k^{(r)}$ &  \begin{tabular}[l]{@{}l@{}}Stationary distribution of offered waiting time in the $k$-th waitlist and  its finite  \\approximation. \end{tabular} \\ \hline \\[-10pt]
${\phi}^{(r,m)}_{k,l}$ & Approximate piecewise linear constraint function for the $k$-th waitlist. \\ \hline \\[-10pt]
$\pmb{\Psi}_k$, $\pmb{\Psi}^{(r)}_k$ & \begin{tabular}[l]{@{}l@{}}Markov chain defined by offered waiting time of the $k$-th waitlist and its\\ finite approximation. \end{tabular}  \\ \hline \\[-10pt]
${\pmb{Q}_k^{(r)}}$, ${q_{k_{i,j}}^{(r)}}$ & Transition matrix of  $\pmb{\Psi}^{(r)}_k$ and its elemenets. \\ \hline \\[-10pt]
$(r)$ & Index of approximate finite-state Markov chain sequence.  \\ \hline \\[-10pt]
${s^n_k}$ & Service time required by patient $n$ in the $k$-th waitlist. \\ \hline \\[-10pt]
${\tau_k}$ & Transition kernel associated with $\pmb{\Psi}_k$. \\ \hline \\[-10pt]
${t^n_k}$ & Inter-arrival time of patient $n$ in the $k$-th waitlist. \\ \hline \\[-10pt]
$\pmb{v}^{(r)}$ & Stationary stochastic vector under transition matrix $\pmb{Q}^{(r)}$. \\ \hline \\[-10pt]
$V(\cdot)$ ($V_u(\cdot)$) & Total variation of a function ($w.r.t.$ variable $u$). \\ \hline \\[-10pt]
${\xi^n_k}$ & Offered waiting time of patient $n$ in the $k$-th waitlist. \\ \hline \\[-10pt]
${y^n_k},\bar{y}_k$ & Patience time of patient $n$ in the $k$-th waitlist and its upper bound. \\ \hline\\[-10pt]
$\zeta_l$ & Parameter of Assumption \ref{A:input distributions 3} on $g_l$.
\\ \hline 
\end{longtable}
\normalsize

\section{Hypothesis Tests for Justifying the Queueing Models} \label{sec:hypothesistest}
\subsection{Exponetiality of Service Times}\label{sec:exp of ST}
Using our OPTN dataset, we now discuss the service time distributions, $i.e.$, the donor liver arrival process is Poisson.  
We use exponential distributions to fit liver inter-arrival durations for the 50 largest donor hospitals and observe that, for both blood types O and A, each hospital has a $p$-value greater than 0.05 in a Kolmogorov–Smirnov test, leading us to accept the hypothesis that donor livers follow a Poisson process. 
Across all donor hospitals, the average Kolmogorov–Smirnov distance between the fitted exponential distribution and the empirical liver inter-arrival duration distribution is 0.07. As an illustration, Figure \ref{fig:exp test service time} shows the exponential distribution fits for livers in two instances of hospitals with blood types O and A, respectively, and annual liver donations of 25 and 18, respectively. The Kolmogorov–Smirnov test has a large $p$-value ($p>0.4$).

\begin{figure}[!htb]  
\centering 
{\includegraphics[width=1\textwidth, height =115pt]{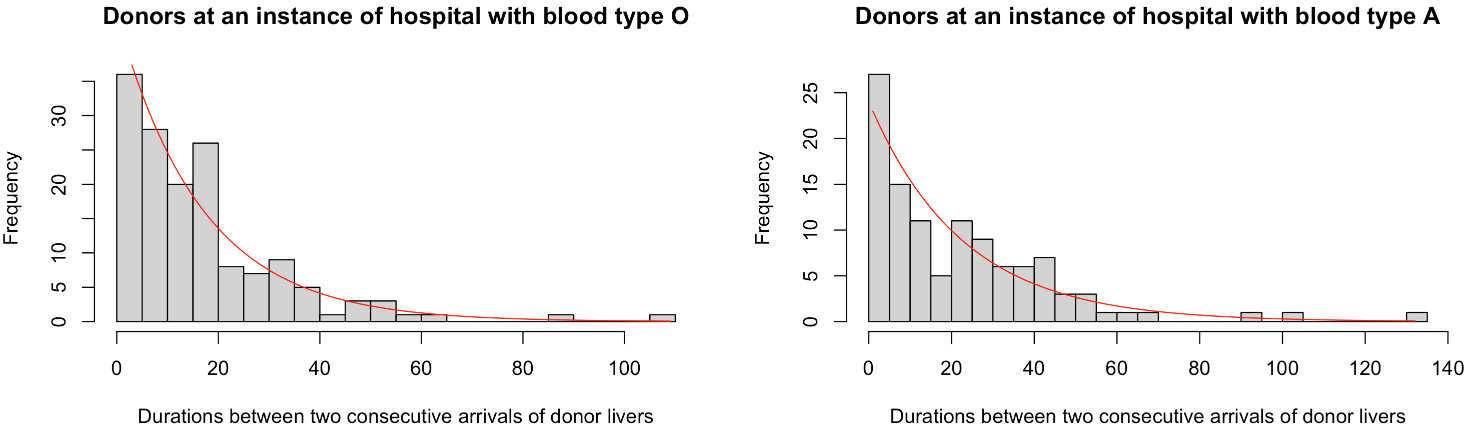}}
\vspace{-20pt}
\caption{Gap between the histogram of liver inter-arrival durations (days) and the fitted exponential distribution in two example donor hosptials. 
}
\label{fig:exp test service time}
\end{figure}

\begin{figure}[!htb]  
\centering 
{\includegraphics[width=1\textwidth, height =155pt]{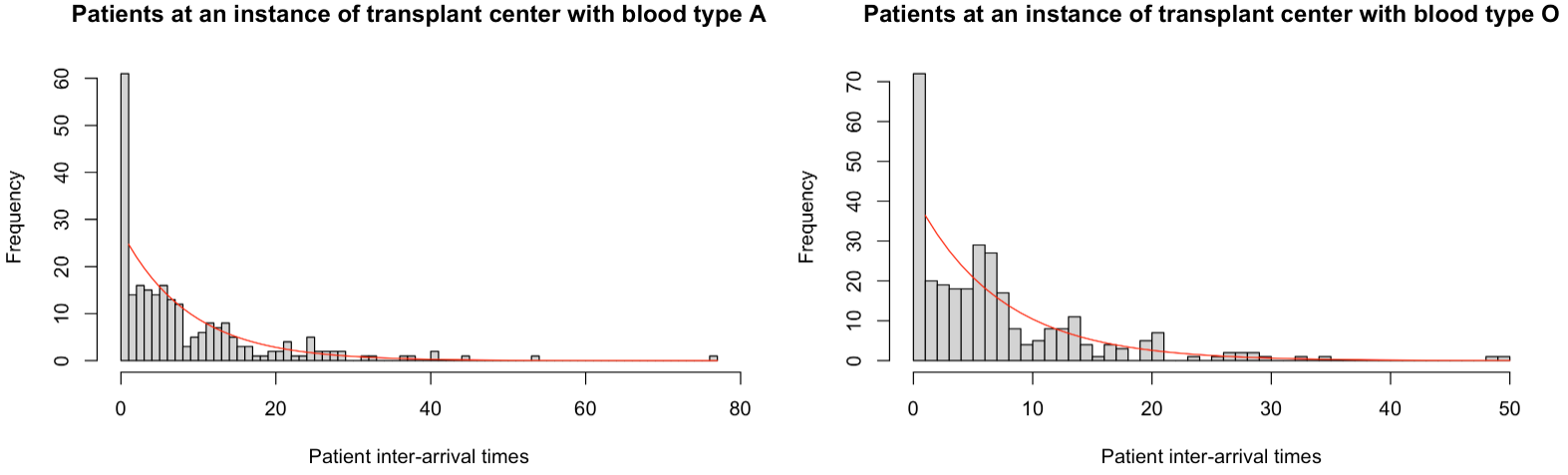}}
\vspace{-20pt}
\caption{Gap between the histogram of patient inter-arrival times (days) and the fitted exponential distribution in two example patient waitlists. 
}
\label{fig:exp test inter-arrival}
\end{figure}

\subsection{Non-Exponentiality of Patient Inter-Arrival Times}\label{sec:non-exp of IAT}
We now discuss that the patient inter-arrival time distributions of transplant queues are not exponential. 
For each instance of patient waitlist, we use moment matching to fit the inter-arrival times with an exponential distribution and a mixture of two exponential distributions, respectively. 
We performed this analysis to the 50 largest waitlists and observe that, for blood type A ($resp.$ type O), when fitting patient inter-arrival times using an exponential distribution, only one waitlist ($resp.$ no waitlist) has a $p$-value greater than 0.05 in a Kolmogorov–Smirnov test, leading us to reject the hypothesis that inter-arrival times follow an exponential distribution. 
Across all waitlists, the average Kolmogorov–Smirnov distance between the fitted exponential distribution and the empirical inter-arrival time distribution is 0.14. 
When using a mixture of two exponential distributions to fit the patient inter-arrival times, 
across all waitlists, the average Kolmogorov–Smirnov distance between the fitted mixture of exponential distributions and the empirical inter-arrival time distribution is reduced to 0.07, which is 50\% less than that obtained using exponential distributions.
In addition, more than 40\% of waitlists have a $p$-value greater than $0.05$ in a Kolmogorov–Smirnov test.
This suggests that a mixture of exponential distributions can more accurately model patient inter-arrival times.
As an illustration, consider two patient waitlist instances (Figure \ref{fig:exp test inter-arrival}) with blood types A and O, respectively, and annual waitlist additions of 41 and 50, respectively. The Kolmogorov–Smirnov test indicates a small $p$-value ($p<10^{-4}$) when fitting exponential distributions. 


\subsection{Error of Markovian Simplifcations}\label{sec:error of MM1M}
We demonstrate that the performance evaluation errors are significant for transplant waitlists if we approximate the original GI/MI/1+GI$_S$ queueing model with a MI/MI/1+MI$_S$ model.
Specifically, we consider a patient waitlist in which the inter-arrival times follow a mixture of two exponential distributions Exp(10) and Exp(40) with the first component weighted by a parameter $p$. 
This suggests an annual waitlist addition of $\frac{1}{\frac{p}{10}+\frac{1-p}{40}}$, which falls within an appropriate range as the average annual waitlist addition is 23.6 across all waitlists in our data.
We assume that: $(i)$ the liver supply ($i.e.$, the service rate) matches 87.9\% of the demand; and $(ii)$ patients have an average patience time of 0.7 years, which is consistent with our data, and the patience times follow a mixture of two exponential distributions that are truncated to a maximum of 25 years to satisfy Assumption \ref{A:input distributions 1}.
We compute the average offered waiting time and the average waiting-line abandonment probability for the original system and an approximate system. 
For the later, the inter-arrival time distribution is simplified into Exp$\big(\frac{1}{\frac{p}{10}+\frac{1-p}{40}}\big)$ and the patience time distribution is simplified into Exp$(\frac{1}{0.7})$.
The relative errors are significant (Table \ref{table}).

\begin{table}[!htb]
\centering
\caption{Relative errors due to approximating the original GI/MI/1+GI$_S$ queue with a MI/MI/1+MI$_S$ queue.}
\label{tbale:exp error}
\label{table}
\scriptsize
\begin{tabular}{|l|ccccc|}
\hline
Weight parameter $p$                               & $p=0.1$     & $p=0.3$     & $p=0.5$     & $p=0.7$     & $p=0.9$     \\ \hline
Relative error in average offered waiting time                 & 17.99\% & 21.30\% & 17.30\% & 13.70\% & 8.06\% \\
Relative error in average waiting-line abandonment probability & 8.66\% & 1.37\% & 5.34\% & 10.30\% & 17.78\% \\ \hline
\end{tabular}
\end{table}

\section{Tables of Inter-Arrival Time Distributions}\label{Append:Arrivial distribution}

Table \ref{table:arrivial distribution} summarizes the common distributions that satisfy (i) of Assumption \ref{A:input distributions 1}.

\begin{table}[!htb]
\centering
\scriptsize
\caption{Inter-arrival time distributions that satisfy (i) of Assumption 1.}
\label{table:arrivial distribution}
\begin{tabular}{lll}
\hline
\multicolumn{1}{l|}{\begin{tabular}[c]{@{}l@{}}\textbf{Inter-arrival time distribution}\end{tabular}} & \multicolumn{1}{l|}{\textbf{Notation}} & \begin{tabular}[c]{@{}l@{}}\textbf{Conditions to satisfy (i) of Assumption 1}\end{tabular} \\ \hline
\multicolumn{3}{l}{\textit{Supported on a bounded interval}} \\ \hline
\multicolumn{1}{l|}{Beta distribution} & \multicolumn{1}{l|}{Beta$(\alpha,\beta)$} & $\alpha,\beta$ are integers, or $\alpha,\beta\geqslant2.$ \\ \hline 
\multicolumn{1}{l|}{Uniform distribution} & \multicolumn{1}{l|}{-} & No extra conditions needed. \\ \hline 
\multicolumn{1}{l|}{Log-normal distribution} & \multicolumn{1}{l|}{-} & No extra conditions needed. \\ \hline
\multicolumn{1}{l|}{Triangular distribution} & \multicolumn{1}{l|}{-} & No extra conditions needed. \\ \hline 
\multicolumn{1}{l|}{Truncated normal distribution} & \multicolumn{1}{l|}{-} &  No extra conditions needed. \\ \hline
\multicolumn{3}{l}{\textit{Supported on $[0,\infty)$}} \\ \hline 
\multicolumn{1}{l|}{Erlang distribution} & \multicolumn{1}{l|}{-} &  No extra conditions needed. \\ \hline 
\multicolumn{1}{l|}{Exponential distribution} & \multicolumn{1}{l|}{-} &  No extra conditions needed. \\ \hline 
\multicolumn{1}{l|}{A finite mixture of exponential distributions} & \multicolumn{1}{l|}{-} &  No extra conditions needed. \\ \hline 
\multicolumn{1}{l|}{Exponential-logarithmic distribution} & \multicolumn{1}{l|}{-} &  No extra conditions needed. \\ \hline 
\multicolumn{1}{l|}{Gamma distribution} & \multicolumn{1}{l|}{$\Gamma(k,\theta)$} & Shape parameter $k\geqslant 2$ or $k=1$. \\ \hline 
\multicolumn{1}{l|}{Generalized gamma distribution} & \multicolumn{1}{l|}{GG$(\alpha,d,p)$} & Shape parameter $d\geqslant 2$, or $d=1$ and $p\geqslant 1$. \\ \hline 
\multicolumn{1}{l|}{Generalized Pareto distribution} & \multicolumn{1}{l|}{GPD$(\mu,\sigma,\xi)$} & Shape parameter $\xi\geqslant 0$ or $\xi\in (-\frac{1}{2},0) $. \\ \hline 
\multicolumn{1}{l|}{Inverse-gamma distribution} & \multicolumn{1}{l|}{-} & No extra conditions needed.\\ \hline 
\multicolumn{1}{l|}{Inverse-Guassian distribution} & \multicolumn{1}{l|}{-} & No extra conditions needed.  \\ \hline  
\multicolumn{1}{l|}{Pareto distribution} &  \multicolumn{1}{l|}{-} & No extra conditions needed. \\ \hline  
\multicolumn{1}{l|}{Weibull distribution} & \multicolumn{1}{l|}{Weibull$(\lambda,k)$} & Shape parameter $k\geqslant 2$ or $k=1$. \\ \hline
\end{tabular}
\end{table}

\section{Preliminaries}
\label{append:prelim}
\subsection{Banach Space Settings}
We first construct a Banach space of distribution functions. 
Let $\mathcal{B}$ be the Borel algebra for $\mathbb{R}_+$ and $\mathbf{D}$ be the collection of probability distribution functions supported on $\mathbb{R}_+$. 
Define $\mathbf{X}:=$\rm span$(\mathbf{D})= \left\{\sum_{k=1}^{n} a_k f_k\,|\, n\in \mathbb{N}, n\geqslant 1, a_k \in \mathbb{R}, f_k \in \mathbf{D} \right\}$, which is also the linear space of distribution functions of all finite signed measures supported on $\mathbb{R}_+$. 
Let $\bar{\mathbf{X}}$ be the closure of $\mathbf{X}$ with norm $\|\cdot\|_\infty$.
Then, $(\bar{\mathbf{X}}, \|\cdot\|_\infty)$ is a Banach space (Theorem 2.5 in \citealt{li2022numerical}).

We next define operations on space $\bar{\mathbf{X}}$.
We write ``$f_k$ converging to $f$ on $\bar{\mathbf{X}}$'' as ``$f_k\rightrightarrows f$'' to emphasize the uniform convergence. 
Let $\mathcal{I}$ be the identity operator.
For an operator $\mathcal{L}$ on $\bar{\mathbf{X}}$, let
$\mathcal{L}^{-1}$ be the inverse and $\|\mathcal{L}\|_O$ be the operator norm, $i.e.$, $\|\mathcal{L}\|_O := \sup_{f\in \bar{\mathbf{X}}, \|f\|_\infty = 1} \|\mathcal{L}f\|_\infty.$

\subsection{Cited Theorems} \label{append:cited theorems}
\begin{lemma}[Theorem 21.67 in \citealt{hewitt2013real}]
\label{lemma:integral by parts}
Let $\alpha$ and $\beta$ be any two real-valued  functions on $\mathbb{R}$ of finite variation,  and let $\lambda_\alpha$ and $\lambda_\beta$ be their corresponding Lebesgue-Stieltjes measures. Then $a<b$ in $\mathbb{R}$ implies 
\begin{align*}
    \int_{[a,b]}\frac{\beta(x+)+\beta(x-)}{2}\hat{\mathrm{d}} \lambda_\alpha(x)+\int_{[a,b]}\frac{\alpha(x+)+\alpha(x-)}{2}\hat{\mathrm{d}} \lambda_\beta(x) = \alpha(b+)\beta(b+)-\alpha(a-)\beta(a-).
\end{align*}
\end{lemma}
Here $\hat d$ denotes a Lebesgue-Stieltjes measure.

The finite approximation for a fixed continuous-state Markov chain supported on $\mathbb{R}$ is analyzed in \cite{li2023new}.
They provide conditions for \emph{pointwise} convergence of finite approximate performance measures, which, 
in our context of analyzing  $\pmb{\Psi}_k(\mu)$, can be written as: 
\begin{lemma}[Corollary 1, Theorem 4, and Lemma 2 in \citealt{li2023new}]
\label{lemma:li pointwise convergence}
Consider $k\in\mathscr{K}, l\in\mathscr{L}$, and $\mu\in[\mu_{\min},\mu_{\max}]$. We have 
$
    \lim_{r\rightarrow \infty}\mathbb{E}_{\pi^{(r)}_k(\xi;\mu)}g_l(\xi,\mu) = \mathbb{E}_{\pi_k(\xi;\mu)}g_l(\xi,\mu)
$ if:
\begin{itemize}
    \item[\textnormal{(i)}] Markov chain $\pmb{\Psi}_k(\mu)$ in \eqref{eqn:original MC} is supported on $[0,\bar{y}_k]$ and $\sup_{x\in[0,\bar{y}_k]}{V_u(\tau_k(x,u;\mu))}<\infty$.
    \item[\textnormal{(ii)}] For all $\varepsilon>0$, there exists a finite split of $[0,\bar{y}_k]$ denoted by knots $0=s_1<s_2<\dots<s_{N_{\varepsilon}}=\bar{y}_k$ and intervals $E_1=[s_1,s_2),E_2=[s_2,s_3),\dots, E_{N_\varepsilon}=\{s_{N_\varepsilon}\}$ such that for all $x_1,x_2\in E_i, i=1,2,\dots,N_\varepsilon$, we have $\Upsilon_1(x_1,x_2):= V_u(\tau_k(x_2,u;\mu) - \tau_k(x_1,u;\mu))\leqslant \varepsilon$ and $ \Upsilon_2(x_1,x_2):= |\tau_k(x_2,\bar{y}_k ;\mu)-\tau_k(x_1,\bar{y}_k;\mu)|\leqslant \varepsilon$.
    \item[\textnormal{(iii)}] There exist $\eta\in\mathbb{N}$ and $\delta>0$ such that  $\mathbb{P}_{u} [\mathcal{T}_{\{0\}}\leqslant \eta]\geqslant \delta$, $\forall u\in(0,\bar{y}_k]$. Here $\mathcal{T}_{\{0\}}$ represents the first return time to state $0$ and $\mathbb{P}_{u}$ represents the probability associated with the transition kernel $\tau_k(x,u;\mu)$ and an initial state $u$.
    In other words, for the Markov chain $\pmb{\Psi}_k(\mu)$ described in \eqref{eqn:original MC}, from any state $u\in(0,\bar{y}_k]$, the probability of returning to state $0$ within $\eta$ transitions is at least $\delta$.
    \item[\textnormal{(iv)}] There exist $\{\varepsilon^{(r)}\}_{r\in\mathbb{N}}\subseteq \mathbb{R}_{+}$ and $\lim_{r\rightarrow \infty} \varepsilon^{(r)}=0$ such that $\forall r \in \mathbb{N}$, $\forall x_1,x_2\in [c_{k_i}^{(r)},c_{k_{i+1}}^{(r)})$, $i=1,2,\dots,J^{(r)}-1$, we have $ \Upsilon_1(x_1,x_2)\leqslant \varepsilon^{(r)}$ and $ \Upsilon_2(x_1,x_2)\leqslant \varepsilon^{(r)}$. 
    \item[\textnormal{(v)}] $\{c_{k_i}^{(r_1)}\}_{i=1}^{J^{(r_1)}}\subseteq \{c_{k_i}^{(r_2)}\}_{i=1}^{J^{(r_2)}}$, $\forall r_1<r_2, r_1,r_2\in \mathbb{N}$. 
    \item[\textnormal{(vi)}] $q^{(r)}_{k_{i,j}}(\mu)={\tau}_k(c^{(r)}_{k_j},c^{(r)}_{k_i};\mu)- \tau_k(c^{(r)}_{k_{j-1}},c^{(r)}_{k_i};\mu)$,  $\forall i,j=1,2,\dots,J^{(r)}$, $r\in\mathbb{N}$.
    \item[\textnormal{(vii)}]  $
        \lim_{r\rightarrow\infty} {\tau}_k^{(r)}(x,u;\mu) = {\tau}_k(x,u;\mu),\forall u,x\in[0,\bar{y}_k],$ where \begin{align}
            &{\tau}_k^{(r)}(x,u;\mu):=  \sum_{i=1}^{J^{(r)}}\sum_{j=1}^{J^{(r)}}q^{(r)}_{k_{i,j}}\mathds{1}\{x\geqslant c_{k_j}^{(r)},u\in(c^{(r)}_{k_{i-1}},c^{(r)}_{k_i}]\},\quad    x,u\in[0,\bar{y}_k],r\in\mathbb{N}.  \label{eqn:def approximate kernel}
        \end{align}
    \item[\textnormal{(viii)}] $V_\xi(g_l(\xi,\mu))<\infty.$
\end{itemize}
\end{lemma}
Here items (i)-(iii) require that the original Markov chain's transition kernel has a proper variation. 
Items (iv)-(vii) require that the finite approximate Markov chain is properly defined via truncation. 
Lastly, item (viii) requires proper performance functions. 
(In \citealt{li2023new}, our items (i)-(iii) are denoted as their Condition 1-3, respectively. 
Particularly, according to Theorem 4 and Lemma 2 in \citealt{li2023new}, their Condition 3 is ensured by our item (iii).
Our items (iv)-(vii) are denoted as their Condition 4.
Our item (viii) is required in their Corollary 1.)

The next lemma provides tractable error bounds for finite approximate performance measures. 

\begin{lemma}[Theorems 1--2 in \citealt{li2023new}]
\label{lemma:li error bound 1}
Consider a function $g(\xi):\mathbb{R}\mapsto\mathbb{R}$ such that $V(g(\xi))<\infty$.
For all $k\in\mathscr{K}$, $\mu\in[\mu_{\min},\mu_{\max}]$, and $r\in\mathbb{N}$, we have the following error bound if condition \textnormal{(i)} in Lemma \ref{lemma:li pointwise convergence} holds and the finite-state Markov chain $\pmb{\Psi}^{(r)}_k(\mu)$ described in \eqref{eqn:def approximate MC} has an absorbing communicating class.
\begin{align}
    & \Big | \mathbb{E}_{\pi_k(\xi;\mu)}g(\xi) -  \mathbb{E}_{\pi^{(r)}_k(\xi;\mu)}g(\xi) \Big | 
    \leqslant  2V(g(\xi)) \cdot \sup_{x,u\in[0,\bar{y}_k]}  \big| {\tau}_k(x,u;\mu)- \tau^{(r)}_k(x,u;\mu)\big|\cdot  \|(\mathcal{I}-\mathcal{K}^{(r)}_{k,\mu})^{-1}\|_O. \label{eqn:li error bound 1}
\end{align}
Here $\mathcal{K}_{k,\mu}^{(r)}$ is a $($finite-rank$)$ operator on $\bar{\mathbf{X}}$ defined as 
\begin{align}
    \mathcal{K}_{k,\mu}^{(r)} f(x) := \int_0^{\bar{y}_k} \tau^{(r)}_k(x,u;\mu) df(u)  =
    \sum_{i=1}^{J^{(r)}} \tau^{(r)}_k(x,c^{(r)}_{k_i};\mu) \cdot [f(c^{(r)}_{k_i}) - f(c^{(r)}_{k_{i-1}})],\quad  x\in \Omega, f\in\bar{\mathbf{X}}. \label{Def:K_1}
\end{align}
Moreover, $(\mathcal{I}-\mathcal{K}^{(r)}_{k,\mu})^{-1}$ is well-defined and $\|(\mathcal{I}-\mathcal{K}^{(r)}_{k,\mu})^{-1}\|_O= e(k,r,\mu)$ as defined in Theorem \ref{theorem:FAO optimality gap}.
\end{lemma}
Here an absorbing communicating class refers to a subset of the Markov chain's states such that (a) the states in this class are mutually accessible and (b) starting from any other state outside this class, the Markov chain will be eventually absorbed into this class \citep{bremaud2013markov}.
The performance measure error bound \eqref{eqn:li error bound 1} has three factors: 
the variation of performance function $g$, the difference of transition kernels $\tau_k$ vs. $\tau^{(r)}_k$, and the norm of the inverse.
The last factor is explained by \cite{li2023new} as the sensitivity of kernel $\tau^{(r)}_k$.
(In \citealt{li2023new}, Theorem 1 derives the last two factors and Theorem 2 derives the first factor.)

The last lemma shows that given $k$ and $\mu$, the factor $e(k,r,\mu)$ is ``eventually" bounded over $r$.
(In \citealt{li2023new}, this boundedness result is provided by Lemma 6, whose sufficient conditions are our conditions (i)-(vii), as proven in their Lemmas 7--8.) 

\begin{lemma}[Lemmas 6--8 in \citealt{li2023new}]  
\label{lemma:li bounded inverse}
For all $k\in\mathscr{K}$ and $\mu\in[\mu_{\min},\mu_{\max}]$, if conditions \textnormal{(i)}-\textnormal{(vii)} in Lemma \ref{lemma:li pointwise convergence} hold, then there exists $r'\in\mathbb{N}$ such that for all $r>r'$, $(\mathcal{I}-\mathcal{K}^{(r)}_{k,\mu})^{-1}$ is well-defined and $\big\{\|(\mathcal{I}-\mathcal{K}^{(r)}_{k,\mu})^{-1}\|_O\big\}_{r>r'}=\big\{e(k,r,\mu)\big\}_{r>r'}$ is bounded.
\end{lemma}

\subsection{Proof of Useful Inequalities}
\label{append:useful inequalities}
We prove important inequalities that will be used to establish our main theorems.
\begin{lemma}\label{lemma:b primes}
Let $b':=\frac{1}{\mu_{\min}e}$ and $b'':=\frac{\mu_{\max}+\mu_{\min}}{\mu_{\min}^2} $. Then $(i)$ $|B(x;\mu)-B(x;\mu')|\leqslant b'|\mu-\mu'|$ for all $x\in\mathbb{R}_+$, $\{\mu,\mu'\}\subseteq[\mu_{\min},\mu_{\max}]$; and 
$(ii)$ $V_x\big(B(x;\mu)-B(x;\mu')\big)\leqslant b''|\mu-\mu'|$ for all $\{\mu,\mu'\}\subseteq[\mu_{\min},\mu_{\max}]$.   
\end{lemma}
\begin{proof}{Proof:}
For all $x\in\mathbb{R_+}$ and $z\in[\mu_{\min},\mu_{\max}]$, 
we have $\left|\frac{dB(x;z)}{dz}\right|= \left|\frac{d(1-e^{-zx}) }{dz}\right|=xe^{-zx }\leqslant \frac{1}{z}e^{-1} \leqslant \frac{1}{\mu_{\min}e}=b'$. 
Thus, $|B(x;\mu)-B(x;\mu')|\leqslant b'|\mu-\mu'|$ for all $x\in\mathbb{R}_+$, $\{\mu,\mu'\}\subseteq[\mu_{\min},\mu_{\max}].$
We also have 
\begin{align*}
    V_x(B(x;\mu)-B(x;\mu')) = &\int_{\mathbb{R}_+} \left|\frac{d\big(B(x;\mu)-B(x;\mu')\big)}{dx} \right|dx = \int_{\mathbb{R}_+} \left|\frac{dB(x;\mu)}{dx}-\frac{dB(x;\mu')}{dx} \right|dx \\
    = & \int_{\mathbb{R}_+}\left|\mu e^{-\mu x}- \mu' e^{-\mu' x}\right|dx \leqslant \int_{\mathbb{R}_+}\left|\mu - \mu'\right| e^{-\mu x}+\mu'\left| e^{-\mu x}-  e^{-\mu' x}\right| dx\\
    \leqslant & \int_{\mathbb{R}_+}\left|\mu - \mu'\right| e^{-\mu x}+\mu' xe^{-\mu_{\min} x}\left|\mu - \mu'\right| dx \\
    = & |\mu-\mu'| \cdot \frac{1}{\mu}+|\mu-\mu'| \cdot \frac{\mu'}{\mu_{\min}^2}   \leqslant  |\mu - \mu'| \cdot\frac{\mu_{\max}+\mu_{\min}}{\mu_{\min}^2}
    \leqslant b''|\mu-\mu'|. \qquad \qedwhite
\end{align*}
\end{proof}
\begin{lemma}
Under Assumption \ref{A:input distributions 1}, for all $k\in\mathscr{K}$ and $\mu\in[\mu_{\min},\mu_{\max}]$, we have 
\begin{align}
    &\sup_{x,u\in[0,\bar{y}_k]}\Big|\frac{d\tau_k(x,u;\mu)}{du}\Big|\leqslant a_k', \label{ieq:dtau du}\\
    &\sup_{x,u\in[0,\bar{y}_k]}\Big|\frac{d\tau_k(x,u;\mu)}{dx}\Big|\leqslant a_k', \label{ieq:dtau dx} \\
    &\sup_{x,u\in[0,\bar{y}_k]}\Big|\frac{d^2\tau_k(x,u;\mu)}{dxdu}\Big|\leqslant a_k'', \label{ieq:d2tau dudx}\\
    &\sup_{x,u\in[0,\bar{y}_k]}|\tau^{(r)}_k(x,u;\mu)-\tau_k(x,u;\mu)|\leqslant  \frac{\bar{y}_ka'_k}{2^{r-1}}, \quad r\in\mathbb{N}, \label{eqn:approximate tau gap} \\
    & \sup_{x,u\in[0,\bar{y}_k]}  \big| {\tau}_k(x,u;\mu')- \tau_k(x,u;\mu)\big|\leqslant  b'|\mu'-\mu|,\quad \mu'\in[\mu_{\min},\mu_{\max}].\label{eqn:tau mu Lip cont}
\end{align}
\end{lemma}
\begin{proof}{Proof:}
We first rewrite $\tau_k(x,u;\mu)$. 
Recall its expression \eqref{eqn:def kernel}, which can be rewritten as $\tau_k(x,u;\mu) = 
\int_{\mathbb{R}_+} \int_{\mathbb{R}_+} A^*_k(\min\{u+s, \max\{y,u\}\}-x)  dG_k(y)dB(s;\mu)$.

\eqref{ieq:dtau du} and \eqref{ieq:dtau dx}. 
Because $|\frac{dA_k(x)}{dx}|\leqslant a_k', \forall x\in\mathbb{R}_+$ according to item (i) of Assumption \ref{A:input distributions 1}, we have that 
$|\frac{d\tau_k(x,u;\mu)}{du}|\leqslant \int_\mathbb{R_+}\int_\mathbb{R_+} |\frac{dA^*_k(\min\{u+s, \max\{y,u\}\}-x) }{du}|dG_k(y)dB(s;\mu)\leqslant a_k',\forall x,u\in[0,\bar{y}_k]$ and $|\frac{d\tau_k(x,u;\mu)}{dx}|\leqslant \int_\mathbb{R_+}\int_\mathbb{R_+} |\frac{dA^*_k(\min\{u+s, \max\{y,u\}\}-x)}{dx}|dG_k(y)dB(s;\mu)\leqslant a_k',\forall x,u\in[0,\bar{y}_k]$.

\eqref{ieq:d2tau dudx}. 
Because $|\frac{d^2A_k(x)}{dx^2}|\leqslant a''_k, \forall x\in\mathbb{R}_+$ according to item (i) of Assumption \ref{A:input distributions 1}, we have $|\frac{d^2\tau_k(x,u;\mu)}{dxdu}|\leqslant \int_\mathbb{R_+}\int_\mathbb{R_+} |\frac{d^2A^*_k(\min\{u+s, \max\{y,u\}\}-x)}{dxdu}|dG_k(y)dB(s;\mu)\leqslant a''_k,\forall x,u\in[0,\bar{y}_k]$.

\eqref{eqn:approximate tau gap}. 
By definition \eqref{eqn:def approximate kernel}, $\tau^{(r)}_k(x,u;\mu)$ is a piecewise function on $[0,\bar{y}_k]^2$ formed by a grid $\{c_{k_i}^{(r)}\}_{i=1}^{J^{(r)}}\times \{c_{k_i}^{(r)}\}_{i=1}^{J^{(r)}}$. 
(Note here $\mu$ is a fixed service rate. 
)
The grid unit length is $\frac{\bar{y}_k}{2^r}$ by definition \eqref{eqn:def approximate MC}. 
If $(x,u)$ is a grid point, $i.e.$,  $(x,u)\in\{c_{k_i}^{(r)}\}_{i=1}^{J^{(r)}}\times \{c_{k_i}^{(r)}\}_{i=1}^{J^{(r)}}$, we have $\tau_k^{(r)}(x,u;\mu)=\tau_{k}(x,u;\mu)$.
Recall that $\tau_k(x,u;\mu)$ has bounded partial derivatives as in \eqref{ieq:dtau du} and \eqref{ieq:dtau dx}.
Then we have $\sup_{x,u\in[0,\bar{y}_k]}|\tau^{(r)}_k(x,u;\mu)-\tau_k(x,u;\mu)|\leqslant 2 \cdot \frac{\bar{y}_k}{2^r}\cdot a_k' = \frac{\bar{y}_ka'_k}{2^{r-1}}$.

\eqref{eqn:tau mu Lip cont}. 
We rewrite $\tau_k(x,u;\mu)$: 
\begin{align}
    \tau_k(x,u;\mu) = 
    &\int_{\mathbb{R}_+} \int_{\mathbb{R}_+} A^*_k(\min\{u+s, \max\{y,u\}\}-x)  dG_k(y)dB(s;\mu) \nonumber \\
    = &\int_{\mathbb{R}_+} \int_{\mathbb{R}_+} \int_{\mathbb{R}_+} \mathds{1}\{[\min\{u+s, \max\{y,u\}\}-t]_+\leqslant x \} dA_k(t) dG_k(y)dB(s;\mu)\nonumber \\
    & \text{(by writing $A_k^*(z)$ into $\int_{\mathbb{R}_+} \mathds{1}\{t \geqslant z \} dA_k(t)$)}\nonumber\\
    = & \int_{\mathbb{R}_+} \int_{\mathbb{R}_+} \big[\mathds{1}\{y,u \leqslant x+t \} + \mathds{1}\{y> x+t,u \leqslant x+t \} B(x+t-u;\mu)\big] dA_k(t) dG_k(y), \nonumber \\
    & \text{(by writing $\int_{\mathbb{R}_+} \mathds{1}\{s \leqslant z \} dB(s;\mu)$ into $B(z;\mu)$)}, \label{eqn:def kernel 2}
\end{align}
for all $(x,u)\in   \mathbb{R}_+^2,\mu\in[\mu_{\min},\mu_{\max}].$
Thus, we have 
\begin{align}
    & \big| {\tau}_k(x,u;\mu')- \tau_k(x,u;\mu)\big| \nonumber \\
    = &  \Big| \int_{\mathbb{R}_+} \int_{\mathbb{R}_+} \big[\mathds{1}\{y,u \leqslant x+t \} + \mathds{1}\{y> x+t,u \leqslant x+t \} B(x+t-u;\mu')\big] dA_k(t) dG_k(y)- \nonumber \\
    & \int_{\mathbb{R}_+} \int_{\mathbb{R}_+} \big[\mathds{1}\{y,u \leqslant x+t \} + \mathds{1}\{y> x+t,u \leqslant x+t \} B(x+t-u;\mu)\big] dA_k(t) dG_k(y) \Big|\nonumber\\
    = & \Big| \int_{\mathbb{R}_+} \int_{\mathbb{R}_+} \mathds{1}\{y> x+t,u \leqslant x+t \} \big[B(x+t-u;\mu') - B(x+t-u;\mu)\big] dA_k(t) dG_k(y) \Big|\nonumber\\
    \leqslant &  \int_{\mathbb{R}_+} \int_{\mathbb{R}_+}  \big|B(x+t-u;\mu') - B(x+t-u;\mu)\big| dA_k(t) dG_k(y) \nonumber\\
    \leqslant &  \int_{\mathbb{R}_+} \int_{\mathbb{R}_+}  b'|\mu'-\mu| dA_k(t) dG_k(y) \quad \text{(due to Lemma \ref{lemma:b primes})} \nonumber\\ 
    = &   b'|\mu'-\mu|,\quad (x,u)\in   \mathbb{R}_+^2,\{\mu,\mu'\}\subseteq[\mu_{\min},\mu_{\max}]. \nonumber
\end{align}
Therefore, we have $\sup_{x,u\in[0,\bar{y}_k]}  \big| {\tau}_k(x,u;\mu')- \tau_k(x,u;\mu)\big|\leqslant  b'|\mu'-\mu|,\forall  \mu'\in[\mu_{\min},\mu_{\max}].$
\qedwhite
\end{proof}

\begin{lemma}
Under Assumptions \ref{A:input distributions 1} and \ref{A:input distributions 3}, for all $k\in\mathscr{K}$, $\{\mu,\mu'\}\subseteq[\mu_{\min},\mu_{\max}]$, and $r\in\mathbb{N}$, 
\begin{align}
    \left | \inf_{f\in \mathbf{X},\|f\|_\infty=1 }{\|(\mathcal{I}-\mathcal{K}^{(r)}_{k,\mu})f\|_\infty} - \inf_{f\in \mathbf{X},\|f\|_\infty=1 }{\|(\mathcal{I}-\mathcal{K}^{(r)}_{k,\mu'})f\|_\infty} \right | \leqslant (3b' +  2b'')|\mu'-\mu|.\label{eqn:diff image ieq7}
\end{align}
\end{lemma}
\begin{proof}{Proof:}

Consider an arbitrary fixed $f\in\mathbf{X}$ such that $\|f\|_\infty=1$. 
We have
\begin{align}
    &\Big|\|(\mathcal{I}-\mathcal{K}^{(r)}_{k,\mu})f\|_\infty-\|(\mathcal{I}-\mathcal{K}^{(r)}_{k,\mu'})f\|_\infty\Big|
    \leqslant \|(\mathcal{I}-\mathcal{K}^{(r)}_{k,\mu})f-(\mathcal{I}-\mathcal{K}^{(r)}_{k,\mu'})f\|_\infty
    =  \|(\mathcal{K}^{(r)}_{k,\mu'}-\mathcal{K}^{(r)}_{k,\mu})f\|_\infty 
    \nonumber \\
    = & \sup_{x\in[0,\bar{y}_k]} \Big| \int_0^{\bar{y}_k} \tau^{(r)}_k(x,u;\mu') df(u)- \int_0^{\bar{y}_k} \tau^{(r)}_k(x,u;\mu) df(u) \Big| \quad \text{(by definitions of $\mathcal{K}^{(r)}_{k,\mu'}$ and $\mathcal{K}^{(r)}_{k,\mu}$ in \eqref{Def:K_1})}  \nonumber \\
    = & \sup_{x\in[0,\bar{y}_k]} \Big| \int_0^{\bar{y}_k} \big[\tau^{(r)}_k(x,u;\mu')- \tau^{(r)}_k(x,u;\mu)\big] df(u) \Big|.\label{eqn:diff image}
\end{align}
To derive a useful bound for \eqref{eqn:diff image}, we consider an arbitrary fixed $x\in[0,\bar{y}_k]$.
We define 
\begin{align*}
    h(u):= \begin{cases}
        \tau^{(r)}_k(x,u;\mu')- \tau^{(r)}_k(x,u;\mu), \quad \text{if $u\in[0,\bar{y}_k]$},\\
        h(\bar{y}_k),\quad \text{if $u>\bar{y}_k$},\\
        h(0),\quad \text{if $u<0$}.
    \end{cases} 
\end{align*}
Let $\hat{\mathrm{d}} p(x)$ be the Lebesgue–Stieltjes measure associated with any function $p:\mathbb{R}\mapsto \mathbb{R}$. 
Define $h_C(u) = \frac{h(u{-})+h(u{+})}{2},\, h_R(x) = h(u+),$ and $f_C(u) = \frac{f(u-)+f(u+)}{2}$ for all $u\in[0,\bar{y}_k]$. 
Then due to integral by parts (Theorem 21.67 in \citealt{hewitt2013real}; see Lemma \ref{lemma:integral by parts} in Appendix \ref{append:cited theorems}), we have  
\begin{align}
    \int_{[0,\bar{y}_k]} {h}_C(u) \hat{\mathrm{d}}    f(u) 
    =& -  \int_{[0,\bar{y}_k]} {f}_C(u) \hat{\mathrm{d}}  h_R(u)  + {h}_C(\bar{y}_k+){f}_C(\bar{y}_k+)  -  {h}_C(0-){f}_C(0-)\nonumber \\
    =& -\int_{[0,\bar{y}_k]} {f}_C(u) \hat{\mathrm{d}}  h_R(u) +h_C(\bar{y}_k){f}_C(\bar{y}_k).\nonumber 
\end{align}
Thus, we have the following bound. 
\begin{align}
    & \Big | \int_{0}^{\bar{y}_k} {h}_C(u)  d f(u) \Big|= \Big | \int_{[0,\bar{y}_k]} {h}_C(u) \hat{\mathrm{d}}    f(u) \Big| 
    = \Big| -\int_{[0,\bar{y}_k]} {f}_C(u) \hat{\mathrm{d}}  h_R(u) +h_C(\bar{y}_k){f}_C(\bar{y}_k) \Big| \nonumber\\
    \leqslant & \Big| \int_{[0,\bar{y}_k]} {f}_C(u) \hat{\mathrm{d}}  h_R(u)\Big|  + \Big| h_C(\bar{y}_k){f}_C(\bar{y}_k) \Big| 
    \leqslant \|f_C\|_\infty\cdot V(h_R) + \|h_C\|_\infty \cdot \|f_C\|_\infty \nonumber \\
    \leqslant &  \|f\|_\infty\cdot V(h) + \|h\|_\infty \cdot \|f\|_\infty  
    =  V(h) + \|h\|_\infty.\quad \text{(due to $\|f\|_\infty=1$)}
    \label{eqn:diff image ieq1}
\end{align}
Now we consider the difference between $\int_{0}^{\bar{y}_k} {h}_C(u) d   f(u)$ and $\int_0^{\bar{y}_k} h(u)df(u).$
Define $S_1= \{u\in [0,\bar{y}_k]\,|\,h(u{-})\neq h(u{+})\}$ and $S_2=\{u\in [0,\bar{y}_k]\,|\,h(u{-}) = h(u{+})\neq h(u) \}$.
Because $V(h)<\infty$, sets $S_1$ and $S_2$ are both countable, $i.e.$, $h_C$ differs from $h$ only on a countable set.
We have 
\begin{align}
    & \Big|\int_0^{\bar{y}_k} {h}_C(u) d    f(u) -  \int_0^{\bar{y}_k} h(u)df(u) \Big| 
    = \Big|\int_0^{\bar{y}_k}  [h_C(u)-h(u)] df(u) \Big|  \nonumber \\
    = & \Big| \sum_{u\in S_1\cup S_2} [h_C(u)-h(u)]\cdot [f(u)-f(u-)] \Big|\leqslant \sum_{u\in S_1\cup S_2} |h_C(u)-h(u)|\cdot 2\|f\|_\infty \nonumber\\
    & \text{(due to $df$ is a finite signed measure, and $h_C$ only differs from $h$  on $S_1\cup S_2$)} \nonumber\\
    = & 2  \sum_{u\in S_1\cup S_2} |h_C(u)-h(u)|
    \leqslant 2 \cdot \frac{V(h)}{2}=  V(h). \quad \text{(due to $\|f\|_\infty =1$)} 
    \label{eqn:diff image ieq2}
\end{align}
The last row is because every unit of difference in $|h_C(u)-h(u)|$ yields at least two units of total variation in $V(h)$.
Combining \eqref{eqn:diff image ieq1} and \eqref{eqn:diff image ieq2}, we have
\begin{align}
    \Big|\int_0^{\bar{y}_k} h(u)df(u) \Big| \leqslant & \Big|\int_{0}^{\bar{y}_k} {h}_C(u) d    f(u) -  \int_0^{\bar{y}_k} h(u)df(u) \Big| + \Big | \int_{0}^{\bar{y}_k} {h}_C(u) d    f(u) \Big|
    \leqslant  2V(h)+\|h\|_\infty. \label{eqn:diff image ieq3}
\end{align}
Using this bound and inequality \eqref{eqn:diff image}, we obtain that 
\begin{align}
    &\Big|\|(\mathcal{I}-\mathcal{K}^{(r)}_{k,\mu})f\|_\infty-\|(\mathcal{I}-\mathcal{K}^{(r)}_{k,\mu'})f\|_\infty\Big|
    \leqslant  \sup_{x\in[0,\bar{y}_k]} \Big| \int_0^{\bar{y}_k} \big[\tau^{(r)}_k(x,u;\mu')- \tau^{(r)}_k(x,u;\mu)\big] df(u) \Big|\quad \text{(due to \eqref{eqn:diff image})} \nonumber\\
    \leqslant & \sup_{x\in[0,\bar{y}_k]} 2V_u\Big(\tau^{(r)}_k(x,u;\mu')- \tau^{(r)}_k(x,u;\mu)\Big)+ \sup_{x,u\in[0,\bar{y}_k]} \Big|\tau^{(r)}_k(x,u;\mu')- \tau^{(r)}_k(x,u;\mu)\Big| \quad \text{(due to \eqref{eqn:diff image ieq3})} \nonumber \\
    \leqslant & \sup_{x\in[0,\bar{y}_k]} 2V_u\Big(\tau_k(x,u;\mu')- \tau_k(x,u;\mu)\Big)+ \sup_{x,u\in[0,\bar{y}_k]} \Big|\tau_k(x,u;\mu')- \tau_k(x,u;\mu)\Big|\nonumber  \\
    & \text{(due to the definition of $\tau^{(r)}_k(x,u;\mu)$ in \eqref{eqn:def approximate kernel})} \nonumber \\
    \leqslant & \sup_{x\in[0,\bar{y}_k]} 2V_u\Big(\tau_k(x,u;\mu')- \tau_k(x,u;\mu)\Big)+ b'|\mu'-\mu|.  \label{eqn:diff image ieq4} \quad \text{(due \eqref{eqn:tau mu Lip cont})}.
\end{align}
Consider the first term $V_u\big(\tau_k(x,u;\mu')- \tau_k(x,u;\mu)\big)$. We have
\begin{align}
    & V_u\Big( {\tau}_k(x,u;\mu')- \tau_k(x,u;\mu)\Big) \nonumber \\
    = &  V_u\Big( \int_{\mathbb{R}_+} \int_{\mathbb{R}_+} \big[\mathds{1}\{y,u \leqslant x+t \} + \mathds{1}\{y> x+t,u \leqslant x+t \} B(x+t-u;\mu')\big] dA_k(t) dG_k(y)- \nonumber \\
    & \int_{\mathbb{R}_+} \int_{\mathbb{R}_+} \big[\mathds{1}\{y,u \leqslant x+t \} + \mathds{1}\{y> x+t,u \leqslant x+t \} B(x+t-u;\mu)\big] dA_k(t) dG_k(y) \Big)\quad\text{(due to \eqref{eqn:def kernel 2})} \nonumber \\
    = & V_u\Big( \int_{\mathbb{R}_+} \int_{\mathbb{R}_+} \mathds{1}\{y> x+t,u \leqslant x+t \} \big[B(x+t-u;\mu') - B(x+t-u;\mu)\big] dA_k(t) dG_k(y) \Big)\nonumber\\
    \leqslant & \int_{\mathbb{R}_+} \int_{\mathbb{R}_+} V_u\Big(  \mathds{1}\{y> x+t,u \leqslant x+t \} \big[B(x+t-u;\mu') - B(x+t-u;\mu)\big] \Big) dA_k(t) dG_k(y) \nonumber\\
    \leqslant & \int_{\mathbb{R}_+} \int_{\mathbb{R}_+} \sup_{u\in\mathbb{R} }\Big|B(x+t-u;\mu') - B(x+t-u;\mu) \Big| +  V_u\Big(B(x+t-u;\mu')- B(x+t-u;\mu) \Big) dA_k(t) dG_k(y) \nonumber\\
    \leqslant & \int_{\mathbb{R}_+} \int_{\mathbb{R}_+} b'|\mu'-\mu| +  b''|\mu'-\mu| dA_k(t) dG_k(y) 
    \leqslant   (b' +  b'')|\mu'-\mu|.\quad \text{(due to Lemma \ref{lemma:b primes})}\label{eqn:diff image ieq6} 
\end{align}
Combining \eqref{eqn:diff image ieq4} and \eqref{eqn:diff image ieq6}, we have 
$\Big|\|(\mathcal{I}-\mathcal{K}^{(r)}_{k,\mu})f\|_\infty-\|(\mathcal{I}-\mathcal{K}^{(r)}_{k,\mu'})f\|_\infty\Big|\leqslant (3b' +  2b'')|\mu'-\mu|.$
Note this inequality is for an arbitrary fixed $f\in\mathbf{X}$ such that $\|f\|_\infty=1$.
Therefore, it indicates \eqref{eqn:diff image ieq7}. \qedwhite
\end{proof}

\section{Proof of Theorems and Lemmas}
\label{append:proof}

\subsection{Proof of Lemma \ref{lemma:fa uniform convergence}}
Our proof for finite approximate performance measures' uniform convergence has three steps. 
First, we prove a weak pointwise convergence.
Second, we show that the finite approximate performance measures are equicontinuous. (Recall that a sequence of functions $\{f_n(x):\mathbb{R}\mapsto\mathbb{R}\}_{n\in\mathbb{N}}$ are equicontinuous if  $\forall x\in\mathbb{R}, \forall \varepsilon>0,  \exists \delta>0$ such that $|f_n(x')-f_n(x)|\leqslant \varepsilon, \forall n\in\mathbb{N}$ as long as $|x'-x|\leqslant \delta$.)
Lastly, we prove uniform convergence using pointwise convergence and equicontinuity.

\noindent \textbf{Step 1. Prove pointwise convergence.}\quad  
\begin{lemma}\label{lemma:fa pointwise convergence}
Under Assumptions \ref{A:input distributions 1}--\ref{A:input distributions 2}, we have
\begin{align}
    \lim_{r\rightarrow \infty}\mathbb{E}_{\pi^{(r)}_k(\xi;\mu)}g_l(\xi,\mu) = \mathbb{E}_{\pi_k(\xi;\mu)}g_l(\xi,\mu),\quad   k\in\mathscr{K},l\in\mathscr{L},  \mu\in[\mu_{\min},\mu_{\max}].\label{eqn:fa pointwise convergence}
\end{align}
\end{lemma}
\begin{proof}{Proof:}
By Lemma \ref{lemma:li pointwise convergence}, \eqref{eqn:fa pointwise convergence} holds if the conditions listed in Lemma \ref{lemma:li pointwise convergence} hold.
We only need to verify these conditions. Consider arbitrary $k\in\mathscr{K}$, $l\in\mathscr{L}$, and $\mu\in[\mu_{\min},\mu_{\max}].$
We have the following.

(i). 
Under Assumption \ref{A:input distributions 1}, all patients' patience time is bounded by $\bar{y}_k$. 
For each patient, the server must be available after $\bar{y}_k$ from the time of their arrival.
Thus, the offered waiting time is bounded by $\bar{y}_k$ and the Markov chain $\pmb{\Psi}_k(\mu)$ is supported on $[0,\bar{y}_k]$. 
By the definition of total variation and \eqref{ieq:dtau du}, we have $
    V_u(\tau_k(x,u;\mu))= \int_{[0,\bar{y}_k]} \Big| \frac{d\tau_k(x,u;\mu)}{du} \Big|du \leqslant a_k'\bar{y}_k<\infty
$ for all $x\in[0,\bar{y}_k].$

(ii).
We will show that for $\varepsilon>0$, we can define the split knots $\{s_i\}_{i=1}^{N_\varepsilon}$ exactly as $\{c^{(r_\varepsilon)}_{k_i}\}_{i=1}^{J^{(r_\varepsilon)}}$ with $r_\varepsilon=\min\{r\in\mathbb{N}\,|\, \frac{\bar{y}_k a_k'}{2^r} \leqslant \varepsilon, \frac{\bar{y}^2_k a_k''}{2^r} \leqslant \varepsilon \}$.
Then $\forall x_1,x_2\in E_i, i=1,2,\dots,N_\varepsilon$, we have that $|x_2-x_1|\leqslant \frac{\bar{y}_k}{2^{r_\varepsilon}}$,
\begin{align}
    \Upsilon_1(x_1,x_2)= & V_u(\tau_k(x_2,u;\mu) - \tau_k(x_1,u;\mu))= \int_{0}^{\bar{y}_k} \Big|\frac{d\big[\tau_k(x_2,u;\mu) - \tau_k(x_1,u;\mu)\big]}{du}\Big| du \nonumber \\
    = & \int_{0}^{\bar{y}_k} \Big|\frac{d\tau_k(x_2,u;\mu) }{du}-\frac{d \tau_k(x_1,u;\mu)}{du}\Big| du \leqslant \int_{0}^{\bar{y}_k} \sup_{x,u\in[0,\bar{y}_k]}\Big|\frac{d^2\tau_k(x,u;\mu) }{dxdu}\Big| \cdot |x_2-x_1| du \nonumber \\
    = &\bar{y}_k\cdot a''_k\cdot |x_2-x_1|\leqslant \frac{\bar{y}^2_k a_k''}{2^{r_\varepsilon}} \leqslant \varepsilon, \quad \text{(due to \eqref{ieq:d2tau dudx})} \nonumber 
\end{align}
and, due to \eqref{ieq:dtau dx},
\begin{align}
   \Upsilon_2(x_1,x_2) = & |\tau_k(x_2,\bar{y}_k ;\mu)-\tau_k(x_1,\bar{y}_k;\mu)| \leqslant \int_{\min\{x_1,x_2\}}^{\max\{x_1,x_2\}} \Big| \frac{d\tau_k(x,\bar{y}_k;\mu) }{dx} \Big| dx   
   \leqslant  |x_2-x_1|\cdot a'_k\leqslant  \frac{\bar{y}_k a_k'}{2^{r_\varepsilon}} \leqslant \varepsilon. \nonumber 
\end{align}

(iii). 
Assumption \ref{A:input distributions 1} indicates that there exist $\varepsilon,\delta>0$ such that  $\mathbb{P}[t_k^{n+1}-s_k^n\geqslant \varepsilon]\geqslant \delta$.  
The event $\{t_k^{n+1}-s_k^n\geqslant \varepsilon\}$ indicates 
$\{ \xi_k^{n+1} \leqslant [\xi_k^{n} -\varepsilon]_+\}$, according to the system dynamics in \eqref{eqn:original MC}.
Thus, we have $\mathbb{P}\big[\xi_k^{n+1} \leqslant [\xi_k^{n} -\varepsilon]_+|\xi_k^{n} =u \big]\geqslant \delta$ for all $u\in[0,\bar{y}_k]$.
By defining $\eta:=\min\{x\in\mathbb{N}|x\geqslant \frac{\bar{y}_k}{\varepsilon}\}$, we have $\mathbb{P}_{u} [\mathcal{T}_{\{0\}}\leqslant \eta]\geqslant \delta^{\eta}$, $\forall u\in(0,\bar{y}_k]$. 
(In other words, if the independent events $\{ \xi_k^{n+1} \leqslant [\xi_k^{n} -\varepsilon]_+\}$ happen consecutively for $\eta$ times, the Markov chain $\pmb{\Psi}_k(\mu)$ will return to state $0$ within $\eta$ transitions from any state $u\in(0,\bar{y}_k]$.)

(iv). 
Define $\varepsilon^{(r)}:=\max\{\frac{\bar{y}_k a_k'}{2^r}, \frac{\bar{y}^2_k a_k''}{2^r}\},\forall r\in\mathbb{N}$.
Then we have $\lim_{r\rightarrow \infty} \varepsilon^{(r)}=0$.
Moreover, $\forall r \in \mathbb{N}$, $\forall x_1,x_2\in [c_{k_i}^{(r)},c_{k_{i+1}}^{(r)})$, $i=1,2,\dots,J^{(r)}-1$, we have
\begin{align}
    \Upsilon_1(x_1,x_2)= & V_u(\tau_k(x_2,u;\mu) - \tau_k(x_1,u;\mu))= \int_{0}^{\bar{y}_k} \Big|\frac{d\big[\tau_k(x_2,u;\mu) - \tau_k(x_1,u;\mu)\big]}{du}\Big| du \nonumber \\
    = & \int_{0}^{\bar{y}_k} \Big|\frac{d\tau_k(x_2,u;\mu) }{du}-\frac{d \tau_k(x_1,u;\mu)}{du}\Big| du \leqslant \int_{0}^{\bar{y}_k} \sup_{x,u\in[0,\bar{y}_k]}\Big|\frac{d^2\tau_k(x,u;\mu) }{dxdu}\Big| \cdot |x_2-x_1| du \nonumber \\
    = &\bar{y}_k\cdot a''_k\cdot |x_2-x_1|\leqslant \frac{\bar{y}^2_k a_k''}{2^{r}} \leqslant \varepsilon^{(r)}, \nonumber 
\end{align}
and
\begin{align}
   \Upsilon_2(x_1,x_2) = & |\tau_k(x_2,\bar{y}_k ;\mu)-\tau_k(x_1,\bar{y}_k;\mu)| \leqslant \int_{\min\{x_1,x_2\}}^{\max\{x_1,x_2\}} \Big| \frac{d\tau_k(x,\bar{y}_k;\mu) }{dx} \Big| dx  
   \leqslant  |x_2-x_1|\cdot a'_k\leqslant  \frac{\bar{y}_k a_k'}{2^{r} } \leqslant \varepsilon^{(r)}. \nonumber 
\end{align}

(v) and (vi) hold because they are in accordance with the definitions in \eqref{eqn:def approximate MC}.

(vii). By \eqref{eqn:approximate tau gap}, $\lim_{r\rightarrow\infty}\tau^{(r)}_k(x,u;\mu)= \tau_k(x,u;\mu),\forall x,u\in[0,\bar{y}_k]$.  

(viii) holds because it is exactly item (ii) of Assumption \ref{A:input distributions 2}.
\qedwhite
\end{proof}

We also obtain a useful corollary of Lemmas  \ref{lemma:li error bound 1}--\ref{lemma:li bounded inverse}. 
Consider an arbitrary set of $k\in\mathscr{K}$, $\mu\in[\mu_{\min},\mu_{\max}]$, and $r>r^*$.
Since (a) conditions (i)-(viii) in Lemma \ref{lemma:li pointwise convergence} hold, and (b) we have proven that the finite-state Markov chain $\pmb{\Psi}^{(r)}_k(\mu)$ has an absorbing communicating class in the proof of Theorem \ref{theorem:MC compute}, the results of Lemma  \ref{lemma:li error bound 1} and Lemma \ref{lemma:li bounded inverse} will hold.
Therefore, the error factor $\|(\mathcal{I}-\mathcal{K}^{(r)}_{k,\mu})^{-1}\|_O=e(k,r,\mu)$ is well-defined and uniformly bounded over $r$. 
These can be summarized as:
\begin{corollary}
\label{corollary:e3 weak bound}
Under Assumptions \ref{A:input distributions 1}--\ref{A:input distributions 2}, we have
\begin{align}
    & \Big | \mathbb{E}_{\pi_k(\xi;\mu)}g_l(\xi,\mu) -  \mathbb{E}_{\pi^{(r)}_k(\xi;\mu)}g_l(\xi,\mu) \Big | 
    \leqslant  2V_\xi(g_l(\xi,\mu)) \cdot \sup_{x,u\in[0,\bar{y}_k]}  \big| {\tau}_k(x,u;\mu)- \tau^{(r)}_k(x,u;\mu)\big|\cdot  \|(\mathcal{I}-\mathcal{K}^{(r)}_{k,\mu})^{-1}\|_O\nonumber \\
    = & 2V_\xi(g_l(\xi,\mu)) \cdot \sup_{x,u\in[0,\bar{y}_k]}  \big| {\tau}_k(x,u;\mu)- \tau^{(r)}_k(x,u;\mu)\big|\cdot e(k,r,\mu), \quad k\in\mathscr{K},l\in\mathscr{L},\mu\in[\mu_{\min},\mu_{\max}],r>r^*, \label{eqn:li error bound 2}
\end{align}
and 
\begin{align}
    \mathcal{E}_k(\mu):=\sup_{r>r^*}e(k,r,\mu)=\sup_{r>r^*}\|(\mathcal{I}-\mathcal{K}^{(r)}_{k,\mu})^{-1}\|_O<\infty,\quad  k\in\mathscr{K},\mu\in[\mu_{\min},\mu_{\max}].\label{eqn:e3 weak bound}
\end{align}
\end{corollary}

\noindent\textbf{Step 2. Equicontinuity Proof.} 
\begin{lemma}\label{lemma:fa equicontinuity}
Under Assumptions \ref{A:input distributions 1}--\ref{A:input distributions 2}, we have that for all $k\in\mathscr{K},l\in\mathscr{L}, \mu\in[\mu_{\min},\mu_{\max}], \varepsilon>0$, there exists $\delta>0$ such that for all $r>r^*$, $\mu'\in[\mu_{\min},\mu_{\max}]$ such that $|\mu'-\mu|\leqslant \delta$, we have $\big|\mathbb{E}_{\pi^{(r)}_k(\xi;\mu')}g_l(\xi,\mu') - \mathbb{E}_{\pi^{(r)}_k(\xi;\mu)}g_l(\xi,\mu) \big| \leqslant \varepsilon.$
\end{lemma}
\begin{proof}{Proof:}
Consider a fixed set of $k\in\mathscr{K},l\in\mathscr{L}$, $\{\mu,\mu'\}\subseteq[\mu_{\min},\mu_{\max}]$, and $r>r^*$.
Recall $\pi_k(\xi;\mu)$ (resp. $\pi_k(\xi;\mu')$) is the stationary distribution of the Markov chain $\pmb{\Psi}_k(\mu)$ (resp. $\pmb{\Psi}_k(\mu')$).
$\pi^{(r)}_k(\xi;\mu)$ (resp. $\pi^{(r)}_k(\xi;\mu')$) is the stationary distribution of the approximate finite-state Markov chain $\pmb{\Psi}^{(r)}_k(\mu)$ (resp. $\pmb{\Psi}^{(r)}_k(\mu')$), which can be used to approximate $\pi_k(\xi;\mu)$ (resp. $\pi_k(\xi;\mu')$).
The performance measure error bounds are given in Lemma \ref{lemma:li error bound 1}. 
Here we apply Lemma \ref{lemma:li error bound 1} but use $\pi^{(r)}_k(\xi;\mu)$ to approximate $\pi^{(r)}_k(\xi;\mu')$ instead, $i.e.$, using an approximate finite-state Markov chain to approximate another approximate finite-state Markov chain under a different service rate.
Then according to Lemma \ref{lemma:li error bound 1}, 
\begin{align}
    \Big | \mathbb{E}_{\pi^{(r)}_k(\xi;\mu')}g_l(\xi,\mu) -  \mathbb{E}_{\pi^{(r)}_k(\xi;\mu)}g_l(\xi,\mu) \Big | 
    \leqslant &  2V_\xi(g_l(\xi,\mu)) \cdot \sup_{x,u\in[0,\bar{y}_k]}  \big| {\tau}^{(r)}_k(x,u;\mu')- \tau^{(r)}_k(x,u;\mu)\big|\cdot  \|(\mathcal{I}-\mathcal{K}^{(r)}_{k,\mu})^{-1}\|_O  \nonumber  \\
    = &  2V_\xi(g_l(\xi,\mu)) \cdot \sup_{x,u\in \{c^{(r)}_{k_i}\}_{i=1}^{J^{(r)}} }  \big| {\tau}^{(r)}_k(x,u;\mu')- \tau^{(r)}_k(x,u;\mu)\big|\cdot  \|(\mathcal{I}-\mathcal{K}^{(r)}_{k,\mu})^{-1}\|_O \nonumber \\
    = &  2V_\xi(g_l(\xi,\mu)) \cdot \sup_{x,u\in \{c^{(r)}_{k_i}\}_{i=1}^{J^{(r)}} }  \big| {\tau}_k(x,u;\mu')- \tau_k(x,u;\mu)\big|\cdot  \|(\mathcal{I}-\mathcal{K}^{(r)}_{k,\mu})^{-1}\|_O. \nonumber  
\end{align}
By the triangle inequality, we have  
\begin{align}
    & \Big | \mathbb{E}_{\pi^{(r)}_k(\xi;\mu')}g_l(\xi,\mu') -  \mathbb{E}_{\pi^{(r)}_k(\xi;\mu)}g_l(\xi,\mu) \Big |  \nonumber \\
    \leqslant & \Big | \mathbb{E}_{\pi^{(r)}_k(\xi;\mu')}g_l(\xi,\mu') -  \mathbb{E}_{\pi^{(r)}_k(\xi;\mu')}g_l(\xi,\mu) \Big | + \Big | \mathbb{E}_{\pi^{(r)}_k(\xi;\mu')}g_l(\xi,\mu) -  \mathbb{E}_{\pi^{(r)}_k(\xi;\mu)}g_l(\xi,\mu) \Big | \nonumber \\
    \leqslant & \mathbb{E}_{\pi^{(r)}_k(\xi;\mu')} \big | g_l(\xi,\mu') -  g_l(\xi,\mu) \big |  + 2V_\xi(g_l(\xi,\mu)) \cdot \sup_{x,u\in \{c^{(r)}_{k_i}\}_{i=1}^{J^{(r)}} }  \big| {\tau}_k(x,u;\mu')- \tau_k(x,u;\mu)\big|\cdot  \|(\mathcal{I}-\mathcal{K}^{(r)}_{k,\mu})^{-1}\|_O.\nonumber \\
    \leqslant & \sup_{\xi\in\mathbb{R}_+} \big | g_l(\xi,\mu') -  g_l(\xi,\mu) \big |  + 2V_\xi(g_l(\xi,\mu)) \cdot \sup_{x,u\in[0,\bar{y}_k]}  \big| {\tau}_k(x,u;\mu')- \tau_k(x,u;\mu)\big|\cdot  \|(\mathcal{I}-\mathcal{K}^{(r)}_{k,\mu})^{-1}\|_O.  \label{ieq:equi-con 1} 
\end{align}

We now consider a fixed set of $k\in\mathscr{K},l\in\mathscr{L}$, $\mu\subseteq[\mu_{\min},\mu_{\max}]$, and $\varepsilon>0$.
Recall that $\lim_{\mu'\rightarrow \mu}\sup_{\xi\in\mathbb{R}_+} \big | g_l(\xi,\mu') -  g_l(\xi,\mu) \big |=0$ according to (iii) of Assumption \ref{A:input distributions 2}.
There exists $\delta_\varepsilon'>0$ such that for all $\mu'\in[\mu_{\min},\mu_{\max}]$, 
\begin{align}
    |\mu'-\mu|\leqslant \delta_\varepsilon' \quad \Rightarrow \quad \sup_{\xi\in\mathbb{R}_+} \big | g_l(\xi,\mu') -  g_l(\xi,\mu) \big |\leqslant \frac{\varepsilon}{2}. \label{ieq:equi factor 1}
\end{align}
According to (ii) of Assumption \ref{A:input distributions 2}, we have
\begin{align}
    z_1:=2V_\xi(g_l(\xi,\mu))<\infty.  \label{ieq:equi factor 2}
\end{align}
According to Corollary \ref{corollary:e3 weak bound}, 
we have 
\begin{align}
    z_2:=\mathcal{E}_k(\mu)\geqslant \|(\mathcal{I}-\mathcal{K}^{(r)}_{k,\mu})^{-1}\|_O.  \label{ieq:equi factor 3}
\end{align}
Let $\delta''_\varepsilon:=\frac{\varepsilon}{2z_1z_2b'}$. 
Then for all $r> r^*$ and $\mu'\in[\mu_{\min},\mu_{\max}]$, we have
\begin{align}
    &|\mu'-\mu|\leqslant \delta_\varepsilon'' \quad \Rightarrow \quad 2V_\xi(g_l(\xi,\mu)) \cdot \sup_{x,u\in[0,\bar{y}_k]} \big| {\tau}_k(x,u;\mu')- \tau_k(x,u;\mu)\big|\cdot  \|(\mathcal{I}-\mathcal{K}^{(r)}_{k,\mu})^{-1}\|_O\leqslant \frac{\varepsilon}{2}. \label{ieq:equi factor 4}
\end{align}
Here the inequality on the right hand side uses  \eqref{eqn:tau mu Lip cont}, \eqref{ieq:equi factor 2}, \eqref{ieq:equi factor 3}, and $|\mu'-\mu|\leqslant \delta_\varepsilon''=\frac{\varepsilon}{2z_1z_2b'}.$
Then for all $r> r^*$ and $\mu'\in[\mu_{\min},\mu_{\max}]$ such that $|\mu'-\mu|\leqslant \delta_\varepsilon := \min\{\delta_\varepsilon',\delta_\varepsilon''\}$, we have 
\begin{align}
    & \Big | \mathbb{E}_{\pi^{(r)}_k(\xi;\mu')}g_l(\xi,\mu') -  \mathbb{E}_{\pi^{(r)}_k(\xi;\mu)}g_l(\xi,\mu) \Big |  \nonumber \\
    \leqslant & \sup_{\xi\in\mathbb{R}_+} \big | g_l(\xi,\mu') -  g_l(\xi,\mu) \big |  + 2V_\xi(g_l(\xi,\mu)) \cdot \sup_{x,u\in[0,\bar{y}_k]} \big| {\tau}_k(x,u;\mu')- \tau_k(x,u;\mu)\big|\cdot \|(\mathcal{I}-\mathcal{K}^{(r)}_{k,\mu})^{-1}\|_O  \nonumber \\ 
    & \text{(due to \eqref{ieq:equi-con 1})}  \nonumber \\
    \leqslant & \frac{\varepsilon}{2} + \frac{\varepsilon}{2} = \varepsilon. \quad  \text{(due to $|\mu'-\mu|\leqslant \min\{\delta_\varepsilon',\delta_\varepsilon''\}$, \eqref{ieq:equi factor 1}, and \eqref{ieq:equi factor 4})} \nonumber
\end{align}
Thus, this $\delta_\varepsilon$ depending on $\varepsilon$ is as required by Lemma \ref{lemma:fa equicontinuity}.
\qedwhite
\end{proof}

\noindent\textbf{Step 3. Uniform convergence proof.} 
\begin{proof}{Proof of Lemma \ref{lemma:fa uniform convergence}:}
We only need to show that for all $k\in\mathscr{K}$ and $l\in\mathscr{L}$, we have $\lim_{r\rightarrow \infty}\sup_{\mu\in[\mu_{\min},\mu_{\max}]}\Big|\mathbb{E}_{\pi_k(\xi;\mu)}g_l(\xi,\mu)- \mathbb{E}_{\pi^{(r)}_k(\xi;\mu)}g_l(\xi,\mu)\Big|=0$. 
Consider a fixed set of $k\in\mathscr{K}$ and $l\in\mathscr{L}$. 
We simplify the notations by writing  $\mathbb{E}_{\pi_k(\xi;\mu)}g_l(\xi,\mu)$ as $f(\mu)$
and $\mathbb{E}_{\pi^{(r)}_k(\xi;\mu)}g_l(\xi,\mu)$ as $f_r(\mu)$.
Then we only need to prove $f_r\rightrightarrows f$.

Recall the pointwise convergence in Lemma \ref{lemma:fa pointwise convergence}, which can be written as $\lim_{r\rightarrow \infty}f_r(\mu)=f(\mu)$ for all $\mu\in[\mu_{\min},\mu_{\max}].$
We define the open interval $E(\mu,\delta):=(\mu-\delta,\mu+\delta).$
The equicontinuity in Lemma \ref{lemma:fa equicontinuity} can be written as: for all $\varepsilon>0, \mu\in[\mu_{\min},\mu_{\max}]$, there exists $\delta>0$ such that for all $r>r^*$ and $\mu'\in E(\mu,\delta)\cap[\mu_{\min},\mu_{\max}]$, we have $|f_r(\mu')-f_r(\mu)| \leqslant \varepsilon.$
Consider an arbitrary $\varepsilon >0$.
By the equicontinuity of $\{f_r\}_{r>r^*}$, for all $\mu\in [\mu_{\min},\mu_{\max}]$, we can find $\delta_{\mu}$ such that for all $r>r^*$ and $\mu'\in E(\mu,\delta_{\mu})\cap [\mu_{\min},\mu_{\max}]$, we have $|f_r(\mu')-f_r(\mu)| \leqslant \frac{\varepsilon}{3}.$
Note $\{ E(\mu,\delta_\mu)\}_{\mu\in [\mu_{\min},\mu_{\max}]}$ forms a cover of $[\mu_{\min},\mu_{\max}]$, $i.e.$, $[\mu_{\min},\mu_{\max}]\subseteq \cup_{\mu\in [\mu_{\min},\mu_{\max}]} E(\mu,\delta_\mu)$.
Because $E(\mu,\delta_\mu)$ is open for all $\mu\in[\mu_{\min},\mu_{\max}]$, we can find a finite cover such that $[\mu_{\min},\mu_{\max}]\subseteq \cup_{\mu\in S} E(\mu,\delta_\mu)$, $S\subseteq [\mu_{\min},\mu_{\max}]$ and $|S|<\infty$.
Because $S$ is a finite set and  $\{f_r\}_{r>r^*}$ has pointwise convergence, we can find $r_0>r^*$ such that $|f_{r_1}(\mu)-f_{r_2}(\mu)|\leqslant \frac{\varepsilon}{3}$ for all $\mu\in S$ and $r_1,r_2>r_0$. 
Consider an arbitrary set of $r_1,r_2>r_0$ and $\mu\in [\mu_{\min},\mu_{\max}]$.
We can find $\mu'\in S$ such that $\mu \in E(\mu',\delta_{\mu'})$ because $\mu\in [\mu_{\min},\mu_{\max}]\subseteq \cup_{\mu'\in S} E(\mu',\delta_{\mu'})$.
Because $\mu'\in S$ and $r_1,r_2>r_0$, we have $|f_{r_1}(\mu')-f_{r_2}(\mu')|\leqslant \frac{\varepsilon}{3}$.
Because $\mu\in E(\mu',\delta_{\mu'})$, we have $|f_{r_1}(\mu)-f_{r_1}(\mu')|\leqslant \frac{\varepsilon}{3}$ and $|f_{r_2}(\mu)-f_{r_2}(\mu')|\leqslant \frac{\varepsilon}{3}$.
Thus, $|f_{r_1}(\mu)-f_{r_2}(\mu)|\leqslant  |f_{r_1}(\mu)-f_{r_1}(\mu')| + |f_{r_1}(\mu')-f_{r_2}(\mu')| + |f_{r_2}(\mu')-f_{r_2}(\mu)|\leqslant \frac{\varepsilon}{3} + \frac{\varepsilon}{3} + \frac{\varepsilon}{3} \leqslant \varepsilon.$
Thus, $f_r\rightrightarrows f$. 
\qedwhite
\end{proof}

\subsection{Proof of Lemma \ref{lemma:pm continuity}}
\begin{proof}{Proof:}
Consider $k\in\mathscr{K}$ and $l\in\mathscr{L}$.
Lemma \ref{lemma:fa equicontinuity} indicates that $\mathbb{E}_{\pi^{(r)}_k(\xi;\mu)}g_l(\xi,\mu)$ is continuous w.r.t. $\mu\in [\mu_{\min},\mu_{\max}]$ for all $r>r^*$.
Lemma \ref{lemma:fa uniform convergence} indicates that $\mathbb{E}_{\pi^{(r)}_k(\xi;\mu)}g_l(\xi,\mu)$ is uniformly (over $\mu$) converging to $\mathbb{E}_{\pi_k(\xi;\mu)}g_l(\xi,\mu)$ as $r\rightarrow \infty$. 
Thus, $\mathbb{E}_{\pi_k(\xi;\mu)}g_l(\xi,\mu)$ as the limit function of a continuous and uniformly converging sequence is also continuous w.r.t $\mu$. \qedwhite
\end{proof}

\subsection{Proof of Theorem \ref{theorem: FAO convergence rate}}
Theorem \ref{theorem: FAO convergence rate} is an enhanced version of Theorem \ref{theorem: FAO consistency} by specifying the convergence rate w.r.t. $J^{(r)}$ and $N^{(m)}$ with the help of Assumption \ref{A:input distributions 3}.
The proof is based on the following two lemmas.
Lemma \ref{lemma:fa uniform convergence rate} below shows that the uniform convergence in Lemma \ref{lemma:fa uniform convergence} is subject to a rate of $O(\frac{1}{J^{(r)}})$.

\begin{lemma}[Finite approximation convergence rate]\label{lemma:fa uniform convergence rate}
Under Assumptions \ref{A:input distributions 1}--\ref{A:input distributions 3}, there exists $C_A>0$ such that $\sup_{\mu\in[\mu_{\min},\mu_{\max}]}\Big|\mathbb{E}_{\pi_k(\xi;\mu)}g_l(\xi,\mu)- \mathbb{E}_{\pi^{(r)}_k(\xi;\mu)}g_l(\xi,\mu)\Big|\leqslant \frac{C_A}{J^{(r)}},$ $\forall k\in\mathscr{K},\,\forall l\in\mathscr{L},\, \forall r> r^*.$ 
\end{lemma}
Lemma \ref{lemma:pm L continuity} below shows that the continuity in Lemma \ref{lemma:pm continuity} is subject to a Lipschitz constant. (We recall that $f(x)$ is $C$-Lipschitz continuous if $|f(x_2)-f(x_1)|\leqslant C |x_2-x_1|$ for all $x_1,x_2$.)
\begin{lemma}[Lipschitz continuous performance measures]\label{lemma:pm L continuity}
Under Assumptions \ref{A:input distributions 1}--\ref{A:input distributions 2}, there exists a constant $C_M>0$ such that for all $k\in\mathscr{K}, l\in\mathscr{L}, r>r^*$, the approximate performance measure $\mathbb{E}_{\pi^{(r)}_k(\xi;\mu)}g_l(\xi,\mu)$ as a function of $\mu\in[\mu_{\min},\mu_{\max}]$ is $C_M$-Lipschitz continuous, and 
the exact performance measure $\mathbb{E}_{\pi_k(\xi;\mu)}g_l(\xi,\mu)$ as a function of $\mu\in[\mu_{\min},\mu_{\max}]$ is $C_M$-Lipschitz continuous.
\end{lemma}

We use four steps.  
In Step 1, we prove Theorem \ref{theorem: FAO convergence rate} using Lemma \ref{lemma:fa uniform convergence rate} and Lemma \ref{lemma:pm L continuity}. 
In Step 2, we prove a useful bound for factor $e(k,r,\mu)$.
In Steps 3 and 4, we respectively prove Lemma \ref{lemma:fa uniform convergence rate} and Lemma \ref{lemma:pm L continuity} using the results established in Step 2.

\noindent \textbf{Step 1. Proof of Theorem \ref{theorem: FAO convergence rate} using Lemma \ref{lemma:fa uniform convergence rate} and Lemma \ref{lemma:pm L continuity}.}\quad

\begin{proof}{Proof of Theorem \ref{theorem: FAO convergence rate}:}
We define $C_F := 4C_A$ and $C_P := {8C_M(\mu_{\max}-\mu_{\min})}$.
We will show that these two constants are as required in Theorem \ref{theorem: FAO convergence rate}.
The proof uses similar ideas to Theorem \ref{theorem: FAO consistency}.

Consider an arbitrarily small $\varepsilon>0$, an arbitrary $r>r^*$ such that $J^{(r)}\geqslant \frac{C_F}{\varepsilon}$, and an arbitrary $m$ such that $N^{(m)}\geqslant \frac{C_P}{\varepsilon}$.
By Lemma \ref{lemma:fa uniform convergence rate}, for all $k\in\mathscr{K}, l\in\mathscr{L}$, 
\begin{align}
    \sup_{\mu\in[\mu_{\min},\mu_{\max}]}\Big|\mathbb{E}_{\pi_k(\xi;\mu)}g_l(\xi,\mu)-\mathbb{E}_{\pi^{(r)}_k(\xi;\mu)}g_l(\xi,\mu)\Big| \leqslant \frac{C_A}{J^{(r)}} \leqslant \frac{\varepsilon C_A}{C_F} = \frac{\varepsilon}{4}.
    \label{eqn:uniform convergence prime}
\end{align}
By Lemma \ref{lemma:pm L continuity}, for all $k\in\mathscr{K},l\in\mathscr{L}$, $\{\mu,\mu'\}\subseteq [\mu_{\min},\mu_{\max}]$, we have
\begin{align}
    &|\mu'-\mu|\leqslant \frac{\mu_{\max}-\mu_{\min}}{2^m} \quad   \nonumber \\ 
    \Rightarrow \quad &\Big|\mathbb{E}_{\pi_k(\xi;\mu')}g_l(\xi,\mu')-\mathbb{E}_{\pi_k(\xi;\mu)}g_l(\xi,\mu)\Big| \leqslant \frac{\mu_{\max}-\mu_{\min}}{2^m} C_M = \frac{\mu_{\max}-\mu_{\min}}{N^{(m)}-1} C_M\leqslant  \frac{\mu_{\max}-\mu_{\min}}{\frac{1}{2}N^{(m)}}C_M  \nonumber\\
    & \qquad \qquad \qquad \qquad \qquad \qquad \qquad \,\,\,\;  \leqslant  \frac{\varepsilon C_M(\mu_{\max}-\mu_{\min})}{\frac{1}{2}C_P}  =  \frac{\varepsilon}{4}.
    \label{eqn:uniform continuity prime}
\end{align}
Then, similar to \eqref{eqn:FA PWLA  error crude}, we have the following  for all $k\in\mathscr{K}, l\in\mathscr{L}, \mu\in [\mu^{(m)}_{i-1},\mu^{(m)}_{i}], i\in \{1,2,\dots,2^m\}$.
\begin{align}
    &\Big| \phi_{k,l}^{(r,m)}(\mu) - \mathbb{E}_{\pi_k(\xi;\mu)}g_l(\xi,\mu) \Big| 
    \leqslant \max \Bigg\{\Big|\phi_{k,l}^{(r,m)}(\mu^{(m)}_{i-1}) - \mathbb{E}_{\pi_k(\xi;\mu)}g_l(\xi,\mu)\Big|, \Big|\phi_{k,l}^{(r,m)}(\mu^{(m)}_i) - \mathbb{E}_{\pi_k(\xi;\mu)}g_l(\xi,\mu)\Big|\Bigg\} \nonumber\\
    = &  \max\Bigg\{\Big|\mathbb{E}_{\pi^{(r)}_k(\xi;\mu^{(m)}_{i-1})}g_l(\xi,\mu^{(m)}_{i-1}) - \mathbb{E}_{\pi_k(\xi;\mu)}g_l(\xi,\mu)\Big|, \Big|\mathbb{E}_{\pi^{(r)}_k(\xi;\mu^{(m)}_i)}g_l(\xi,\mu^{(m)}_i) - \mathbb{E}_{\pi_k(\xi;\mu)}g_l(\xi,\mu)\Big|\Bigg\}  \nonumber\\ 
    \leqslant &  \max\Bigg\{\Big|\mathbb{E}_{\pi_k(\xi;\mu^{(m)}_{i-1})}g_l(\xi,\mu^{(m)}_{i-1}) - \mathbb{E}_{\pi_k(\xi;\mu)}g_l(\xi,\mu)\Big|+\Big|\mathbb{E}_{\pi^{(r)}_k(\xi;\mu^{(m)}_{i-1})}g_l(\xi,\mu^{(m)}_{i-1}) - \mathbb{E}_{\pi_k(\xi;\mu^{(m)}_{i-1})}g_l(\xi,\mu^{(m)}_{i-1}) \Big|, \nonumber\\ 
    &  \Big|\mathbb{E}_{\pi_k(\xi;\mu^{(m)}_i)}g_l(\xi,\mu^{(m)}_i)  - \mathbb{E}_{\pi_k(\xi;\mu)}g_l(\xi,\mu)\Big|+\Big|\mathbb{E}_{\pi^{(r)}_k(\xi;\mu^{(m)}_i)}g_l(\xi,\mu^{(m)}_i)  - \mathbb{E}_{\pi_k(\xi;\mu^{(m)}_i)}g_l(\xi,\mu^{(m)}_i)\Big|\Bigg\} \nonumber \\
    \leqslant &  \max\Bigg\{\Big|\mathbb{E}_{\pi_k(\xi;\mu^{(m)}_{i-1})}g_l(\xi,\mu^{(m)}_{i-1}) - \mathbb{E}_{\pi_k(\xi;\mu)}g_l(\xi,\mu)\Big|+\frac{\varepsilon}{4}, \Big|\mathbb{E}_{\pi_k(\xi;\mu^{(m)}_i)}g_l(\xi,\mu^{(m)}_i)  - \mathbb{E}_{\pi_k(\xi;\mu)}g_l(\xi,\mu)\Big|+\frac{\varepsilon}{4}\Bigg\} \nonumber\\ 
    & \text{(due to inequality \eqref{eqn:uniform convergence prime})} \nonumber \\
    \leqslant &  \max \{\frac{\varepsilon}{4} + \frac{\varepsilon}{4}, \frac{\varepsilon}{4} + \frac{\varepsilon}{4} \} =\frac{\varepsilon}{2} \label{eqn:FA PWLA  error crude prime} \\
    &\text{(due to inequality \eqref{eqn:uniform continuity prime}, $|\mu-\mu^{(m)}_{i-1}|\leqslant \frac{\mu_{\max}-\mu_{\min}}{2^m}$, and $|\mu-\mu^{(m)}_{i}|\leqslant \frac{\mu_{\max}-\mu_{\min}}{2^m}$).}\nonumber
\end{align}
Because $[\mu_{\min},\mu_{\max}]=\cup_{i=1,2,\dots,2^m}[\mu^{(m)}_{i-1},\mu^{(m)}_{i}]$, inequality \eqref{eqn:FA PWLA  error crude prime} can also be stated as
\begin{align}
    \sup_{\mu\in[\mu_{\min},\mu_{\max}]} \Big| \phi_{k,l}^{(r,m)}(\mu) - \mathbb{E}_{\pi_k(\xi;\mu)}g_l(\xi,\mu) \Big| \leqslant \frac{\varepsilon}{2}, \quad  k\in\mathscr{K}, l\in\mathscr{L}. \label{eqn:FA PWLA  error prime}
\end{align}
Now we use inequality \eqref{eqn:FA PWLA  error prime} to prove the $\varepsilon$-optimality of solution $(\pmb{\mu}^{(r,m,\varepsilon)},\pmb{w}^{(r,m,\varepsilon)})$, $i.e.$, $(\pmb{\mu}^{(r,m,\varepsilon)},\pmb{w}^{(r,m,\varepsilon)})$ satisfies conditions (i)-(iii) of Definition \ref{def:near optimal solution}.

(i). Let $(\pmb{\mu}^*,\pmb{w}^*)=\big(({\mu}_k^*)_{k\in\mathscr{K}},({w}^*_{k,l})_{k\in\mathscr{K},l\in\mathscr{L}}\big)$ be the optimal solution of problem \eqref{opt:objective} and $(\pmb{\mu}^{(r,m,\varepsilon)},\pmb{w}^{(r,m,\varepsilon)})=\big(({\mu}_k^{(r,m,\varepsilon)})_{k\in\mathscr{K}},({w}^{(r,m,\varepsilon)}_{k,l})_{k\in\mathscr{K},l\in\mathscr{L}}\big)$ be the  solution of problem \eqref{opt:pwl objective}.  
Then $(\pmb{\mu}^*,\pmb{w}^*)$ satisfies \eqref{opt:constraint relax}, and $\mathbb{E}_{\pi_k(\xi;\mu^*_k)}g_l(\xi,\mu^*_k) \leqslant w^*_{k,l},\, (k,l)\in\mathscr{K}\times\mathscr{L}.$
Recall \eqref{eqn:FA PWLA  error prime}. 
We further have $\phi_{k,l}^{(r,m)}(\mu^*_k)  \leqslant \mathbb{E}_{\pi_k(\xi;\mu^*_k)}g_l(\xi,\mu^*_k) +\frac{\varepsilon}{2} \leqslant w_{k,l}^*+\frac{\varepsilon}{2},\, (k,l)\in\mathscr{K}\times\mathscr{L}.$
Thus, $(\pmb{\mu}^*,\pmb{w}^*)$ also satisfies \eqref{opt:pwl constraint relax}. 
Because, by definition, $(\pmb{\mu}^*,\pmb{w}^*)$ satisfies \eqref{opt:constraint} and \eqref{opt:constraint lu bound} in problem \eqref{opt:objective}, it also satisfies \eqref{opt:pwl constraint} and \eqref{opt:pwl constraint lu bound} in problem \eqref{opt:pwl objective} as these constraints are identical.
Thus, $(\pmb{\mu}^*,\pmb{w}^*)$ is a feasible solution to problem \eqref{opt:pwl objective}.
Since $(\pmb{\mu}^{(r,m,\varepsilon)},\pmb{w}^{(r,m,\varepsilon)})$ is the optimal solution to problem \eqref{opt:pwl objective}, we have $\pmb{c}^{\text{T}} \pmb{w}^{(r,m,\varepsilon)} \leqslant \pmb{c}^{\text{T}} \pmb{w}^*$, $i.e.$, (i) holds.

(ii). Because $(\pmb{\mu}^{(r,m,\varepsilon)},\pmb{w}^{(r,m,\varepsilon)})$ is the optimal solution to problem \eqref{opt:pwl objective}, it satisfies \eqref{opt:pwl constraint relax} and we have that $\phi_{k,l}^{(r,m)}(\mu^{(r,m,\varepsilon)}_k)  \leqslant w_{k,l}^{(r,m,\varepsilon)}+\frac{\varepsilon}{2},\, (k,l)\in\mathscr{K}\times\mathscr{L}.$
Recall \eqref{eqn:FA PWLA  error prime}. 
We further have
$\mathbb{E}_{\pi_k(\xi;\mu^{(r,m,\varepsilon)}_k)}g_l(\xi,\mu^{(r,m,\varepsilon)}_k)  \leqslant \phi_{k,l}^{(r,m)}(\mu^{(r,m,\varepsilon)}_k)  +\frac{\varepsilon}{2} \leqslant w_{k,l}^{(r,m,\varepsilon)}+\frac{\varepsilon}{2}+\frac{\varepsilon}{2} = w_{k,l}^{(r,m,\varepsilon)}+\varepsilon,\, (k,l)\in\mathscr{K}\times\mathscr{L}.$

(iii) also holds because it is exactly \eqref{opt:pwl constraint}. 
\qedwhite
\end{proof}

\noindent \textbf{Step 2. Uniform bound of the error factors.}\quad
\begin{lemma}
Under Assumptions \ref{A:input distributions 1}--\ref{A:input distributions 3}, we have 
$\mathcal{\bar E}:=\sup_{r>r^*,k\in\mathscr{K},\mu\in[\mu_{\min},\mu_{\max}]}e(k,r,\mu)<\infty$. 
\label{lemma:bar mathcal E}
\end{lemma}
\begin{proof}{Proof:}
Recall the definition of $\mathcal{E}_k(\mu)$ in \eqref{eqn:e3 weak bound}. 
$\mathcal{\bar{E}} = \sup_{k\in\mathscr{K},\mu\in[\mu_{\min},\mu_{\max}]}\mathcal{E}_k(\mu)$.
Thus, we only need to prove $\sup_{\mu\in[\mu_{\min},\mu_{\max}]}\mathcal{E}_k(\mu)<\infty$ for all $k\in\mathscr{K}$.
One sufficient condition is that $\mathcal{E}_k(\mu)$ is continuous on $[\mu_{\min},\mu_{\max}]$ by the boundedness theorem.
Then we only need to prove $\mathcal{E}_k(\mu)$ is continuous.

For all $k\in\mathscr{K}$, $\mu\in[\mu_{\min},\mu_{\max}]$, and $r>r^*$, we have the following equalities. 
\begin{align}
    e(k,r,\mu)=& \|(\mathcal{I}-\mathcal{K}^{(r)}_{k,\mu})^{-1}\|_O=\sup_{f\in \bar{\mathbf{X}},f\neq 0}\frac{\|(\mathcal{I}-\mathcal{K}^{(r)}_{k,\mu})^{-1}f\|_\infty}{\|f\|_\infty}  \quad \text{(by Lemma \ref{lemma:li error bound 1} and definition of operator norms)}\nonumber\\
    =& \Big\{ \inf_{f\in \bar{\mathbf{X}},f\neq 0}\frac{\|f\|_\infty}{\|(\mathcal{I}-\mathcal{K}^{(r)}_{k,\mu})^{-1}f\|_\infty} \Big\}^{-1} = \Big\{ \inf_{f\in \bar{\mathbf{X}},f\neq 0}\frac{\|(\mathcal{I}-\mathcal{K}^{(r)}_{k,\mu})f\|_\infty}{\|f\|_\infty} \Big\}^{-1} \nonumber \\
    =& \Big\{ \inf_{f\in \bar{\mathbf{X}},\|f\|_\infty=1 }{\|(\mathcal{I}-\mathcal{K}^{(r)}_{k,\mu})f\|_\infty} \Big\}^{-1}=\Big\{ \inf_{f\in \mathbf{X},\|f\|_\infty=1 }{\|(\mathcal{I}-\mathcal{K}^{(r)}_{k,\mu})f\|_\infty} \Big\}^{-1}.\label{eqn:e3 expression}
\end{align}
The last equality is due to that $\bar{\mathbf{X}}$ is the closure of $\mathbf{X}$.
Now we use \eqref{eqn:diff image ieq7}, \eqref{eqn:e3 weak bound}, and \eqref{eqn:e3 expression} to prove the continuity of $\mathcal{E}_k(\mu)$. 
For all $\varepsilon\in (0,\mathcal{E}_k(\mu))$, define $\delta:=\frac{\varepsilon}{4(3b' +  2b'')\mathcal{E}^2_k(\mu)}.$
Then for all $\mu'\in[\mu_{\min},\mu_{\max}]$ such that $|\mu'-\mu|\leqslant \delta$, we have 
\begin{align}
    &|e^{-1}(k,r,\mu)-e^{-1}(k,r,\mu')|
    = \left | \inf_{f\in \mathbf{X},\|f\|_\infty=1 }{\|(\mathcal{I}-\mathcal{K}^{(r)}_{k,\mu})f\|_\infty} - \inf_{f\in \mathbf{X},\|f\|_\infty=1 }{\|(\mathcal{I}-\mathcal{K}^{(r)}_{k,\mu'})f\|_\infty} \right |\quad \text{(by \eqref{eqn:e3 expression})}\nonumber \\
    \leqslant & (3b' +  2b'')|\mu'-\mu|\leqslant \frac{\varepsilon}{4\mathcal{E}^2_k(\mu)} .\qquad \text{(by \eqref{eqn:diff image ieq7})} 
    \label{eqn:diff image ieq8}
\end{align}
Because $\varepsilon\in(0,\mathcal{E}_k(\mu))$, $\frac{\varepsilon}{4\mathcal{E}^2_k(\mu)} \leqslant \frac{1}{4}\mathcal{E}^{-1}_k(\mu)\leqslant \frac{1}{2}\mathcal{E}^{-1}_k(\mu)$.
Consider the upper bound of $\mathcal{E}_k(\mu')$. We have 
\begin{align}
    \mathcal{E}_k(\mu') 
    = & \max_{r>r^*}e(k,r,\mu')=\big\{\min_{r>r^*}e^{-1}(k,r,\mu')\big\}^{-1}  \quad \text{(by the definition in \eqref{eqn:e3 weak bound})}\nonumber \\
    \leqslant & \Big\{\min_{r>r^*}e^{-1}(k,r,\mu)-\frac{\varepsilon}{4\mathcal{E}^2_k(\mu)}\Big\}^{-1}
    = \Big\{\mathcal{E}^{-1}_k(\mu) -\frac{\varepsilon}{4\mathcal{E}^2_k(\mu)}\Big\}^{-1}\quad \text{(by \eqref{eqn:diff image ieq8} and the definition in \eqref{eqn:e3 weak bound})}\nonumber \\
    \leqslant & \Big\{\mathcal{E}^{-1}_k(\mu) \Big\}^{-1}+4\mathcal{E}^2_k(\mu)\cdot \frac{\varepsilon}{4\mathcal{E}^2_k(\mu)}
    \leqslant  \mathcal{E}_k(\mu)+\varepsilon.\quad \text{\big(due to $\frac{\varepsilon}{4\mathcal{E}^2_k(\mu)} \leqslant \frac{1}{2}\mathcal{E}^{-1}_k(\mu)$\big)} 
    \label{eqn:diff image ieq9}
\end{align}
Consider the lower bound of $\mathcal{E}_k(\mu')$. We have 
\begin{align}
    \mathcal{E}_k(\mu') 
    = & \max_{r>r^*}e(k,r,\mu') = \big\{\min_{r>r^*}e^{-1}(k,r,\mu')\big\}^{-1} \quad \text{(by the definition in \eqref{eqn:e3 weak bound})}\nonumber \\
    \geqslant & \Big\{\min_{r>r^*}e^{-1}(k,r,\mu)+\frac{\varepsilon}{4\mathcal{E}^2_k(\mu)}\Big\}^{-1}
    =  \Big\{\mathcal{E}^{-1}_k(\mu) +\frac{\varepsilon}{4\mathcal{E}^2_k(\mu)}\Big\}^{-1}\quad \text{(by \eqref{eqn:diff image ieq8} and the definition in \eqref{eqn:e3 weak bound})}\nonumber \\
    \geqslant & \Big\{\mathcal{E}^{-1}_k(\mu) \Big\}^{-1}-4\mathcal{E}^2_k(\mu)\cdot \frac{\varepsilon}{4\mathcal{E}^2_k(\mu)} 
    \geqslant \mathcal{E}_k(\mu)-\varepsilon.\quad \text{\big(due to $\frac{\varepsilon}{4\mathcal{E}^2_k(\mu)} \leqslant \frac{1}{2}\mathcal{E}^{-1}_k(\mu)$\big)}
    \label{eqn:diff image ieq10}
\end{align}
Combining \eqref{eqn:diff image ieq9} and \eqref{eqn:diff image ieq10}, for all $\mu'\in[\mu_{\min},\mu_{\max}]$ such that $ |\mu'-\mu|\leqslant\delta$, we have $|\mathcal{E}_k(\mu')-\mathcal{E}_k(\mu)|\leqslant \varepsilon$.
Therefore, $\mathcal{E}_k(\mu)$ is continuous w.r.t. $\mu$.
\qedwhite
\end{proof}

\noindent \textbf{Step 3. Proof of Lemma \ref{lemma:fa uniform convergence rate}.}\quad
\begin{proof}{Proof:}
Let $
    C_A := 8\bar{\mathcal{E}}\cdot\sup_{l\in\mathscr{L},\mu\in[\mu_{\min},\mu_{\max}]} V_\xi(g_l(\xi,\mu)) \cdot \sup_{k\in\mathscr{K}} \bar{y}_k a_k'.$
We will prove this constant is as required in Lemma \ref{lemma:fa uniform convergence rate}.
Recall \eqref{eqn:li error bound 2}. 
For all $k\in\mathscr{K},l\in\mathscr{L},\mu\in[\mu_{\min},\mu_{\max}], r>r^*$, we have
\begin{align}
    & \Big | \mathbb{E}_{\pi_k(\xi;\mu)}g_l(\xi,\mu) -  \mathbb{E}_{\pi^{(r)}_k(\xi;\mu)}g_l(\xi,\mu) \Big | 
    \leqslant 2V_\xi(g_l(\xi,\mu)) \cdot \sup_{x,u\in[0,\bar{y}_k]}  \big| {\tau}_k(x,u;\mu)- \tau^{(r)}_k(x,u;\mu)\big|\cdot e(k,r,\mu), \nonumber \\
    \leqslant &  2V_\xi(g_l(\xi,\mu)) \cdot \sup_{x,u\in[0,\bar{y}_k]}  \big| {\tau}_k(x,u;\mu)- \tau^{(r)}_k(x,u;\mu)\big|\cdot \bar{\mathcal{E}}\quad \text{(due to Lemma \ref{lemma:bar mathcal E})} \nonumber \\
    \leqslant &  2V_\xi(g_l(\xi,\mu)) \cdot \frac{\bar{y}_ka'_k}{2^{r-1}} \cdot \bar{\mathcal{E}}
    \leqslant 2V_\xi(g_l(\xi,\mu)) \cdot \frac{4\bar{y}_ka'_k}{J^{(r)}} \cdot \bar{\mathcal{E}} \leqslant \frac{C_A}{J^{(r)}}. \quad \text{(due to \eqref{eqn:approximate tau gap} and $J^{(r)}=2^r+1$)} \nonumber
\end{align}
Thus, $\sup_{\mu\in[\mu_{\min},\mu_{\max}]}\Big|\mathbb{E}_{\pi_k(\xi;\mu)}g_l(\xi,\mu)- \mathbb{E}_{\pi^{(r)}_k(\xi;\mu)}g_l(\xi,\mu)\Big|\leqslant \frac{C_A}{J^{(r)}}.$
\qedwhite
\end{proof}

\noindent \textbf{Step 4. Proof of Lemma \ref{lemma:pm L continuity}.}\quad
\begin{proof}{Proof:}
Define $
    C_M:= \sup_{l\in\mathscr{L}}\zeta_l+2\bar{\mathcal{E}}b'\cdot\sup_{l\in\mathscr{L},\mu\in[\mu_{\min},\mu_{\max}]} V_\xi(g_l(\xi,\mu)).$
We will prove this constant is as required in Lemma \ref{lemma:pm L continuity}.
Recall \eqref{ieq:equi-con 1}. 
For all $k\in\mathscr{K},l\in\mathscr{L},\{\mu,\mu'\}\in[\mu_{\min},\mu_{\max}], r>r^*$, 
\begin{align}
    &\Big | \mathbb{E}_{\pi^{(r)}_k(\xi;\mu')}g_l(\xi,\mu') -  \mathbb{E}_{\pi^{(r)}_k(\xi;\mu)}g_l(\xi,\mu) \Big |  \nonumber \\
    \leqslant & \sup_{\xi\in\mathbb{R}_+} \big | g_l(\xi,\mu') -  g_l(\xi,\mu) \big |  + 2V_\xi(g_l(\xi,\mu)) \cdot \sup_{x,u\in [0,\Bar{y}_k] }  \big| {\tau}_k(x,u;\mu')- \tau_k(x,u;\mu)\big|  \cdot \|(\mathcal{I}-\mathcal{K}^{(r)}_{k,\mu})^{-1}\|_O  \nonumber \\
    \leqslant & \zeta_l |\mu'-\mu| + 2V_\xi(g_l(\xi,\mu))\cdot b'|\mu'-\mu| \cdot  \bar{\mathcal{E}}\leqslant C_M\cdot |\mu'-\mu|\nonumber. \quad \text{(due to  Assumption \ref{A:input distributions 3}, \eqref{eqn:tau mu Lip cont}, and Lemma \ref{lemma:bar mathcal E})} 
\end{align}
Thus, for all $k\in\mathscr{K}, l\in\mathscr{L}, r>r^*$, the approximate performance measure $\mathbb{E}_{\pi^{(r)}_k(\xi;\mu)}g_l(\xi,\mu)$ as a function of $\mu$ is $C_M$-Lipschitz continuous on $[\mu_{\min},\mu_{\max}]$. 
Moreover, because $\mathbb{E}_{\pi^{(r)}_k(\xi;\mu)}g_l(\xi,\mu)$ uniformly converges to the exact performance measure $\mathbb{E}_{\pi_k(\xi;\mu)}g_l(\xi,\mu)$ as $r\rightarrow \infty$, $\mathbb{E}_{\pi_k(\xi;\mu)}g_l(\xi,\mu)$ as a function of $\mu$ is also $C_M$-Lipschitz continuous on $[\mu_{\min},\mu_{\max}]$.
\qedwhite
\end{proof}

\subsection{Proof of Theorem \ref{theorem:FAO optimality gap}}
\begin{proof}{Proof:}
For all $k\in\mathscr{K},l\in\mathscr{L},\mu\in[\mu_{\min},\mu_{\max}], r>r^*$, we have
\begin{align}
    & \Big | \mathbb{E}_{\pi_k(\xi;\mu)}g_l(\xi,\mu) -  \mathbb{E}_{\pi^{(r)}_k(\xi;\mu)}g_l(\xi,\mu) \Big | 
    \leqslant  2V_\xi(g_l(\xi,\mu)) \cdot \sup_{x,u\in[0,\bar{y}_k]}  \big| {\tau}_k(x,u;\mu)- \tau^{(r)}_k(x,u;\mu)\big|\cdot e(k,r,\mu) \quad \text{(by \eqref{eqn:li error bound 2})} \nonumber \\
    \leqslant &  2V_\xi(g_l(\xi,\mu)) \cdot \frac{\bar{y}_ka'_k}{2^{r-1}} \cdot e(k,r,\mu).\quad \text{(due to \eqref{eqn:approximate tau gap})}\hspace{232pt}\qedwhite\nonumber
\end{align}

\end{proof}
\end{APPENDICES}

\end{document}